\documentclass[reqno]{article}

\usepackage[mathscr]{eucal}
\usepackage[all]{xy}
\usepackage{mathrsfs}
\usepackage{xypic}
\usepackage{amsfonts}
\usepackage{amsmath}
\usepackage{amsthm}
\usepackage{amssymb}
\usepackage{latexsym}
\usepackage{eufrak}

\theoremstyle{plain}
\newtheorem{thm}[equation]{Theorem}
\newtheorem{prop}[equation]{Proposition}
\newtheorem{lem}[equation]{Lemma}
\newtheorem{cor}[equation]{Corollary}
\newtheorem{claim}[equation]{Claim}

\theoremstyle{definition}
\newtheorem{defn}[equation]{Definition}

\theoremstyle{remark}
\newtheorem{examp}[equation]{Example}
\newtheorem{rem}[equation]{Remark}

\makeatletter
\renewcommand{\subsection}{\@startsection{subsection}{2}{0pt}{-3ex
plus -1ex minus -0.2ex}{-2mm plus -0pt minus
-2pt}{\normalfont\bfseries}} \makeatother

\numberwithin{equation}{subsection}

\newcommand{\erem}{\hfill$\lozenge$\end{rem}}

\newcommand{\bimod}[1]{#1\text{-}{\sf{bimod}}}
\newcommand{\alg}[1]{#1\text{-}{\sf{alg}}}

\DeclareMathOperator{\Ext}{\mathrm{Ext}}

\DeclareMathOperator{\modu}{\mathrm{mod}\,}
\DeclareMathOperator{\limp}{\underset{n\to\infty}{\textsl{lim}\,\textsl{proj}\,}}
\DeclareMathOperator{\Lie}{\mathrm{Lie}}
\DeclareMathOperator{\Tr}{\mathrm{Tr}}
\DeclareMathOperator{\Rep}{\mathrm{Rep}}
\DeclareMathOperator{\Ad}{\mathrm{Ad}}
\def\map{\longrightarrow}
\newcommand{\dis}{\displaystyle}

\newcommand{\beq}{\begin{equation}\label}
\newcommand{\eeq}{\end{equation}}

\DeclareMathOperator{\Spec}{\mathrm{Spec}}

\newcommand{\iso}{{\;\stackrel{_\sim}{\to}\;}}
\newcommand{\cd}{\!\cdot\!}

\DeclareMathOperator{\GL}{\mathrm{GL}}

\newcommand{\g}[1]{\mathfrak{#1}}
\newcommand{\scr}[1]{\mathscr{#1}}

\DeclareMathOperator{\End}{\mathrm{End}}

\newcommand{\op}{\operatorname}
\DeclareMathOperator{\Hom}{\mathrm{Hom}}

\DeclareMathOperator{\Der}{\mathrm{Der}}
\DeclareMathOperator{\dder}{\mathbb{D}\mathbf{er}}

\newcommand{\exact}{_{\op{exact}}}
\newcommand{\closed}{_{\op{closed}}}

\newcommand{\D}{T^*}

\newcommand{\id}{\mathrm{id}}
\newcommand{\Id}{\mathrm{Id}}
\newcommand{\super}{{\text{super}}}
\newcommand{\Mat}{\mathrm{Mat}}
\newcommand{\bbr}{{\mathbf{e}}}

\newcommand{\ms}{\mapsto}

\newcommand{\y}{\diamond}

\newcommand{\bpsi}{{\boldsymbol{\psi}}}
\newcommand{\G}{\Gamma}
\newcommand{\eex}{\hfill$\lozenge$\end{examp}}

\newcommand{\eps}{\varepsilon}

\newcommand{\bi}{\imath}

\newcommand{\bbi}{{\mathbf{i}}}
\newcommand{\II}{{\mathbf{I}}}

\newcommand{\BB}{{\mathsf{B}}}
\renewcommand{\lll}{{[\![}}
\newcommand{\rrr}{{]\!]}}
\renewcommand{\b}{\bar}
\newcommand{\bl}{{\scr L}}

\newcommand{\rep}{{\Rep(A,V)}}
\newcommand{\opp}{{\operatorname{op}}}
\DeclareMathOperator{\Inn}{\mathrm{Inn}}

\newcommand{\fg}{{\g g}}
\newcommand{\GG}{{\g G}}
\DeclareMathOperator{\Sym}{\mathrm{Sym}}

\DeclareMathOperator{\ad}{\mathrm{ad}}

\DeclareMathOperator{\Ker}{\mathrm{Ker}}
\DeclareMathOperator{\ev}{{\mathrm{ev}}}

\def\AA{{A_\omh}}
\def\KK{K}
\def\RR{R}
\def\Ups{\wp}

\def\calT{{\mathscr{T}}}

\newcommand{\BL}{{\mathbf{L}}}
\newcommand{\pmat}{P_{\text{matrix}}}
\newcommand{\chimat}{\chi_{\text{matrix}}}

\renewcommand{\o}{\otimes}

\newcommand{\sminus}{\smallsetminus}
\newcommand{\mto}{\longmapsto}

\newcommand{\inv}{^{-1}}
\newcommand{\vi}{${\sf {(i)}}\;$}
\newcommand{\vii}{${\sf {(ii)}}\;$}
\newcommand{\viii}{${\sf {(iii)}}\;$}
\newcommand{\iv}{${\sf {(iv)}}\;$}

\newcommand{\sset}{\subset}

\newcommand{\into}{{}^{\,}\hookrightarrow^{\,}}
\newcommand{\too}{\,\longrightarrow\,}
\newcommand{\onto}{\twoheadrightarrow}
\newcommand{\tooo}{{\;{-\!\!\!-\!\!\!-\!\!\!-\!\!\!\longrightarrow}\;}}
\newcommand{\eqq}{{=\!\!\!=\!\!\!=\!\!\!=\!\!\!=\!\!\!=}}
\newcommand{\oper}{\operatorname}

\newcommand{\N}{\mathbb N}
\newcommand{\Z}{\mathbb Z}

\newcommand{\La}{\Lambda}

\newcommand{\HH}{{H\!H}}
\newcommand{\bm}{{\mathbf{m}}}

\newcommand{\SH}{{\mathsf{H}}}
\newcommand{\DR}{{\oper{DR}}}

\renewcommand{\th}{\theta}

\newcommand{\aar}{A^R}
\newcommand{\munc}{{\mu_{\mathrm{nc}}}}
\newcommand{\bmu}{{\boldsymbol{\mu}}}

\newcommand{\bth}{\bar{\theta}}
\newcommand{\deri}{{\Der_R(A,I)}}

\newcommand{\ncO}{\Omega}

\newcommand{\Th}{\Theta}

\newcommand{\AR}{_RA}

\newcommand{\bbm}{{\mathbf{m}}}

\DeclareMathOperator{\gr}{\mathrm{gr}}

\def\ip<#1,#2>{\left\langle#1,#2\right\rangle}
\def\sp<#1>{\left\langle#1\right\rangle}

\newcommand{\eu}{\mathsf{eu}}
\newcommand{\Eu}{\mathsf{E}}

\def\ip<#1,#2>{\left\langle#1,#2\right\rangle}

\def\npb{\noindent$\bullet\quad$\parbox[t]{115mm}}

\def\hp{\hphantom{x}}

\newcommand{\bd}{{\mathbf{d}}}
\newcommand{\la}{\lambda}
\newcommand{\om}{\omega}
\newcommand{\Om}{\Omega}

\newcommand{\omh}{{\mathbf{w}}}

\newcommand{\al}{{\alpha}}
\newcommand{\be}{\beta}
\newcommand{\en}{{\enspace}}
\newcommand{\br}{_\natural}

\newcommand{\wh}{\widehat}
\newcommand{\wt}{\widetilde}
\newcommand{\pa}{\widetilde}

\newcommand{\PP}{{\mathbf{P}}}
\newcommand{\BI}{{\mathfrak{A}}}

\newcommand{\pif}{{\mathbf{\widehat\Pi}}}

\newcommand{\bt}{{\mathbf{t}}}

\def\oo{{\mathcal O}}

\def\gl{{\mathfrak{g}\mathfrak{l}}}

\def\k{{\Bbbk}}

\def\ccirc{{{}_{\,{}^{^\circ}}}}

\begin{document}


\centerline{\huge{\textbf{\Large{Noncommutative Geometry and Quiver algebras}}}}

\bigskip

\centerline{\sc William Crawley-Boevey, Pavel Etingof,
and Victor Ginzburg}
\medskip
\begin{abstract} We develop a new framework for noncommutative 
differential geometry based on {\em double} derivations.
This leads to the notion of  moment map and of
Hamiltonian reduction in noncommutative symplectic geometry.
For any smooth associative algebra $B$,
we define its noncommutative cotangent bundle $T^*B$, which is
a basic
example of noncommutative symplectic manifold. Applying
Hamiltonian reduction to noncommutative cotangent bundles
gives an interesting class of associative algebras, $\Pi=\Pi(B),$
that includes  preprojective
algebras  associated with quivers. Our formalism of
noncommutative  Hamiltonian reduction provides the
space $\Pi/[\Pi,\Pi]$ with a Lie algebra structure,
analogous to the Poisson bracket on the zero
fiber of the moment map.
In the special case where $\Pi$ is the preprojective
algebra  associated with a quiver of non-Dynkin type,
we give a complete description of the  Gerstenhaber algebra
structure on the   Hochschild cohomology of $\Pi$
in
terms of the  Lie algebra $\Pi/[\Pi,\Pi]$.
\end{abstract}.

\bigskip

\centerline{\sf Table of Contents}
\vskip -1mm

$\hspace{20mm}$ {\footnotesize \parbox[t]{115mm}{
\hp${}_{}$\!\hp1.{ $\;\,$} {\tt Introduction}\newline
\hp2.{ $\;\,$} {\tt Calculus of double derivations}\newline
\hp3.{ $\;\,$} {\tt The derivation $\Delta$}\newline
\hp4.{ $\;\,$} {\tt Hamilton operators}\newline
\hp5.{ $\;\,$} {\tt Noncommutative cotangent bundle}\newline
\hp6.{ $\;\,$} {\tt The representation functor}\newline
\hp7.{ $\;\,$} {\tt Hamiltonian reduction in Noncommutative Geometry}\newline
\hp8.{ $\;\,$} {\tt The necklace Lie algebra}\newline
\hp9.{ $\;\,$} {\tt Hochschild cohomology of preprojective algebras}\newline
10.{ $\;\,\,$} {\tt Deformations of preprojective algebras}\newline
11.{ $\;\,\,$} {\tt Representation schemes}
}}

\section{Introduction}
\subsection{} This paper is devoted
to the general formalism of
noncommutative differential and symplectic geometry.
Noncommutative symplectic geometry was
introduced by Kontsevich \cite{Ko}
and further studied in  \cite{BLB}, \cite{Gi},
where connections with quivers were discovered.

In the present paper, we develop a {\em different}
`double version' of noncommutative  geometry.
This version has much richer structure; it allows, in particular,
to deal with noncommutative analogues of the Gerstenhaber
algebra structure on polyvector
fields. Roughly speaking,  Kontsevich's version is obtained  from the
double version 
by applying the commutator quotient construction.
This procedure is  somewhat analogous to replacing a chain complex by
its homology, so a lot of information is lost in this process.

It may also be mentioned that the  theory 
developed below
plays  an important role in
the mathematical formalism  of `open string' theory,
cf. \cite{La}.

\subsection{} Let $A$ be  an arbitrary
 associative, not necessarily commutative, algebra.
We extend  the standard constructions of differential
geometry like Lie derivative or contraction
of a differential form with
respect to a vector field, to noncommutative geometry.
The role of differential forms
is played by the so-called Karoubi-de Rham complex
$\DR^\bullet A$. Derivations $A\to A$ act
naturally on $\DR^\bullet A$ via Lie derivative and contraction
operators,
 cf. \cite{Lo}.
In our `double version' approach, the role of vector fields 
is played by  double  derivations $A\to A\otimes A$
rather than ordinary derivations $A\to A.$ One of the key technical
points of the present paper is a construction of 
Lie derivative and contraction
operators for double  derivations.

For any integer $n\geq 1$, one has a
{\em representation functor} 
$$\Rep_n:\
 \text{Associative algebras}
\map \text{Affine schemes},
$$ that assigns
to an associative algebra $A$ the scheme $\Rep_nA$ of
all $n$-dimensional $A$-modules, that is, of all
homomorphisms of $A$ into the associative algebra of
$n\times n$-matrices.
A general philosophy, due in particular to Kontsevich,
says that meaningful concepts of  noncommutative
geometry should go, under the  representation functor,
to their commutative counterparts. Thus,
elements of the Karoubi-de Rham complex $\DR^\bullet A$
go to differential forms on the scheme $\Rep_nA$, and derivations
$A\to A$ go to vector fields on  $\Rep_nA$.
The meaning of double
derivations is more subtle, it is discussed in Sect. \ref{KR}.

In Section 3 we develop a
 `double' version of noncommutative
 symplectic geometry. Thus, our notion
of symplectic 2-form involves double derivations
and is different from the
one introduced by Kontsevich \cite{Ko} and used in \cite{Gi}, \cite{BLB}.
Given an algebra $A$ equipped with a symplectic 2-form,
 one has a natural
Lie algebra structure on the vector space $A/[A,A]$, analogous
to the Poisson bracket on the space of functions on
a symplectic manifold. If $A$ is smooth,
then  $\Rep_nA$ turns out to be a symplectic manifold.
Elements of  $A/[A,A]$
go, under the  representation functor,
to regular functions on  $\Rep_nA$,
and this map is compatible with the Lie brackets.

For any algebra $A$, there is a distinguished
 double  derivation $\Delta: a\mto$
$ a\o 1-1\o a$.
Contraction with  $\Delta$ gives rise to a very interesting
new operation on the  Karoubi-de Rham complex $\DR^\bullet A$.
This operation  is
closely related to the cyclic homology of $A$, see \cite{Gi2}. Using
contraction with  $\Delta$, we also define
a map  sending closed 2-forms in $\DR^2A$
to  elements of $A$ (up to a constant summand from the ground field).
This map is referred to as  {\em noncommutative
moment map}. 

Given a symplectic 2-form
$\om\in \DR^2A$, we let $\omh$ be the image of $\om$
under the  noncommutative
moment map. We show 
that the function on $\Rep_nA$ corresponding
to the element  $\omh$ gives rise to a moment map
$\Rep_nA\to (\gl_n)^*,$ see Theorem \ref{moment_thm}.
 The algebra  $A_\omh=A/A\omh A$,
the quotient of $A$ by the two-sided
ideal generated by $\omh$, is called a
Hamiltonian reduction of $A$ at $\omh$
since the scheme $\Rep_n(A_\omh)$ may be identified
with the zero fiber of the  moment map.

In general, the symplectic 2-form on $A$ does not
give a symplectic form on the Hamiltonian reduction.
The  resulting structure on $A_\omh$ is weaker;
it is conveniently expressed in terms of the notion
of {\em Hamilton operator}, see Sect. 4.3, borrowed
from the works of Gelfand and Dorfman \cite{GD}.
One of the main results of our paper  says that the Lie bracket on
$A/[A,A]$ descends to $A_\omh/[A_\omh,A_\omh]$,
see Proposition \ref{br_cor} and
Theorem \ref{br_thm}.
It should be emphasized that  the `double version'
approach  seems to be absolutely indispensable
for the above constructions and results.

For any smooth associative algebra $B$,
we define its noncommutative cotangent bundle.
This is an associative algebra $T^*B$ that
comes equipped with a canonical
symplectic 2-form $\om\in \DR^2(T^*B)$. 
Each element of $T^*B$ gives rise
to a regular function on $T^*(\Rep_nB)$, the cotangent
bundle on the scheme $\Rep_nB$. Furthermore,
the 2-form $\om$ goes, under the  representation functor,
to the standard symplectic 2-form on $T^*(\Rep_nB)$.

Applying
Hamiltonian reduction to noncommutative cotangent bundles
gives an interesting class of associative algebras, $\Pi=\Pi(B),$
that includes  preprojective
algebras  associated with quivers. This special case is considered
in more detail in section 8.
By our general result, the space $\Pi/[\Pi,\Pi]$
acquires a canonical Lie algebra structure.
Theorem \ref{main}  , which is the second main result of the paper,
gives, for quivers of non-Dynkin and not extended Dynkin type,
 a complete  description  of  the Gerstenhaber algebra
structure on the   Hochschild cohomology of $\Pi$
in
terms of the  Lie algebra $\Pi/[\Pi,\Pi]$.

Our last important result is a computation of
 the center of the corresponding 
Poisson algebra $\Sym(\Pi/[\Pi,\Pi])$,
see Theorem \ref{cen}. This should be viewed as an
`additive analogue' of a similar result for
Goldman's Lie algebra associated to a Riemann surface,
see \cite{E}. The proof depends heavily
on the results from \cite{CB1} about {\em quiver varieties}.

Our approach  to noncommutative geometry 
was partly motivated by
\cite{CB2}, and is  quite close to Van den Bergh's
work \cite{VB2}. Specifically, Van den Bergh considered
a noncommutative analogue of Poisson geometry
while in the present paper we mostly deal with a
noncommutative analogue of symplectic geometry.
Although  symplectic geometry is a special
case of  Poisson geometry, the two papers
have almost no overlap and (apart from ideological
motivation originating from quiver theory) are quite 
independent. 

For more details about the connection
between Van den Bergh's work and the present paper
the reader is referred to the Appendix to \cite{VB2}.

\subsection{Acknowledgements.} We are very grateful to Michel Van den Bergh for
useful discussions and for a careful reading of the manuscript.
 The second and third author are
partially supported by the NSF grants DMS-9988796
and  DMS-0303465, respectively, and also by
the CRDF grant RM1-2545-MO-03.

\section{Calculus of double derivations} 
\subsection{Notation.}\label{notation}
Throughout, we fix a
field $\k$ of characteristic 0 and write $\otimes=\otimes_\k,
\,\Hom=\Hom_\k,$ etc.

Given an associative $\k$-algebra
$A$, let $A^\opp$ denote
  the opposite algebra, and
write $A^e=A\otimes A^\opp$. There is a canonical
isomorphism $(A^e)^\opp\cong A^e$. Thus,
an $A$-bimodule is the same thing as a left
$A^e$-module, and also the same thing as a right $A^e$-module.

The space $A\otimes A$  has two
{\em commuting} $A$-bimodule structures,
called  the {\em outer}, resp., {\em inner}, bimodule structure.
These two   bimodule structures are given by
$$ b(a'\otimes a'')c:=(ba')\otimes (a''c),\quad\text{resp.},\quad
b(a'\otimes a'')c:=(a'c)\otimes(ba''),\,\en
a',a'',b,c\in A.
$$

Fix  a unital associative $\k$-algebra $R$.
Throughout, by an $R$-algebra we mean an
associative unital $\k$-algebra
equipped with  a
unit preserving 
$\k$-algebra imbedding $R\to A$.
A morphism of $R$-algebras is meant
to be an algebra homomorphism compatible with the
identity map $R\to R$; in particular, any 
$R$-algebra morphism
is unit preserving.

Let  $A$ be an $R$-algebra and
$A\otimes_R A\stackrel{\bm}\to A, \, a'\otimes a''\mapsto a'a'',$
the multiplication map. Let $\ncO^1\AR:=\Ker(\bm)$
 be  the $A$-bimodule of noncommutative
relative  1-forms on $A$
(with respect to the subalgebra $R$), see \cite[Sect.~2]{CQ}.
If $A$ is finitely generated as an $R$-algebra
then  $\ncO^1\AR$ is finitely generated as a
 left $A^e$-module.
\begin{defn}\label{smooth}
An $R$-algebra $A$ is called {\em smooth over} $R$
if it is finitely generated as an $R$-algebra and
 $\ncO^1\AR$ is projective as a left $A^e$-module.
\end{defn}

Path algebras associated to quivers,
to be studied later in this paper,
are typical examples of smooth algebras.

Given an $A$-bimodule $M$,
write  $\Der_R(A,M)$ for the space of $R$-linear
derivations, that is,   derivations $\th:A\to M$
 such that $\th(R)=0$. 
There is a canonical `universal'
derivation $d:A\to \ncO^1\AR,\, a\mapsto da:=a\o 1-1\o a,$
such that, for any $A$-bimodule $M$, we have a bijection
\beq{der_def}
\Der_R(A,M)\iso\Hom_{A^e}(\ncO^1\AR, M),\en
\th\mapsto i_\th\enspace\text{where}\enspace
i_\th(u\,dv):=u\cdot\th(v).
\eeq

In the `absolute' case $R=\k$, we will use unadorned notation
$\Der A=\Der_\k A, \Om^1A=\Om^1\AR$, etc. We 
 say that $A$ is {\em smooth} if $A$ is smooth over $\k$.

\subsection{Separability elements.}\label{Rsep}
In this subsection, as well as in  most sections  of the
paper, we assume  $R$ to be 
a finite dimensional semisimple $\k$-algebra, e.g., a direct
sum of matrix algebras.

For such an algebra,
there exists a canonical
{\em symmetric separability element}
$\bbr=\sum_i\,e_i\otimes e^i\in R\otimes R$
such that the following holds, cf. \cite[Prop. 4.2]{CQ}:
\begin{equation}\label{sep_def}
{\mathsf {(i)}}\en r\cd\bbr=\bbr\cd r,\en
\forall r\in R;\en\quad{\mathsf {(ii)}}\en\sum\nolimits_i\,e_i\cd e^i=1;\en\quad
{\mathsf {(iii)}}\en\bbr=\bbr{^\opp},
\eeq
where in equation (i) we write
$r\cdot \bbr:= \sum_i\,(r\cdot e_i)\otimes e^i,$
resp.,
$\bbr\cdot r:=\sum_i\,e_i\otimes (e^i\cdot r),$
and in  equation (iii)
we let $(x\otimes y)^\opp:= y\otimes x$
denote the flip involution on $R\otimes R$.
We will often use Sweedler's notation and 
write $\bbr=\bbr'\otimes\bbr''$ for
the sum $\sum_i\,e_i\otimes e^i.$

Given an $R$-bimodule $M$ and a (not necessarlily $R$-stable)
$\k$-vector subspace $V\sset M$, we
write $V^R=\{v\in V\mid rv=vr,\,\forall r\in R\}$ for the {\em centralizer} of $R$ in $V$,
and   put $\bbr\cdot V:=\{\bbr\cdot v\mid v\in V\},$
where $\bbr\cdot v$ stands for the  left $R^e$-action.
Also, let
$[R,M]\sset M$ be the
$\k$-linear span of the set $\{rm-mr\mid r\in R,m\in M\}$.
Thus, $[R,M]$ is a vector subspace in $M$.

Equation 
\eqref{sep_def}(i) means that
$\bbr\in (R\otimes R)^R,$
where $R\otimes R$ is viewed as an $R$-bimodule with respect
to the outer bimodule structure. The symmetry property
\eqref{sep_def}(iii) insures that the same holds for
the inner bimodule structure as well.

Properties \eqref{sep_def} yield the following result, cf. \cite[(43)]{CQ}.

\begin{prop}\label{sep} \vi  For any  $R$-bimodule $M$
we have
$$M=M^R\oplus [R,M],
\quad\op{and}\quad
\bbr\cd M=M^R\iso M/[R,M],
$$
where the isomorphism on the right is the composite
$M^R\into M\onto  M/[R,M].$

\vii
The projection  $M\to M^R$ along $[R,M]$ is given by the formula:
$$
M\onto \bbr\cd M,\quad m\mto
\sum\nolimits_i\,e_i\cdot m\cdot e^i=\bbr'\cd m\cd \bbr''=\bbr''\cd m\cd \bbr'.
\qquad\Box
$$
\end{prop}

Let $A$ be an $R$-algebra  and $M$ an $A$-bimodule.
Clearly, $[R,M]\sset [A,M]$. Applying Proposition \ref{sep}
we deduce 
\beq{MR}
M=M^R+[A,M],\quad\text{and}\quad\bbr\cd[A,M]= [A,M]^R.
\eeq

Observe further that, for $R$ as above,
any $R$-algebra $A$ that is smooth over $\k$ is also smooth
over $R$. To see this, note that the map
$r\mapsto r\cdot\bbr=\bbr\cdot r $ provides a section
of
the $R$-bimodule map $R\otimes R \to R,$
of multiplication. Tensoring
over $R$ on each side with $A$, we deduce that the $A$-bimodule
map $A\otimes A\to A\otimes_R A$ has a section. Thus,
$\Omega^1_R A$ is a direct summand of $\Omega^1A$,
and our claim follows.

\subsection{Double derivations.}\label{Double_derivations}
Fix an $R$-algebra  $A$.
We  consider  $R$-linear derivations
 $\Th: A\to A\otimes A$,
where we view $A\otimes A$ as an $A$-bimodule with respect
to the outer  bimodule structure.
Put  $\dder_RA:=\Der_R(A,A\otimes A)$.
The inner  bimodule structure
on  $A\otimes A$ gives  $\dder_RA$
 a natural 
 left  $A^e$-module structure.
Equivalently, we may view a double derivation
$\Th\in\dder_RA$ as a derivation $A\to A^e$,
where $A^e$ is viewed as an $A$-bimodule 
corresponding to left multiplication of
$A^e$ on itself. From this point of view,
the $A^e$-module structure
on  $\dder\AR$ described above
 comes from
right multiplication of
$A^e$ on itself.

Let $M$ be an $A$-bimodule  and $u\in M$.
The map $u_*: A\otimes A \to M, a'\otimes a''$
$\mapsto a'ua''$ is a morphism of $A$-bimodules.
Applying the functor $\Der_R(A,-)$ to this morphism
yields a map $\dder_RA=\Der_R(A,A\otimes A) \to
\Der_R(A,M)$. It is straightforward to verify that this
way one obtains a well-defined map
\beq{m}
(\dder_RA)\otimes_{A^e} M \too \Der_R(A,M),\quad
\Th\otimes u\mto u_*\ccirc \Th,
\eeq
where the tensor product on the left is taken with respect
to the $A$-bimodule structure on $\dder_RA$.

In the special case $M=A$, we put $\Der\AR:=\Der_R(A,A)$.
The 
multiplication map
$\bm: A\otimes A\to A,\,a'\otimes a''\mapsto a'a'',$ 
 is a morphism of $A$-bimodules
with respect to the outer bimodule structure on
$A\otimes A$ (but  {\em not} with respect to
the inner structure), hence,
induces a map $\bm_*:\dder\AR\to\Der\AR$.
Using the notation of \eqref{m}, one can write
$\bm_*: \Th\mapsto 1_*\ccirc \Th.$

\begin{prop}\label{form_smooth}  Let $R$ be a finite dimensional
semisimple
algebra.
For any  smooth $R$-algebra $A$ the map
 \eqref{m} is a bijection, $\dder\AR$ is a finitely generated
projective $A^e$-module,
and there is  a canonical short exact sequence
$$
0\map\Hom_{A^e}(\Om^1\AR,\Om^1A)\too\dder\AR\stackrel{\bm_*}\too\Der\AR\map
0.
$$
\end{prop}
\begin{proof} 
For any
$R$-algebra $A$
 and $A$-bimodule $M$, we have a canonical
map
\beq{form_smooth1}
\Hom_{A^e}(\Om^1\AR,\,A\otimes A)\otimes_{A^e}M
\to \Hom_{A^e}(\Om^1\AR,\,M).
\eeq
Using the universal property of $\Om^1\AR$ we rewrite this
map as
${(\dder_RA)\otimes_{A^e}M}\to$
$\Der_R(A,M)$.
It is easy to check that the latter map is nothing but
\eqref{m}.

Now, if $A$ is smooth over $R$, then
$\Omega^1_R A$ is a finitely generated  projective left $A^e$-module.
Thus, $\Omega^1_R A$ is a direct summand
of a finite rank free $A^e$-module.
It follows that 
$\Hom_{A^e}(\Om^1\AR,\,A\otimes A)=\dder\AR$ is also
a direct summand
of a finite rank free $A^e$-module, hence, projective.
Furthermore, we conclude that  the map
in \eqref{form_smooth1} is a bijection because
a similar map for $\Om^1\AR$ being replaced by
a finite rank free $A^e$-module is clearly
bijective. 

Finally, if $A$ is smooth then, applying   {\em exact}
functor
$\Hom_{A^e}(\Om^1\AR,-)$ to the short exact sequence
of $A$-bimodules 
$\Om^1A\into A\otimes A\stackrel{\bm}\onto A$
yields the short exact sequence  of the Proposition.
\end{proof}

\subsection{Double-derivations for a free algebra.}
Let $R=\k$, fix an integer $n\geq 1$ and set $A=\k\langle
x_1,\ldots,x_n\rangle$,
a free associative $\k$-algebra
on $n$ generators. 
An element $f\in \k\langle x_1,\ldots,x_n\rangle$
may be thought of as a polynomial
function in $n$ non-commuting variables $x_1,\ldots,x_n$.

For each $i=1,\ldots,n$,
we introduce a double derivation $\partial_i\in \dder A$
defined on generators by the formula:
$$\partial_i(x_j)=\begin{cases}
1\otimes 1&\text{if}\enspace i=j\\
0 & \text{if}\enspace i\neq j.
\end{cases}
$$
It is easy to see that $\dder A$ is a free left
$A^e$-module with basis $\{\partial_i,\,i=1,\ldots,n\}.$

Next, let  $f_1,\ldots,f_n\in \k\langle x_1,\ldots,x_n\rangle$
be an $n$-tuple of  elements of $A$.
We  write $F=(f_1,\ldots,f_n)$
and think of $F$ as
a self-map of an $n$-dimensional `noncommutative affine space'.

We may identify the space $A\otimes A$ with $A^e$, and
 define  the {\em Jacobi matrix} for the map $F$ to be
the following $A^e$-valued $n\times n$-matrix
\beq{jac}
DF=\|\partial_i(f_j)\|_{i,j=1,\ldots,n}\in \Mat_n(A^e).
\end{equation}

Now, let
$F=(f_1,\ldots,f_n)$ and $G=(g_1,\ldots,g_n)$
be two $n$-tuples  of  elements of $A$.
Let $G\ccirc F$ be the `composite' $n$-tuple  obtained by
substituting the elements $f_1,\ldots,f_n\in A$
into the arguments of the
noncommutative polynomials $g_1,\ldots,g_n.$
Equivalently put, giving an $n$-tuple $F=(f_1,\ldots,f_n)$
is the same thing as giving an algebra homomorphism
$F: A\to A$ such that $x_i\mapsto f_i, \,i=1,\ldots,n.$
In this interpretation, the composite $G\ccirc F$
corresponds to composing  algebra homomorphisms.

Let $P,Q\mapsto P\star Q$ denote multiplication
in the algebra $\Mat_n(A^e)$; note that it differs
from multiplication in $\Mat_n(A\o A)$.

The following result is proved by a straightforward computation.
\begin{prop}[Chain rule] For any two algebra homomorphisms
$F,G: A\to A,$ in $\Mat_n(A^e)$ one has 
$\dis D(G\ccirc F)= (DG)(F)\star DF.$\qed
\end{prop}
\noindent
In the right-hand side of the equation above,
 the
matrix $(DG)(F)$ is obtained by applying the homomorphism
$F\otimes F: A\otimes
A \to A\otimes A$ to each entry of
the matrix $DG\in\Mat_{n\times n}(A^e)$.

We remark that the above Proposition has no analogue for
ordinary derivations $A\to A$ instead of double derivations;
cf. also \cite{Vo} for closely related constructions.

\subsection{Karoubi-de Rham complex.}
Let $A$ be an $R$-algebra.
The tensor algebra,  $T_A^\bullet(\ncO^1\AR)
=\oplus_{n\geq 0}\,T_A^n(\ncO^1\AR)$ (tensor product over $A$), of the 
$A$-bimodule  $\ncO^1\AR$ is a DG algebra  $(\ncO^\bullet\AR,d),$
called
the  algebra of noncommutative
relative differential forms on $A$. 
We have an
isomorphism of left $A$-modules $\Om^n\AR=A\otimes_R T_R^n(A/R),$
see \cite{CQ}; usually, one writes
$a_0\,da_1\,da_2\ldots da_n\in \Om^n\AR$ for the $n$-form
corresponding to an element $a_0\o(a_1\o\ldots\o a_n)\in
A\otimes_R T_R^n(A/R)$ under this isomorphism.

If $R$ is a finite dimesional semisimple algebra, then
the DG  algebra  $(\ncO^\bullet\AR,d)$
is known to be acyclic in positive degrees.
Indeed, let
 $\pi: A\to A/R$ denote the projection.
The differential $d_n: \Om^n\AR\to\Om^{n+1}\AR $ corresponds, 
under the isomorphism $\Om^\bullet\AR=A\otimes_R T_R^\bullet(A/R),$ to
the composite 
$$ A\otimes_R T_R^n(A/R)\stackrel{\pi\otimes
\Id^{\otimes n}}\too T_R^{n+1}(A/R) \stackrel{_\sim}\to R\otimes_R
T_R^{n+1}(A/R)
\into A\otimes_R T_R^{n+1}(A/R),
$$
(we have used here that
 tensoring  over $R$ is an exact functor
since $R$ is a finite dimensional semisimple algebra).

From this definition of the differential, it is clear that
$\op{Im}(d_n)=\Ker(d_{n+1})$; thus,
we have
\beq{acyclic}
(\Om^k\AR)_{\op{closed}}=(\Om^k\AR)\exact:=d\left(\Om^{k-1}\AR\right),
\quad\forall k\geq 1.
\eeq

Next, following Karoubi, we
define 
$$\DR^\bullet\AR:=\ncO^\bullet\AR/[\ncO^\bullet\AR,\ncO^\bullet\AR]_\super,$$
where $[-,-]_\super$ denotes the $\k$-linear span of
all supercommutators.
This is the noncommutative de Rham complex of
$A$ with de Rham differential $d: \DR^j\AR\to\DR^{j+1}\AR$,
cf. also \cite{Lo} for more details (in the case $R=\k$).

For any  $A$-bimodule $M$, we have the commutator space
 $[A,M]$,  and we let
$M\br:= M/[A,M]= A\otimes_{A^e} M$
denote 
the corresponding commutator quotient.

We have a  natural projection $(\Om^\bullet\AR)\br\to\DR^\bullet\AR$
that becomes an  isomorphism:
$$ 
\DR^j\AR=(\Om^j\AR)\br, \quad\text{for}\quad j=0,1.
$$
For $j\geq 2$,
 the projection $(\Om^j\AR)\br\to\DR^j\AR$ is not 
an isomorphism, in general.

Sometimes, to distinguish between $\Om^\bullet\AR$ and
$\DR^\bullet\AR$ we will write $\lll\om\rrr$ for the image
of an element $\om\in\Om^\bullet\AR$
under the projection $\Om^\bullet\AR\onto\DR^\bullet\AR$.

The following standard result will be proved 
at the end of Sect. \ref{lie}.

\begin{lem}[Homotopy invariance]\label{homotopy}
Let $A=\bigoplus_{k\geq 0} A_k$
be a graded $R$-algebra such that $R\sset A_0$. Then, the 
algebra imbedding $A_0\into A$ induces an isomorphism
of de Rham cohomology
$H^\bullet(\DR^\bullet_RA_0,d)\iso
H^\bullet(\DR^\bullet_RA,d).$
\end{lem}

For  any vector space $V\sset A$, we write
$dV=\{dv,\,v\in V\},$ a subspace in $\Om^1_RA.$
Given an $R$-subalgebra $B\sset A$, there is a canonical
DG algebra isomorphism 
$\dis
\ncO^\bullet_BA\cong
\ncO^\bullet\AR/\ncO^\bullet\AR\cd dB\cd\ncO^\bullet\AR.
$
Similarly, given a two-sided ideal  $I\sset A$, one has a DG algebra
isomorphism, see \cite{CQ}:
\beq{A/I}
 \Om^\bullet_R(A/I)\cong\Om^\bullet\AR\big/
(\Om^\bullet\AR\cd I\cd \Om^\bullet\AR+
\Om^\bullet\AR\cd dI\cd \Om^\bullet\AR).
\eeq
\begin{rem}  In general, the induced map $R\to A/I$ is not
necessarily injective. So, according to our definition,
 $A/I$ may not be an $R$-algebra. Such a situation may arise
when we  consider Hamiltonian reduction in noncommutative geometry
(section 6.4 below). However, the only reason
for insisting that $R$ be 
injectively mapped into any $R$-algebra $A$ is a frequent use
of the notation $A/R$. So, in those rare occasions where
the map $R\to A$ fails to be injective, the symbol
$A/R$ should be understood as the quotient of $A$
by the image of $R$ in $A$.
\erem
\subsection{Contraction with a double  derivation.}
We will use lower case Greek letters $\th,\xi,\ldots$
to denote derivations $A\to A$, and upper case Greek
letters, $\Th,\Xi,\ldots$  to denote double derivations
$A\to A\o A$.

Any  derivation $\th\in \Der\AR$
 gives rise to   contraction (with $\th$) maps
$i_\th: \Om^k\AR\to \Om^{k-1}\AR,$
resp., $\DR^k\AR\to \DR^{k-1}\AR$. 
The  map $i_\th$ is
defined on 1-forms by  formula
 \eqref{der_def} and is
extended to a map
$\Om^\bullet\AR\to\Om^{\bullet-1}\AR$
as a
{\em super}-derivation. 

Now let
 $\Th\in \dder_RA$. 
For any 1-form $\al\in \Om^1\AR$,
 contraction with $\Th$
gives an $A$-bimodule map, see \eqref{der_def}:
\beq{bi_contr}
i_\Th:\
\Om^1\AR\too A\otimes A,\quad\al\mto i_\Th\al= i_\Th'\al \otimes i_\Th''\al.
\eeq
Here and below, we will systematically
use symbolic 
Sweedler's notation to write
$i_\Th\al=i'_\Th\al\otimes i''_\Th\al$ (omitting the summation symbol) for
an element in the tensor product.
Similarly, we write the map $\Th: A\to A\otimes A$ 
as $a\mapsto \Th'(a)\otimes\Th''(a).$

As usual, one may uniquely extend
the map $i_\Th$ in \eqref{bi_contr} to
higher degree differential forms
by requiring that  $i_\Th$ be a 
{\em super-derivation} of degree
$(-1)$. This way, we obtain a map
$i_\Th: \Om^\bullet\AR\to\Om^\bullet\AR\o\Om^\bullet\AR,$
which is a super-derivation of the graded algebra
$\Om^\bullet\AR$ with coefficients
in $\Om^\bullet\AR\o\Om^\bullet\AR,$ viewed
as an $\Om^\bullet\AR$-bimodule with respect to the outer 
bimodule structure.

Explicitly,  for any $n=1,2,\ldots,$ and $\al_1,\ldots,\al_n\in \Om^1\AR,$ we get
\beq{iL}
i_\Th(\al_1\al_2\ldots\al_n)=
\sum_{1\leq k \leq n} 
(-1)^{k-1}\cd
(\al_1\ldots\al_{k-1}\,(i_\Th'\al_k))\otimes((i_\Th''\al_k)\,
\al_{k+1}\ldots\al_n).
\eeq

It is often convenient to view the contraction
map $i_\Th$ as a map $\Om^\bullet\AR\to
T^\bullet_\k(\Om^\bullet\AR),$ and to extend the latter map
further
to get a super-derivation of the tensor
algebra $T^\bullet_\k(\Om^\bullet\AR)$.

With this understood, 
 one has the following standard identities.

\begin{lem}\label{cartan2} We have
\begin{align*}
&i_\Th(\al\beta)=(i_\Th\al)\,\beta+ 
(-1)^{\deg\al}\al\, (i_\Th\beta),
\quad\forall \al,\beta\in\Om^\bullet\AR;\\
&i_\Phi\ccirc i_\Th+i_\Th\ccirc i_\Phi=0,\quad\forall \Phi,\Th\in\dder_RA.
\end{align*}
\end{lem}

The meaning as well as the proof
of the second identity of the Lemma 
is illustrated by the following computation,
\begin{align*}
&i_\Phi\ccirc i_\Th(\al_1\al_2)=
i_\Phi\Bigl(i_\Th'\al_1\otimes (i_\Th''\al_1)\,\al_2-
\al_1\,(i_\Th'\al_2)\otimes i_\Th''\al_2\Bigr)\\
&=i_\Th'\al_1\otimes (i_\Th''\al_1)\,(i_\Phi'\al_2)\otimes i_\Phi''\al_2
-i_\Phi'\al_1\otimes (i_\Phi''\al_1)(i_\Th'\al_2)\otimes i_\Th''\al_2.
\end{align*}
The expression for  $i_\Th\ccirc i_\Phi(\al_1\al_2)$ is obtained from
this by switching the roles of $\Th$ and $\Phi$.
But the flip $\Th\leftrightarrow\Phi$ takes the
the last line in the displayed
formula above to its negative, and \eqref{cartan2} follows.\qed

\subsection{Lie derivative.}\label{lie} For any derivation
$\th\in\Der\AR$, there is a standard 
 {\em Lie derivative} operator  $L_\th$ acting on various
objects associated naturally with the algebra $A$, e.g.,
an operator
$L_\th: \dder\AR\to\dder\AR$. There is also 
 a derivation $L_\th: \Om^\bullet\AR\to\Om^\bullet\AR$,
and the corresponding induced map
 $L_\th: \DR^\bullet\AR\to\DR^\bullet\AR.$
All standard formulas involving de Rham differential, contraction
and Lie derivative operators hold in the Karoubi-de Rham setting. In particular,
one has  the Cartan identity
$L_\th=d\ccirc i_\th+i_\th\ccirc d,$ for any
$\th\in\Der\AR$.

Now, given a double derivation
$\Th\in\dder_RA$, we define
the corresponding Lie derivative map as follows:
\begin{align*}
L_\Th:\
\Om^1\AR &\too (A\otimes\Om^1\AR)\bigoplus
(\Om^1\AR\otimes A),\quad \al\mto L_\Th\al,\en\text{where}\\
L_\Th(x\,dy)&:=
\Th'(x)\otimes\Th''(x)\,dy +(x\,d\Th'(y))\otimes\Th''(y)
+x\,\Th'(y)\otimes d\Th''(y).
\end{align*}

We may naturally extend 
the map $L_\Th$  to a degree preserving derivation
$$
L_\Th:\ 
\Om^n\AR\to\bigoplus_{0\leq k \leq n}
\Om^{k}\AR\otimes\Om^{n-k}\AR,
$$
of the graded algebra $\Om^\bullet\AR$ with coefficients
in $\Om^\bullet\AR\o\Om^\bullet\AR,$ viewed
as an $\Om^\bullet\AR$-bimodule with respect to the outer 
bimodule structure.
Explicitly, we have
\begin{align}\label{L'}
L_\Th(a_0\,da_1 &\ldots \,da_n):=
\Th'(a_0)\otimes\Th''(a_0)\,da_1\ldots da_n+\\
&+\sum_{1\leq k \leq n}\Big(
a_0\,da_1\ldots
da_{k-1}\,d\Th'(a_k)\otimes
\Th''(a_k)
\,da_{k+1}\ldots da_n\Big. \nonumber\\
&\Big.\qquad\quad\en+
a_0\,da_1\ldots
da_{k-1}\,\Th'(a_k)\otimes
d\Th''(a_k)
\,da_{k+1}\ldots da_n\Big) .\nonumber
\end{align}

As in the case of contractions, it is often convenient
to view the Lie derivative
as a map $L_\Th: \Om^\bullet\AR\to
T^\bullet_\k(\Om^\bullet\AR)$, and to extend the
latter map further as a derivation
of the tensor algebra $T^\bullet_\k(\Om^\bullet\AR)$.
Similarly, given $\th\in\Der\AR$,
we  extend the  Lie derivative $L_\th$ as a derivation
of  the tensor algebra $T^\bullet_\k(\Om^\bullet\AR)$;
we may also extend
the differential $d: \Om^\bullet\AR\to
\Om^{\bullet+1}\AR$ as a {\em super}-derivation
of  the tensor algebra $T^\bullet_\k(\Om^\bullet\AR)$.

With these definitions, it is straightforward 
to verify the identities
$$[L_\xi, i_\Th]=i_{L_\xi(\Th)}
\quad\text{and}\quad
[L_\xi,L_\Th]=L_{L_\xi(\Th)},
\quad\forall \Th\in\dder\AR,\,\xi\in\Der\AR.
$$
One also has
the following {\em Cartan  formula} for double derivations:
\beq{cartan}
d\ccirc i_\Th+i_\Th\ccirc d=L_\Th,
\quad \forall \Th\in\dder_RA.
\eeq
It follows in particular that the Lie derivative commutes with
the de Rham differential: 
$$d\ccirc L_\Th=d\ccirc d\ccirc i_\Th+d\ccirc i_\Th\ccirc d=
d\ccirc i_\Th\ccirc d=d\ccirc i_\Th\ccirc d+i_\Th\ccirc d\ccirc d=L_\Th\ccirc
d.
$$

\begin{rem} It seems very likely that, given $\Th,\Phi\in\dder\AR$,
the map 
$$L_\Th\ccirc L_\Phi-L_\Phi\ccirc L_\Th:\
\Om^\bullet\AR\map \Om^\bullet\AR\otimes\Om^\bullet\AR\otimes\Om^\bullet\AR
$$
is equal to the (appropriately defined) Lie derivative
with respect to $\{\Th,\Phi\}$,
the {\em Schouten double  bracket} of $\Th$ and $\Phi,$
introduced by Van den Bergh \cite{VB2}.
\erem

\begin{proof}[Proof of Lemma \ref{homotopy}.] 
The grading on $A$ gives rise to  the {\em Euler derivation}
$\Eu: A\to A$, defined by $\Eu|_{A_k}=k\cdot\Id,\,k=0,1,\ldots.$
The action of the corresponding Lie derivative
operator $L_\Eu: \DR^\bullet\AR\to\DR^\bullet\AR$ has
nonnegative integral eigenvalues. It is clear that
the zero weight subspace is equal to
$\DR^\bullet_RA_0$, so we have a  direct
sum decomposition $\DR^\bullet\AR=
(\DR^\bullet_RA_0)\bigoplus (\DR^\bullet\AR)_{>0},$
where the second summand is spanned by the eigenspaces
corresponding to  strictly positive eigenvalues.

The direct
sum decomposition above is stable under the maps
$d$ and $i_\Eu$. Furthermore,
  the Cartan identity shows that an appropriate rescaling of 
contraction map $i_\Eu$ provides a map $h : (\DR^\bullet\AR)_{>0}
\to(\DR^{\bullet-1}\AR)_{>0}$ such that we have
$\Id= d\ccirc h+h\ccirc d$. Thus, the differential $d:
(\DR^\bullet\AR)_{>0}\to(\DR^\bullet\AR)_{>0}$
is homotopic to zero. It follows that 
the direct summand $(\DR^\bullet\AR)_{>0}$
has trivial de Rham cohomology.
\end{proof}
\subsection{Reduced contraction and Lie derivative.}
Observe that  the sign of the permutation
$(1,\ldots ,k,k+1,\ldots ,k+l)\mapsto
(k+1,\ldots ,k+l,1,\ldots ,k)$ equals $(-1)^{kl}$.
For any  $\al\o \be\in\Om^k\AR\o \Om^l\AR$,
we put 
\beq{si}
(\al\o \be)^\y:=  (-1)^{kl}\be\,\al\in
\Om^{k+l}\AR,
\eeq
 and extend the assignment
$\al\o \be\mto(\al\o \be)^\y$ to a linear
map
$ \Om^\bullet\AR\o\Om^\bullet\AR\to\Om^\bullet\AR.$
It is clear that, in $\DR^{k+l}\AR,$ one has
\beq{mom}
\lll(\alpha\otimes\beta)^\y\rrr=\lll\alpha\,\be\rrr,
\quad\forall \al\in\Om^k\AR,\,\be\in\Om^l\AR.
\eeq

Let $\Th\in\dder_RA$. For  any $n=1,2,\ldots,$
 we define
{\em reduced contraction} $\bi_\Th$, resp.,  {\em reduced Lie derivative}
$\bl_\Th$,
as the following maps, cf. \eqref{iL}, resp., \eqref{L'}:
\begin{align}\label{reduced}
&\bi_\Th:\ 
\Om^n_RA \to\Om^{n-1}_RA,\quad \al\mto \bi_\Th\al=(i_\Th\al)^\y,\quad\nonumber
\text{resp.,}\\
&\bl_\Th:\ 
\Om^n_RA \to\Om^n_RA,\quad\al\mto\bl_\Th\al=(L_\Th\al)^\y.
\end{align}

Explicitly, for any $\al_1,\al_2,\ldots,\al_n\in\Om^1\AR,$ 
using the definition of $i_\Th$,  we find
\beq{bi_formula}
\bi_\Th(\al_1\al_2\ldots\al_n)=
\sum_{k=1}^n
(-1)^{(k-1)(n-k+1)}\cd  (i_\Th''\al_k)\cd
\al_{k+1}\ldots\al_n\,\al_1\ldots\al_{k-1}\cd (i_\Th'\al_k).
\eeq
Similarly, using the definition of $L_\Th$, we find
\begin{align}\label{bl_formula}
&\bl_\Th(a_0\,da_1\,da_2\ldots\,da_n)=\Th''(a_0)\cd da_1\,da_2\ldots\,da_n\cd
\Th'(a_0)\nonumber\\
&+\sum_{k=1}^n
(-1)^{k(n-k)}\cd \Th''(a_k)\cd da_{k+1}\ldots\,da_n\,a_0\,
da_1\,da_2\ldots da_{k-1}\cd d\Th'(a_k)\\
&+\sum_{k=1}^n
(-1)^{(k-1)(n-k+1)}\cd  d\Th''(a_k)\cd da_{k+1}\ldots\,da_n\,a_0\,
da_1\,da_2\ldots da_{k-1}\cd \Th'(a_k).\nonumber
\end{align}

The main properties
of reduced contraction $\bi$  may be  summarized as follows.
\begin{lem}\label{bi_property}\vi
Fix $\Th\in\dder_RA$.
Then, for any  $\om\in\Om^n\AR$,
the element $\bi_\Th\om\in \Om^{n-1}_RA$ depends only
on the image of $\om$ in $\DR^n\AR$;
in other words,
the assignment $\om\mapsto \bi_\Th\om$
descends to a well-defined map $\bi_\Th:
\DR^n_RA\to\Om^{n-1}_RA.$
\vskip 1pt

\vii For  fixed  $\om\in \DR^n\AR$, 
the assignment $\Th\mapsto \bi_\Th\om$ gives
an $A$-bimodule morphism $\bi(\om): \dder_RA\to \Om^{n-1}_RA,$
 where $\dder\AR$ is equipped with the $A$-bimodule
structure
induced from the inner bimodule structure on $A\otimes A$. 
\vskip 1pt

\viii 
For any  $\om\in \Om^n\AR$, 
 the following
diagram commutes:
$$
\xymatrix{
(\dder_RA)\br\ar[rr]^<>(0.5){\bm_*}
\ar[d]_<>(0.5){\bi(\om)\br}&&\Der\AR\ar[d]^<>(0.5){i:\, \th\mapsto
i_\th\lll\om\rrr}&\\
(\Om^{n-1}_RA)\br\ar[rr]^<>(0.5){\op{proj}}&&\DR^{n-1}_RA.&\quad\Box
}
$$
\end{lem}
\begin{proof} Given a 1-form
$\beta\in \Om^1\AR$ and $a\in A$,  using definitions, we compute
$$\bi_\Th(a\,\beta)=i''_\Th(a\,\beta)\cd i'_\Th(a\,\beta)=
(i''_\Th\beta)\cd a\cd(i'_\Th\beta)=
i''_\Th(\beta\,a)\cd i'_\Th(\beta\,a)=
\bi_\Th(\beta\,a).
$$
It is now immediate  that, for any $\om\in\Om^n\AR$,
we have:
$$\bi_\Th(a\,\om)=\bi_\Th(\om\,a),
\quad\text{and}\en
\bi_{a'\cdot\Th\cdot a''}\om=
a'\cd(\bi_\Th\om)\cd a'',
\quad\forall\,a,a',a''\in A.
$$
Further, 
it is clear from \eqref{bi_formula} that 
$\bi_\Th(\al_1\ldots\al_n)$ and
 $(-1)^{n-1}\bi_\Th(\al_n\al_1\ldots\al_{n-1})$
are given by the same formula.
This yields parts (i) and (ii).
Part (iii) follows from \eqref{mom}.
\end{proof}

We have $(\Om^1\AR)\br=\DR^1\AR$.
Hence, for    $\om\in \DR^2\AR$, the
diagram of  Lemma \ref{bi_property}(iii),
in the special case $n=2$, reads
\beq{ii}
\xymatrix{
(\dder_RA)\br\ar[rr]^<>(0.5){\bm\br}
\ar[dr]_<>(0.5){\bi(\om)\br}&&\Der\AR
\ar[dl]^<>(0.5){i:\, \th\mapsto i_\th\lll\om\rrr}\\
&\DR^1\AR&
}
\eeq

As a general rule, any formula involving contraction
$i_\Th$, resp., Lie derivative $L_\Th$, and some other natural map
$\Om^\bullet\AR\to \Om^\bullet\AR$ (as opposed to
 a map $\Om^\bullet\AR\to \Om^\bullet\AR\otimes\Om^\bullet\AR$)
gives rise to an analogous formula  involving reduced contraction
$\bi_\Th$, resp., reduced Lie derivative $\bl_\Th$. In particular,
one has 

\begin{lem}\label{anti_comm} \vi For any $\Th\in\dder\AR$, we have
$$
d\ccirc \bi_\Th+\bi_\Th\ccirc d=\bl_\Th,\quad d\ccirc \bl_\Th=\bl_\Th\ccirc
d.
$$

\vii For any  $\xi\in\Der\AR$,
the  maps $\bi_\Th$ and $i_\xi$ anti-commute, i.e. for any $\om\in\DR^\bullet\AR,$ in
 $\Om^{\bullet-2}\AR,$ we have
$(i_\xi\ccirc\bi_\Th +\bi_\Th\ccirc i_\xi)(\om)=0.$
\end{lem}

\begin{rem}
Given
$\Th,\Phi\in \dder_RA$, it is {\em not}
true that $\bi_\Phi\ccirc\bi_\Th+\bi_\Th\ccirc\bi_\Phi=0,$
in general.
\erem

\begin{proof}[Proof of Lemma.] Both statements will be verified
by direct, somewhat tedious, computations. 
The reader is referred to \cite{Gi2} for a more conceptual
argument.

To prove (i), fix $\Th\in\dder\AR$ and
$a_0,a_1,\ldots,a_n\in A$. Using formula \eqref{bi_formula},
for $d\ccirc\bi_\Th(a_0\,da_1\ldots da_n)$ we get the following 
expression:
\begin{align*}
&d\left(\sum_{k=1}^n
(-1)^{(k-1)(n-k+1)}\cd  (i_\Th''\al_k)\cd
\al_{k+1}\ldots\al_n\,\al_0\ldots\al_{k-1}\cd (i_\Th'\al_k)\right)\\
&=
\sum_{k=1}^n (-1)^{(k-1)(n-k+1)+n-k}
\Th''(a_k)\cd da_{k+1}\ldots da_n\,da_0\,da_1\ldots
da_{k-1}\cd\Th'(a_k)\\
&+\sum_{k=1}^n (-1)^{(k-1)(n-k+1)}
d\Th''(a_k)\cd da_{k+1}\ldots da_n\,
a_0\,da_1\ldots
da_{k-1}\cd\Th'(a_k)\\
&+\sum_{k=1}^n (-1)^{(k-1)(n-k+1)+n-1}\cd
\Th''(a_k)\cd da_{k+1}\ldots da_n\,
a_0\,da_1\ldots
da_{k-1}\cd d\Th'(a_k).
\end{align*}

After simplification of powers of $(-1),$ 
we find that
$d\ccirc\bi_\Th(a_0\,da_1\ldots da_n)$ is a sum of three
terms, $S_0+$
$S_++S_-,$
where
\begin{align*}
S_0&=
\sum_{k=1}^n (-1)^{nk-1}\Th''(a_k)\cd da_{k+1}\ldots da_n\,da_0\,da_1\ldots
da_{k-1}\cd\Th'(a_k);\\
S_+&=\sum_{k=1}^n (-1)^{(k-1)(n-k+1)}d\Th''(a_k)\cd da_{k+1}\ldots da_n\,
a_0\,da_1\ldots
da_{k-1}\cd\Th'(a_k);\\
S_-&=\sum_{k=1}^n (-1)^{k(n-k)}
\cd
\Th''(a_k)\cd da_{k+1}\ldots da_n\,
a_0\,da_1\ldots
da_{k-1}\cd d\Th'(a_k).
\end{align*}

On the other hand, for
$\bi_\Th\ccirc d(a_0\,da_1\ldots da_n)=\bi_\Th(da_0\,da_1\ldots da_n),$
we find
\begin{align*}
&\bi_\Th\ccirc d(a_0\,da_1\ldots da_n)=
\Th''(a_0)\cd da_1\ldots
da_n\cd\Th'(a_0)\\
&\quad\qquad+\sum_{k=1}^n (-1)^{k(n-k+1)}\cd
\Th''(a_k)\cd da_{k+1}\ldots da_n\,da_0\,da_1\ldots
da_{k-1}\cd\Th'(a_k)\\
&=\Th''(a_0)\cd da_1\ldots
da_n\cd\Th'(a_0)\\
&\quad\qquad+\sum_{k=1}^n (-1)^{nk}\cd
\Th''(a_k)\cd da_{k+1}\ldots da_n\,da_0\,da_1\ldots
da_{k-1}\cd\Th'(a_k)\\
&=\Th''(a_0)\cd da_1\ldots
da_n\cd\Th'(a_0)-S_0.
\end{align*}

Thus, we see that 
$$(d\ccirc\bi_\Th+\bi_\Th\ccirc d)(a_0\,da_1\ldots da_n)
=\Th''(a_0)\cd da_1\ldots
da_n\cd\Th'(a_0)+S_++S_-.
$$ 
The expression on the right hand side is exactly  the 
one given by formula \eqref{bl_formula}, and part (i) is proved.

To prove (ii),
fix  $\Th\in\dder\AR$ and $\xi\in\Der\AR$, and also
1-forms $\al_1,\ldots,\al_n\in\Om^1\AR$.
For any pair of integers  $1\leq j,k\leq n$ such that $j\neq k,$
let 
$S_{j,k}$ be the $(n-2)$-form obtained from
the form
$$(i''_\Th\al_k)\cd
\al_{k+1}\ldots\al_n\,\al_{1}\ldots\al_{k-1}\cd (i'_\Th\al_k)\in \Om^{n-1}\AR$$
by replacing the factor $\al_j\in\Om^1\AR$ by the element $i_\xi\al_j\in A$.

With this notation, we
express $i_\xi\ccirc\bi_\Th(\al_1\al_2\ldots\al_n),$
using  formula
\eqref{bi_formula} as follows:
\begin{align*}
&i_\xi\left(\sum_{1\leq k\leq n}
(-1)^{(k-1)(n-k+1)} (i''_\Th\al_k)\cd
\al_{k+1}\ldots\al_n\,\al_{1}\ldots\al_{k-1}\cd (i'_\Th\al_k)\right)
\\
&=\sum_{k=1}^n
(-1)^{(k-1)(n-k+1)}\left(\sum_{j=k+1}^n
(-1)^{j-k+1}S_{j,k}
+\sum_{j=1}^{k-1}(-1)^{n-k+j-1}S_{j,k}\right).
\end{align*}
Simplifying, we obtain
\beq{1for}
i_\xi\ccirc\bi_\Th(\al_1\al_2\ldots\al_n)=
\sum_{1\leq k<j\leq n} (-1)^{nk-n+j}
S_{j,k}+\sum_{1\leq j<k\leq n} (-1)^{nk+j}S_{j,k}.
\eeq

On the other hand, we have
\begin{align}\label{2for}
&\bi_\Th\ccirc i_\xi(\al_1\al_2\ldots\al_n)=
\bi_\Th\left(\sum_{j=1}^n(-1)^{j-1} 
\al_{1}\ldots\al_{j-1}\cd (i_\xi\al_j)\cd\al_{j+1}\ldots\al_{n}\right)
\nonumber\\
&=\sum_{1\leq k<j\leq n} (-1)^{j-1+(k-1)(n-k)}S_{j,k}
+\sum_{1\leq j<k\leq n}
(-1)^{j-1+(k-2)((n-1)-(k-1)+1)}S_{j,k}\nonumber\\
&=\sum_{1\leq k<j\leq n} (-1)^{nk-n+j-1}S_{j,k}+
\sum_{1\leq j<k\leq n}
(-1)^{nk+j-1}S_{j,k}.
\end{align}
Comparison of formulas \eqref{1for} and \eqref{2for}
completes the proof of part~(ii).
\end{proof}

Observe further that, by Lemma \ref{bi_property},  the 
maps $\bi_\Th$ and $\bl_\Th$ descend to 
the following well-defined maps
 (the diagram below does {\em not} commute):
\beq{des}
\xymatrix{
\DR^\bullet\AR\ar[rr]^<>(0.5){d}\ar[d]_<>(0.5){\bi_\Th}
\ar@{.>}[drr]^<>(0.5){\bl_\Th}&&
\DR^{\bullet+1}\AR\ar[d]_<>(0.5){\bi_\Th}\\
\Om^{\bullet-1}\AR\ar[rr]^<>(0.5){d}&&
\Om^\bullet\AR.
}
\eeq

Next, fix $\Th\in\dder\AR$, let $\th=\bm_*(\Th)\in\Der\AR$
be the corresponding derivation $A\to A$,
and $L_\th: \DR^\bullet\AR\to\DR^\bullet\AR$
the Lie derivative operator induced by $\th$.
It is immediate from part (iii) of
 Lemma \ref{bi_property} that one has
a {\em commutative} triangle
\beq{des2}
\xymatrix{
\DR^\bullet\AR\ar[rr]^<>(0.5){\bl_\Th}
\ar[drr]_<>(0.5){L_\th}&&
\Om^\bullet\AR\ar[d]_<>(0.5){\op{proj}}\\
&&\DR^\bullet\AR
}
\eeq
In other words, the diagram says that,
if  $\Th\in\dder\AR$ and  $\th=\bm_*(\Th)\in\Der\AR,$
then  $\bl_\Th=L_\th$ as maps
$\DR^\bullet\AR\to\DR^\bullet\AR.$

\begin{rem} The 
reader should be warned that,
in the above setting,  it is {\em not}
true in general that $\bl_\Th=L_\th$ as maps
$\Om^n\AR\to \Om^n\AR.$
In particular, for $n=0$ we have
two {\em different} maps  $\bl_\Th,\,L_\th:\ A\to A,$
where $L_\th(a)=\Th'(a)\Th''(a)$ and
$\bl_\Th(a)=\Th''(a)\Th'(a)$.
\erem

\section{The derivation $\Delta$}\label{vdb}
\subsection{} From now on, we assume that $R$ is a finite dimensional
semisimple algebra, and $\bbr=\bbr'\otimes\bbr''\in R\o R$
is a fixed separability element, see Definition~\eqref{sep_def}.

Refining slightly an idea of Van den Bergh, see
\cite[\S3.3]{VB2}, we introduce
 the following distinguished
derivation 
$$\Delta: A\to A\otimes A,\en
a\mapsto
\Delta(a)=\sum\nolimits_i(ae_i\otimes e^i-e_i\otimes e^ia)=
a\bbr'\otimes \bbr''-\bbr'\otimes\bbr''a.
$$

We observe that $\Delta(R)=0$ since $\bbr\in (R\otimes R)^R.$ Thus,
$\Delta\in\dder_RA$. 

\begin{lem}\label{biDelta} \vi For any $\om\in\Om^\bullet\AR$,
we have
$$L_\Delta\om=\om \bbr'\otimes \bbr''-\bbr'\otimes\bbr''\om,\quad
\op{and}\quad
 \bl_\Delta\om=0.
$$

\vii Writing $[a,\om]:=a\,\om-\om\,a,$ for any $a_0,a_1,\ldots, a_n\in A,$ we have
\begin{align*}
\bi_\Delta(a_0\,da_1\ldots \,da_n)
=\sum_{k=1}^n(-1)^{(k-1)(n+1)}
\bbr'[a_k,\;da_{k+1}\ldots da_n\,a_0\,da_1\ldots da_{k-1}]\bbr''.
\end{align*}

In particular, we have $\op{Im}(\bi_\Delta)\sset[A,\Om^\bullet\AR]^R.$
\end{lem}
\begin{proof} We prove (i) in the special case
 $\om=x\,dy\,dz\in\Om^2\AR.$
We compute
\begin{align*}
&L_\Delta\om=L_\Delta(x\,dy\,dz)=(x\bbr'\otimes \bbr''-\bbr'\otimes \bbr''x)
\,dy\,dz\\
&+x\,d(y\bbr'\otimes \bbr''-\bbr'\otimes \bbr''y)\,dz
+x\,dy\,d(z\bbr'\otimes \bbr''-\bbr'\otimes \bbr''z)\\
&=x\bbr'\otimes \bbr''\,dy\,dz -\bbr'\otimes \bbr''x\,dy\,dz\\
&+
x\,dy\,\bbr'\otimes \bbr''\,dz -x\bbr'\otimes \bbr''\,dy\,dz+x\,dy\,dz\,
\bbr'\otimes \bbr''-
x\,dy\,\bbr'\otimes \bbr''\,dz \\
&=x\,dy\,dz\,\bbr'\otimes \bbr''-\bbr'\otimes 
\bbr''\,x\,dy\,dz=\om\,\bbr'\otimes \bbr''-\bbr'\otimes 
\bbr''\,\om.
\end{align*}
The computation of $L_\Delta\om$ in the general case is very similar.

Further, it is clear from formula \eqref{si}
that for any $\om\in\Om^n\AR$, we have $(\sum_i\,\om e_i\otimes
e^i)^\y= 
\sum_i\,e^i\om  e_i=(e_i\otimes  e^i\om)^\y$.
We deduce that $\bl_\Delta\om=0.$ This proves
part~(i). 

A similar direct
computation based on \eqref{bi_formula}
 yields an analogue of the formula of part (ii),
with the roles of $\bbr'$ and $\bbr''$ flipped. By the symmetry
of the separability element, we have
$\bbr''\otimes\bbr'=\bbr'\otimes\bbr''$,
and  the formula of part (ii) follows. The last claim of the Lemma follows
from that formula and \eqref{MR}.
\end{proof}

Observe that the projection $\bm_*: \dder\AR\to\Der\AR$
clearly maps the
derivation $\Delta$
to  zero.
It follows that the composite
$\DR^\bullet\AR\stackrel{\bi_\Delta}\map
\Om^{\bullet-1}\AR\stackrel{\op{proj}}\too
\DR^{\bullet-1}\AR,$
$\al\mapsto \lll\bi_\Delta\al\rrr,$
is the zero map.

\begin{cor}\label{iddi} \vi For any $\om\in\Om^\bullet\AR$,
in $\DR^{\bullet-1}\AR$, we have $\lll\bi_\Delta\om\rrr=0.$
Also, in $\Om^{\bullet-2}\AR$,  resp. in  $\DR^{\bullet-1}\AR$, we have
$$ \bi_\Th(\bi_\Delta\om)=0,
\quad\text{resp.,}\quad\bl_\Th(\bi_\Delta\om)=0,\quad\forall\Th\in\dder\AR.
$$

\vii Furthermore,   we have:
$$\bi_\Delta\ccirc d+d\ccirc \bi_\Delta=0,\quad\op{and}\en
L_\th\ccirc\bi_\Delta=\bi_\Delta\ccirc L_\th,
\quad\forall \th\in\Der\AR.$$
\end{cor}
\begin{proof}
By Lemma \ref{bi_property}(i),
the element $\bi_\Th(\bi_\Delta\om)\in \Om^{\bullet-2}\AR$
depends only on $\lll\bi_\Delta\om\rrr$,
the image of $\bi_\Delta\om$ in
$\DR^{\bullet-1}\AR$. Therefore, since $\lll\bi_\Delta\om\rrr=0,$
we conclude that $\bi_\Th(\bi_\Delta\om)=0$.
The proof for $\bl_\Th(\bi_\Delta\om)=0$
is entirely similar.

The first formula of part (ii) follows from
Lemma \ref{biDelta}(i) combined with  the Cartan formula of
Lemma \ref{anti_comm}(i).
The second formula follows from the  first one, using 
the standard Cartan formula $L_\th=d\ccirc i_\th+i_\th\ccirc d$
and Lemma \ref{anti_comm}.
\end{proof}

\subsection{Inner derivations.} 
Given an $A$-bimodule $M$ and $m\in M$, write
$\ad m: A\to M$ for the {\em inner derivation}
$\ad m(a):= am-ma$. If $m\in M^R$
then $\ad m(R)=0$, hence $\ad m\in \Der_R(A,M).$

In the special case $M=A\otimes A$,
the map
$M\to\Der(A, M),\,m\mto \ad m$ is, in effect,
a morphism of $A$-bimodules, provided
$\Der(A, M)=\Der(A, A\otimes A)$ is equipped
with the  bimodule structure induced
by the inner bimodule structure on $A\otimes A$.

Let $\bbr\in R\otimes R\sset
A\otimes A$
be
 the separability element, see \S\ref{Rsep}.
By definition,
we have  $\Delta=\ad \bbr$.

Further, using Proposition \ref{sep} we see
that, for any $p=p'\otimes p''\in
(A\otimes A)^R,$
one has that
$p=p'\otimes p''=\sum_{i\in I}\,e_ip'\otimes p''e^i=\bbr'p'\otimes p''\bbr''.$
The  last expression is nothing but the
result of the {\em inner} action of 
the element $p^{\opp}=p''\otimes p'$ on $\bbr=\bbr'\otimes\bbr''$.
It follows that the inner derivation $\ad p$
corresponding to $p=p'\otimes p''$
can be written, using the
$A$-bimodule structure on $\dder_RA$,
as $p''\cdot\Delta\cdot p'$. 

This way, one obtains the following result.

\begin{lem}\label{inn} \vi For any $p=p'\otimes p''\in
(A\otimes A)^R,$ we have $\ad p=p''\cdot\Delta\cdot p'$;
hence, $\bm_*(p''\cdot\Delta\cdot p')=\ad a,$ where $a=p'p''\in A^R.$

\vii  Inner derivations form
an $A^e$-submodule in $\dder_RA$,
and this $A^e$-submo-{dule} is generated by the
derivation $\Delta$.\qed
\end{lem}

\begin{cor}\label{ipom} For  
 any  $\om\in\DR^n\AR$ and  $p=p'\otimes p''\in (A\otimes A)^R$, in $\Om^\bullet\AR$ we have
$$
\bi_{\ad p}\om =p''\,(\bi_\Delta\om)\,p',
\quad\text{and}\quad
\bl_{\ad p}\om =dp''\,(\bi_\Delta\om)\,p'-(-1)^{\deg\om}
p''\,(\bi_\Delta\om)\,dp'.
$$

In particular, $\bl_{\ad p}(a)=0,$ for any
$a\in A=\Om^0\AR.$
\end{cor}

\begin{proof} The map $\bi(\om): \dder_RA\map\Om^{n-1}\AR$ being
a map of $A$-bimodules, for any $p=p'\otimes p''\in (A\otimes A)^R$,
we deduce 
$$\bi_{\ad p}\om=\bi_{p''\cdot\Delta\cdot p'}\om=
p''\,(\bi_\Delta\om)\,p'.
$$
This proves the first formula.

Next, using  the Cartan formula,
we compute
\begin{align*}
\bl_{\ad p}\om &=\bi_{\ad p}\,d\om+d\,\bi_{\ad p}\om
=p''\,(\bi_\Delta d\om)\,p'+d\bigl(p''\,(\bi_\Delta\om)\,p'\bigr)\\
&=p''\,(\bi_\Delta d\om)\,p'+dp''\,(\bi_\Delta\om)\,p'
+p''\,(d\bi_\Delta \om)\,p'
-(-1)^{\deg\om}
p''\,(\bi_\Delta\om)\,dp'\\
&=p''\,(\bl_\Delta \om)\,p'+dp''\,(\bi_\Delta\om)\,p'
-(-1)^{\deg\om}
p''\,(\bi_\Delta\om)\,dp'.
\end{align*}
But $\bl_\Delta\om=0$ by Lemma \ref{biDelta}(i),
and the result follows.
\end{proof}

Let $\Inn\AR\sset\Der\AR$ be the subspace of 
inner derivations. This is a {\em Lie ideal}
with respect to the commutator of derivations.
From Corollary \ref{ipom}, using 
 the commutative diagram of Lemma \ref{bi_property}(iii)
we deduce
\begin{cor}\label{ipom2}
In $\DR^\bullet\AR$, we have
$$i_{\ad a}\om=a\cd\bi_\Delta\om,\quad\op{and}\quad
L_{\ad a}\om=da\cd\bi_\Delta\om,\quad\forall
\om\in\DR^n\AR,\,a\in A^R.
$$
\end{cor}

\subsection{Geometric example.}
Let  $R=\k$ be an algebraically closed field
(of characteristic zero as usual),
and let $A=\k[X]$ be the coordinate ring of a smooth
affine algebraic variety. Then, $A\otimes A=\k[X\times X]$.
In this case, it is well-known that the Ext-group 
$\Ext_{A^e}^k(A, A\otimes A)$ vanishes for all
$k\neq\dim X$. Thus, if $\dim X >1$ then,
we have $\Der(A, A\otimes A)/\Inn(A, A\otimes A)=\Ext_{A^e}^1(A, A\otimes A)=0.$
We deduce that any derivation $\Th\in\dder A$ is {\em inner}.

Assume now that $X$ is a smooth curve ($\dim X=1$). Then,
 the bimodule $\ncO^1 A\sset A\otimes A$ is the ideal 
of the diagonal divisor $ D \sset X\times X.$ Hence  we have 
\beq{dderD}
\dder A=\Hom_{A^e}(\ncO^1 A, A\otimes A)=
\Gamma(X\times X,\,\oo_{X\times X}( D ))
\eeq
is the space of regular functions on $(X\times X)\sminus D$ with 
at most simple
poles 
along $ D .$ 
We also have algebra isomorphisms
 $\Hom_{A^e}(\Om^1A, \Om^1A)=
\G(X\times X,\oo_{X\times X})$\break
$=A\otimes A$.

It is known that, for a smooth curve $X$,
the algebra $A=\k[X]$ is smooth in the sense of Definition
\ref{smooth}, cf. \cite{CQ}.
This way, the short exact sequence of Proposition \ref{form_smooth}  becomes
the top row of the following diagram
$$
\xymatrix{
0\ar[r]&A\otimes A\ar[r]^<>(0.5){j}\ar@{=}[d]&
\dder A\ar[r]\ar@{=}[d]^<>(0.5){\Psi}&
\Der A\ar[r]\ar@{=}[d]^<>(0.5){\Phi}&0
\\
0\ar[r]&\G(\oo_{X\times X})\ar[r]&
\G(\oo_{X\times X}(D))\ar[r]&
\G(\oo_{X\times X}(D)/\oo_{X\times X})
\ar[r]&0
}
$$

In the  bottom row of the diagram above we have used shorthand notation
$\G(-)$ for $\G(X\times X,-)$; this row is obtained by applying
 the global sections functor to
the natural extension of sheaves on $X\times X$. 
The vertical isomorphism $\Phi$, in the diagram, follows from the
identification
$\Der A=\calT(X)$, with the space of regular
vector fields on $X$. 
The vertical isomorphism $\Psi$ comes from
\eqref{dderD}.

Observe  that the function $1\in A\otimes A$
corresponds under the  above identifications
to the element
$\Id_\Om\in \Hom_{A^e}(\Om^1A,\Om^1A)$. Therefore,  in the diagram we have
$j(1)=\Delta$, and the map  $j$ is nothing but
the imbedding 
$\ad :  A\otimes A\into  \Der(A,A\otimes A),$
of inner derivations.

\section{Hamilton operators.}
\subsection{Noncommuative moment map.}
\label{NCmom} We fix a semisimple finite-dimensional 
$\k$-algebra $R$ and an $R$-algebra $A$. 
Let $(A/R)^R$ be the centralizer of $R$ in the
$R$-bimodule $A/R$.

We are going to define a canonical linear 
map 
\beq{munc}
\munc:
 (\DR^2\AR)\closed\too (A/R)^R,
\quad\text{such that}\quad
d\munc(\om)=\bi_\Delta\om
\eeq
 holds in  $\Om^1\AR$, for any $\om\in(\DR^2\AR)\closed.$

To this end, we observe that the
de Rham differential 
 $d$ anti-commutes  with $\bi_\Delta$, by Corollary
 \ref{iddi}(ii).
Therefore, for any $\om\in\DR^2\AR$ such that $d\om=0$,
we have $d(\bi_\Delta\om)=-\bi_\Delta(d\om)=0$. 
Moreover, 
by the last statement of Lemma \ref{biDelta},
we know that $\bi_\Delta\om\in(\Om^1_RA)^R\closed.$

Recall further that the complex
$(\Om^\bullet\AR,\,d)$ is acyclic in positive
degrees, see \eqref{acyclic}.
 It follows that the de Rham  differential yields an isomorphism
$d: A/R\iso$
$(\Om^1_RA)\closed$. This isomorphism
clearly commutes with the adjoint action of $R$,
hence, induces an isomorphism
$(A/R)^R\iso(\Om^1_RA)^R\closed$.
Thus, we may (and will)
 define
the map $\munc$  as a composite:
\beq{muncb}\munc:\
(\DR^2\AR)\closed\stackrel{\bi_\Delta}\map
(\Om^1_RA)^R\closed=(\Om^1_RA)^R\exact
\stackrel{d\inv}\too (A/R)^R.
\eeq

The map $\munc$ will
play a crucial role in our approach to noncommutative Hamiltonian reduction.
We will see in \S\ref{rep_fun}
that the object corresponding
to $\munc$ under the representation functor
is closely related to the ordinary moment
map used in (commutative) symplectic geometry.

Here are a few basic formulas involving the moment map $\munc$.

\begin{prop}\label{omh} Let $\om$ be a  closed 2-form
and $\omh\in A$ a representative of the class
$\munc(\om)\in A/R$. Then,

\vi We have $
\omh\in R+A^R,$ 
and also $ d\omh=0$ in $\DR^1\AR.$

\vii
For any $\th\in\Der\AR$ and $\Th\in\dder\AR$, in $A/R$, one has
$$
\munc(L_\th\om)=L_\th(\munc(\om)),\quad\op{and}\quad\bl_\Th(\munc(\om))= 0.
$$

\viii For any $p=p'\otimes p''\in (A\otimes A)^R,\, a\in A^R,$ and $u,v\in A$, one has
\begin{align*}
&\bi_{\ad p}\om =p''\cdot d\omh \cdot p'\en\quad\mbox{$\text{holds in}$}\en\Om^1\AR,\\
&i_{\ad a}\om =a\cdot d\omh\;\en\qquad\mbox{$\text{holds
in}$}\en\DR^{1}\AR;\\
&\munc(du\,dv)
=\bbr'\cdot [u,v]\cdot\bbr''\,\op{mod}\, R.
\end{align*}
\end{prop}

\begin{proof} 
By Corollary \ref{ipom2}, we know that $\bi_\Delta\om=0$
in $\DR^1\AR$, hence $d\omh=\bi_\Delta\om=0$.
Observe next that, since
$\munc(\om)\in (A/R)^R$, 
we have $[\omh,R]\sset R$.
Hence the derivation $\ad \omh$ preserves $R$.
But any derivation of the  finite dimensional semisimple 
algebra $R$ is inner.
Hence, there exists $r\in R$ such
that $\ad\omh|_R=\ad r|_R$.
Thus, $\omh-r$ commutes with $R$,
and (i) follows.

Part (ii) and  the first two
equations of part (iii)
 are immediate from the
properties of $\bi_\Delta$ established earlier.

To prove the last formula, using Lemma \ref{biDelta}(ii)
we compute
\begin{align*}
\bi_\Delta(du\,dv)&=\bbr'(-[dv,u]+[du,v])\bbr''
=\bbr'([du,v]+[u,dv])\bbr''
=\bbr'(d[u,v])\bbr''.
\end{align*}
The last expression equals
$d\bigl(\bbr'\cdot [u,v]\cdot\bbr''\bigr).$
Thus, in $\Om^1\AR,$ we obtain
$d\munc(du\,dv)=$
$\bi_\Delta(du\,dv)=d\bigl(\bbr'\cdot [u,v]\cdot\bbr''\bigr),$
hence $\munc(du\,dv)-\bbr'\cdot [u,v]\cdot\bbr''\in R$,
and we are done.
\end{proof}

Recall the notation $V^R:=V\cap M^R$ for any vector
subspace $V$ in an $R$-bimodule $M$.
 
\begin{prop}\label{ham_lemma}
Let  $A$ be an   $R$-algebra such that the following
sequence
is exact
\beq{exact}
0\map R\map \DR^0\AR\stackrel{d}\map\DR^1\AR.
\eeq

Then, the  map $\dis\munc: (\DR^2\AR)\closed\to
(A/R)^R$ can be  lifted
canonically to a map $\widetilde{\munc}$
in the  diagram below
\beq{munctil}
\xymatrix{
\DR^1\AR\ar[rr]^<>(0.5){d}\ar@{->>}[d]_<>(0.5){\bi_\Delta}&&
(\DR^2\AR)\closed
\ar@{.>}[dll]_<>(0.5){\widetilde{\munc}}\ar[d]_<>(0.5){\bi_\Delta}\\
{{[A,A]^R}^{^{}}\;}\ar@{^{(}->}[rr]^<>(0.5){d}&&[A,\Om^1\AR]^R.
}
\eeq

 The map  $\widetilde{\munc}$ has the following properties:\smallskip

\npb{The above diagram commutes and 
we have $\munc(\om)=\widetilde{\munc}(\om)\,\op{mod}\,R;$}

\npb{The map  $\widetilde{\munc}$ commutes with
the Lie derivative $L_\th$, for any $\th\in\Der\AR$.}
\end{prop}

The exactness of \eqref{exact}
 is equivalent to the following two equalities:
\beq{ham_lemma1}
\Ker[d:\ \DR^0\AR\to \DR^1\AR]=R\quad\text{and}
\quad
[A,A]\cap R=0.
\eeq

\begin{proof}[Proof of Proposition.] Let $\om\in\DR^2\AR$ be a closed
2-form and $\omh\in A$ a representative of the class
$\munc(\om)\in A/R$. 
We know that  $d\omh=0$ in $\DR^1\AR$,
by Proposition \ref{omh}(i). 
Hence, the first equality in \eqref{ham_lemma1}
yields $\omh\in R+[A,A]$.
Now,  the
second  equality in \eqref{ham_lemma1}
implies that there is a unique representative of the
class $\munc(\om)$ that belongs to $[A,A]$.
This provides a lift $\widetilde{\munc}: (\DR^2\AR)\closed\to
[A,A]$.

Further, for any derivation $\th\in\Der\AR$,
we have
 $\th([A,A])\sset[A,A]$ and $\th(R)=0.$
 Thus, for any
choice of  $\omh\in R+[A,A]$
we get $\th(\omh)\in\th(R+[A,A])=\th([A,A])\sset[A,A].$
Therefore, we deduce
$L_\th\widetilde{\munc}(\om)-\widetilde{\munc}(L_\th\om)
\in [A,A].$
On the other hand, since $\munc$ commutes with $L_\th$,
we have
$L_\th\widetilde{\munc}(\om)-\widetilde{\munc}(L_\th\om)
\in R.$ It follows that
$L_\th\widetilde{\munc}(\om)-\widetilde{\munc}(L_\th\om)
\in R\cap [A,A]=0$, by the  
second  equality in \eqref{ham_lemma1}.
We conclude that the
map $\widetilde{\munc}$ commutes with~$L_\th.$

Observe finally that the  
second  equality in \eqref{ham_lemma1} implies, in particular,
that $R$ is commutative.
Hence, from Proposition \ref{omh}(i) we deduce $[\omh,R]\sset[R,R]=0.$
Thus, $\omh\in A^R$, and we get $\widetilde{\munc}(\om)\in[A,A]\cap A^R=[A,A]^R$.
 This completes the
proof.
\end{proof}

\subsection{Symplectic 2-forms.}\label{2forms} 
We say that a 2-form $\om\in\DR^2\AR$
is {\em nondegenerate}
if the  map $i(\om):\Der\AR\to\DR^1\AR,\,\th\mapsto i_\th\om,$
is a bijection.
Given a nondegenerate form $\om\in\DR^2\AR$,
we invert the  bijection $i(\om)$ to obtain a bijection
$\SH_\om=i(\om)\inv: \DR^1\AR$ $\iso\Der\AR.$

We say that a 2-form $\om\in\DR^2\AR$
is  {\em bi-nondegenerate} if the  map $\bi(\om):$
$\dder\AR\to\Om^1\AR,\,\Th\mapsto\bi_\Th\om$ is a bijection,
hence, an isomorphism of $A$-bimodules.

Let  $\om\in\DR^2\AR$ be a bi-nondegenerate 2-form. 
Applying the functor $(-)\br$ to  the  map $\bi(\om)$ yields the
 bijection
$\bi(\om)\br: (\dder\AR)\br\iso (\Om^1\AR)\br=\DR^1\AR$.
Further, using the commutativity of diagram \eqref{ii}, we deduce
that the map $\bm\br:
(\dder\AR)\br\to \Der\AR$  is {\em injective}; moreover,
one has the following
diagram
$$
\xymatrix{
&\Der\AR=\op{Im}(\bm\br)\oplus\Ker(i(\om))\ar@<0.5ex>[d]\ar@{->>}[dr]^<>(0.5){i(\om)}&\\
{(\dder\AR)\br}\ar[r]_<>(0.5){\bm\br}^<>(0.5){\sim}\ar@{^{(}->}[ur]^<>(0.5){\bm\br}&
\op{Im}(\bm\br)\ar@<0.5ex>[u]\ar[r]_<>(0.5){i(\om)}^<>(0.5){\sim}&
\DR^1\AR.
}
$$
In this diagram, the vertical maps are the natural imbedding and 
projection (along $\Ker(i(\om))$), respectively, and the horizontal maps 
$\bm\br$ and $i(\om)$ are both bijections
whose composite equals  $\bi(\om)\br$.
It follows in particular that, for
a smooth algebra $A$,
any bi-nondegenerate form is automatically nondegenerate.

Next, we invert the bijection $\bi(\om)\br$ and  introduce
 the following composite
\beq{SH}
\SH_\om:\
\xymatrix{
\DR^1\AR\ar[rr]^<>(0.5){\bi(\om)\br\inv}_<>(0.5){\sim}&&
(\dder_RA)\br\;\ar@{^{(}->}[rr]^<>(0.5){\bm\br}&&
\Der\AR.}
\eeq
We see from the diagram above that the map $\SH_\om$ is injective,
furthermore, it provides a section of the projection
$i(\om): \Der\AR\onto\DR^1\AR.$
It follows in particular that we have 
\beq{eqn}
i_{\SH_\om(\al)}\om=\al,\quad\forall \al\in\DR^1\AR.
\eeq

Both in the nondegenerate and  bi-nondegenerate
cases, we also consider the
 composite
\beq{th_p}
A\br\stackrel{d}\map\DR^1\AR\stackrel{\SH_\om}\too
\Der\AR,\quad p\mto \th_p:=\SH_\om(dp).
\eeq

\begin{prop}\label{thomh}  Let $\om$ be a closed 2-form and
 $\omh\in A$ a representative of the class
$\munc(\om)\in A/R$. Let
$\th\in\Der\AR$ and put
$\al=i_\th\om$. Then

\vi We have
$\dis\en\th(\omh)=-
\bi_\Delta\al.$

\vii The 1-form  $\al\in\DR^1\AR$ is closed if and only if $L_\th\om=0.$
\end{prop}
\begin{proof} Let $\th\in\Der\AR$. 
The maps $i_\th$ and $\bi_\Delta$ anti-commute,
by Lemma \ref{anti_comm}(ii). Hence, using \eqref{eqn} we find
$\th(\omh)=i_\th(\bi_\Delta\om)=-\bi_\Delta(i_\th\om)=-\bi_\Delta\al,$
and (i) follows.
To prove (ii), we write
$L_\th\om=di_\th\om +i_\th
d\om=di_\th\om,$
since $d\om=0$.
We see that   $L_\th\om=0$
 if and only if $di_\th\om=0$.
\end{proof}

Fix a closed 2-form 
 $\om\in \DR^2\AR$.
A derivation $\th\in\Der\AR$ is said to be $\om$-{\em symplectic}
if in $\DR^2\AR$ one has $L_\th\om=0$. Clearly,
$\om$-symplectic derivations form a Lie subalgebra
 $\Der^\om\AR\sset\Der\AR$.

\begin{defn}\label{def}
A closed 2-form $\om\in\DR^2\AR$ is called
{\em symplectic}, resp.,
{\em bi-symplectic}, if it is nondegenerate,
resp., bi-nondegenerate.
\end{defn}

\begin{lem}\label{H_1}   Let $\om$ be a closed 2-form and
 $\omh\in A$ a representative of the class
$\munc(\om)\in A/R$. Then

\vi For any $\om$-symplectic derivation
$\th\in \Der^\om\AR$ we have $\th(\omh)\in R$.

\vii If $\om$ is either symplectic or bi-symplectic then,
for any $a\in A\br$, we have $\th_a\in\Der^\om\AR;$
furthermore, 
$\th_a(\omh)=0$.

\viii If the sequence \eqref{exact} is exact then 
$\dis\en\th\in\Der^\om\AR$ implies
$\th(\omh)=0$.
\end{lem}
\begin{proof}
Let
$\th\in \Der^\om\AR$.
Since $L_\th$ commutes with $\bi_\Delta$,
from Proposition \ref{thomh}(ii)
we deduce
$L_\th\bi_\Delta\om=\bi_\Delta L_\th\om=0.$
The map $d$, hence $d\inv$,
also commutes with $L_\th$.
Thus, in $A/R$, we get 
$L_\th\munc(\om)=L_\th
d\inv(\bi_\Delta\om)=d\inv(L_\th\bi_\Delta\om)=0$.
This means that  $L_\th\omh=0\,\op{mod}\,R,$ hence
$\th(\omh)\in R$, and (i) follows.

Since $da$ is a closed 1-form for any $a\in A$,
we deduce from  Proposition
\ref{thomh}(ii) that $\th_a\in \Der^\om\AR.$
Further, in the special case $n=1$, the  formula of Lemma
\ref{biDelta}(ii) says that the map
$\bi_\Delta: \Om^1\AR\to [A,\Om^0\AR]^R=[A,A]^R$ is given
 by  $\bi_\Delta(a_0\,da_1)=\bbr'[a_1,a_0]\bbr''$;
in particular, we have $\bi_\Delta(da)=0$, for any $a\in A$.
Therefore, since $i_{\th_a}\om=da$, 
we conclude  that $\th_a(\omh)=-\bi_\Delta(da)=0$,
by  Proposition
\ref{thomh}(i).
Part (ii) follows.

Finally, let
$\th\in\Der^\om\AR.$ 
Then $\th(\omh)\in R$ by (i). Thus, in the setting
 of Proposition \ref{ham_lemma}, for $\omh=\widetilde{\munc}(\om),$
we obtain
$\th(\omh)\in R\cap \th([A,A])\sset R\cap [A,A]$ $=0.$
This proves part (iii).
\end{proof}

\begin{rem} Let $\omh\in A$
be a representative of the class
$\munc(\om)\in A/R$. Then,
 for any derivation
$\th\in\Der\AR$, the element
$\th(\omh)=i_\th(\bi_\Delta\om)$ is independent
of the choice of  $\omh$.
\erem

\subsection{Hamilton operators.}\label{sec_GD} We are going to introduce
the notion of {\em Hamilton operator}, motivated
by a similar construction used
by Gelfand and Dorfman in their work on integrable systems, 
cf. \cite{GD} and \cite{Do}.

Let $\GG$ be a (possibly infinite dimensional) Lie algebra
and  $(C^\bullet, d)$ a complex
 with
differential $d: C^\bullet\to C^{\bullet+1}.$
By a  $\GG$-{\em equivariant structure} on $C^\bullet$ we
mean a pair of linear maps
$\GG\map \Hom_\k(C^\bullet,C^\bullet),\, x\mapsto L_x,$
and $\GG\map \Hom_\k(C^\bullet,C^{\bullet-1}),\,
x\mapsto i_x,$ which satisfy the standard commutation relations:
$$[L_x,L_y]=L_{[x,y]},\quad L_x=d\ccirc i_x+i_x\ccirc d,\quad
i_x\ccirc i_y+i_y\ccirc i_x=0,
\quad [L_x,i_y]=i_{[x,y]},
$$
for any $x,y\in\GG.$
In this case, we call  $(C^\bullet, d)$ a
$\GG$-{\em equivariant complex}.

\begin{defn}\label{Hdefn} Let  $(C^\bullet, d)$ be a
$\GG$-{\em equivariant complex}.
A linear map $\SH: C^1\map
\GG$ is said to
be {\em skew-symmetric} if one has
\beq{skew}
i_{\SH(\al_1)}\al_2+ i_{\SH(\al_2)}\al_1=0,\quad\forall \al_1,\al_2\in
C^1.
\eeq

A skew-symmetric map $\SH$ is called a {\em
 Hamilton operator} if  one has
\beq{jacobi}
i_{\SH(\al_3)}L_{\SH(\al_1)}\al_2 +\text{cyclic permutations of}\;
\{1,2,3\}\;=\;0,
\quad\forall \al_1,\al_2,\al_3\in C^1.
\eeq
\end{defn}

Given a  skew-symmetric map $\SH,$ we introduce the following
 bilinear
pairings
\begin{align}\label{Hbracket}
&\mathsf{(0)} \quad\{-,-\}^0_\SH:\ C^0\times C^0\too C^0,\quad
(p,q)\mto \{p,q\}^0_{\SH}:= i_{\SH(dp)}dq;
\nonumber\\
&\mathsf{(1)} \quad\{-,-\}_\SH^1:\ C^1\times C^1\too C^1,\\
&\qquad \qquad(\al,\be)\mto \{\al,\be\}_\SH^1:=
i_{\SH(\al)}\ccirc d\be-
i_{\SH(\be)}\ccirc d\al +d\ccirc i_{\SH(\al)}\be
\nonumber
\end{align}

The main results about  Hamilton operators,
due to Gelfand and Dorfman, may be summarized as follows.
\begin{prop}\label{GD} 
\vi A skew-symmetric map $\SH$ is  a Hamilton operator if and only if 
 one has
$$[\SH(\al),\SH(\be)]= \SH\bigl(\{\al,\be\}_\SH^1\bigr),
\quad\forall \al,\be\in C^1.
$$

\vii For any  Hamilton operator $\SH$, the bilinear pairing given by
the first, resp. second, formula in \eqref{Hbracket}
makes $C^0$, resp., $C^1$, a Lie algebra.

\viii With these Lie algebra structures, the  maps
$\dis C^0\stackrel{d}\too C^1\stackrel{\SH}\too \GG
$
are both  Lie algebra homomorphisms.\qed
\end{prop}

The Lie bracket on $C^0$ given by
the first  formula in \eqref{Hbracket} is usually referred to
as the {\em Poisson bracket} on $C^0$ induced by $\SH$.

A basic example of Hamilton operators is provided
by Poisson bivectors.
Specifically, let $X$ be a  smooth manifold,
and write $\calT(X)$ for  the Lie algebra of vector fields on $X$
and $\Om^\bullet(X)$ for the graded space 
of differential forms on $X$.
The de Rham complex, $\bigr(\Om^\bullet(X), d\bigl)$, has an
obvious structure of $\calT(X)$-equivariant complex.
Now, to each bivector
 $\pi\in\G(X, \wedge^2\calT_X)$ one associates
 a skew-symmetric map $\SH_\pi: \Om^1(X)\to\calT(X),\, \al\mapsto i_\pi\al$.
This map is a Hamilton operator if and only if
the  bivector
 $\pi$ defines a  Poisson structure on~$X$,
which holds  if and only if $\pi$ has vanishing
Schouten bracket:
$[\pi,\pi]=0$.

\subsection{Hamilton operator arising from a bi-symplectic form.}\label{NC}
The  $\calT(X)$-equivariant complex of differential forms on 
a manifold $X$,
considered above, may be generalized to the noncommutative
setup. Specifically, for any  $R$-algebra $A$,
the Karoubi-de Rham complex $\DR^\bullet\AR$
 has an
obvious structure of $\Der\AR$-equivariant complex.

Any symplectic or bi-symplectic 2-form
$\om\in\DR^2\AR$ gives rise to a Hamilton operator
$\SH_\om:  \DR^1\AR\to\Der\AR$. This is the map $\SH_\om$ introduced in Sect.
\ref{2forms}, see \eqref{SH}.
One can verify that the equation $d\om=0$ implies
condition \eqref{jacobi} in the definition of
Hamilton operator.

According to Proposition \ref{GD},
the Hamilton operator $\SH_\om: \DR^1\AR\to\Der\AR$
gives the space $A\br=\DR^0\AR$ a Lie algebra
structure. Thus, we obtain the following result, first
discovered by Kontsevich, cf. \cite{Ko},\cite{BLB},\cite{Gi}:
 
\begin{prop}\label{neck} For any symplectic or  bi-symplectic 2-form $\om\in\DR^2\AR$,
the pairing
$(p,q)\mto \{p,q\}:=\th_p(q)\,\op{mod}\,[A,A]$ gives $A\br$ a Lie algebra
structure.

With this  structure, the map
$A\br\to \Der^\om\AR,\,p\mto \th_p$
becomes a Lie algebra homomorphism.\qed
\end{prop}

Recall that we have  a moment map
$\munc: (\DR^2\AR)\closed\to A/R.$
Let
 $\omh\in A$ be a representative of the class $\munc(\om)\in A/R$, and
write $R[\omh]$ for the subalgebra of $A$ generated by
$R$ and $\omh$.
Notice that the subalgebra $R[\omh]$
is independent of the choice of $\omh$.
Notice further that, by Proposition \ref{omh}(i), we may always
find a representative
$\omh$
that commutes with $R$. This justifies the notation
$R[\omh]$.

\begin{lem}\label{hhh}
The image of the composite $R[\omh]\into A\onto A\br$ 
is contained in the center of the Lie algebra $A\br$.
\end{lem}

\begin{proof}
For any $p\in A\br,\, r',r''\in R,$ and  $m=1,2,\ldots,$ we compute
$$\{p, r'\cd\omh^m\cd r''\}=\th_p(r'\cd\omh^m\cd r'')=
r'\cd\th_p(\omh^m)\cd r''
=r'\cd\left(\sum_{k=0}^m\,
\omh^k\cd\th_p(\omh)\cd\omh^{m-k}\right)\cd r''.$$
Each term in the last sum vanishes since $\th_p(\omh)=0$,
 by Lemma \ref{H_1}(ii).
\end{proof}

Next, write
 $A\omh A\sset A$ for the two-sided ideal generated by~$\omh.$
The algebra $A_\omh:=A/A\omh A$ may be thought of as 
a Hamiltonian reduction of $A$.

The following important result will be proved in Sect. \ref{proof}. 

\begin{prop}\label{br_cor}
\vi For any $p\in A\br$ the derivation
$\th_p$ preserves the ideal $A\omh A,$ hence induces
a well-defined derivation $\bth_p: A_\omh\to A_\omh.$

\vii The Lie bracket 
on $A\br$ descends to a well-defined  Lie bracket
on $(A_\omh)\br$.

\viii There is a Lie algebra map $\b\varphi$
making the following diagram commute
$$
\xymatrix{
A\br\ar[rr]^<>(0.5){p\mapsto \bth_p}
\ar@{->>}[d]^<>(0.5){\op{proj}}&&\Der_R(A_\omh)\ar@{->>}[d]^<>(0.5){\op{proj}}\\
(A_\omh)\br\ar[rr]^<>(0.5){\b\varphi}&&
\Der_R(A_\omh)/\Inn_R(A_\omh).
}
$$
\end{prop}

\subsection{Hamilton operators via Hochshild homology.}\label{ham_hoch} 
Another class of examples of equivariant complexes
arises from Hochschild homology.

In more detail, given
 an $R$-algebra $A$ and an $A$-bimodule $M,$ let $H^k_R(A,M)$,
resp.  $H_k^R(A,M),$ denote the $k$-th
Hochschild cohomology, resp. Hochschild  homology, group  (with coefficients in $M$)
of the algebra $A$ relative to the subalgebra $R$.
In particular, we have
$H^0_R(A,M)=M^A$, the centralizer of $A$ in $M$,
and $H_0^R(A,M)=M/[A,M],$ the commutator quotient of $M$.
Also, we have
 $H^1_R(A,M)=\Der_R(A,M)/\Inn_R(A,M).$

\begin{rem}\label{H_R} If the algebra $R$  is
semisimple, then one has a canonical isomorphism
$H^\bullet_R(A,M)\cong
H^\bullet(A,M),$ between relative and absolute Hochschild cohomology
of $A$, cf. \cite[\S1.2.13]{Lo}. 
However, this isomorphism will  play no role below. 
\erem

In the special case $M=A$, we shall use  simplified notation
$\HH^k_R(A):=H^k_R(A,A),$ resp.  $\HH_k^R(A)=H_k^R(A,A).$

We recall  that Hochschild cohomology of
an associative algebra has a natural  Gerstenhaber bracket
$\HH^{k-1}_R(A)\times \HH^{l-1}_R(A)\map \HH^{k+l-1}_R(A)$.
This bracket, together with cup-product, give $\HH^\bullet_R(A)$
the structure of a  {\em Gerstenhaber algebra}.
The  Gerstenhaber bracket on $\HH^1_R(A)$ is
the one induced by the commutator bracket on
derivations. 

Further, there is a natural contraction-pairing
$i: \HH^k_R(A)\times \HH_l^R(A)\to \HH_{l-k}^R(A)$ 
and Lie derivative pairing
$L:\HH^k_R(A)\times \HH_l^R(A)\to \HH_{l-k+1}^R(A)$ which
make
Hochschild homology a module over the
 Gerstenhaber algebra $\HH^\bullet_R(A)$.
Furthermore, 
the resulting Lie derivative action of the Lie algebra
$\HH^1_R(A)$ on Hochschild homology
 may be promoted to
a structure of  $\HH^1_R(A)$-equivariant complex on
$\HH_\bullet^R(A)$, with the Connes differential
$\BB: \HH_\bullet^R(A)\to\HH_{\bullet+1}^R(A)$
playing the role of  differential, see \cite{GDT}, \cite[Ch.\,3]{Lo}.

In  particular, in the special case $l=0$ we get  an action of
the  Lie algebra
$\HH^1_R(A)$ on $A\br$. It is easy to verify that this action
is induced by the tautological $\Der\AR$-action on $A$;
the  induced  $\Der\AR$-action  on $A\br$ clearly descends to
$\Der\AR/\Inn\AR$ since any inner derivation
maps $A$ to $[A,A]$.

To construct Hamilton operators for this equivariant complex,
we exploit the main idea due to Van den Bergh \cite{VB1}.

To this end, observe that the inner bimodule
structure on $A\o A$ survives in $H^\bullet_R(A,A\o A)$,
hence, makes each cohomology
group $\HH^k_R(A,A\o A)$
an $A$-bimodule.
Now,  fix an integer $d>0$ and
assume that $\HH^k_R(A,A\o A)=0$
unless $k=d$. Assume, in addition, that $A$ has finite
Hochschild dimension, cf. Correction to \cite{VB1}.
Then, according to (a relative version 
of) \cite[Thm.~1]{VB1},
one has  {\em canonical} duality isomorphisms
\beq{vdb_iso}
H_k^R\bigl(A,\,H^d_R(A,A\o A)\bigr)\cong\HH^{d-k}_R(A), \quad
\forall k=0,1,\ldots.
\eeq

For any {\em central} element
$\pi\in H^d_R(A,A\otimes A)^A$,  we introduce the following
composition of $A$-bimodule maps
\beq{unit}
\SH_\pi: \xymatrix{
{\HH_k^R(A)}\ar[r]^<>(0.5){a\mapsto a\cdot\pi}&
{H_k^R\bigl(A,\,H^d_R(A,A\o A)\bigr)\,}
\ar[r]^<>(0.5){\eqref{vdb_iso}}_<>(0.5){\sim}&{\,\HH^{d-k}_R(A)}.
}
\eeq

A {\em central} element $m\in M^A$ of 
an $A$-bimodule $M$ is said to be
 a {\em unit} if
the map $a\mapsto a\cdot m=m\cdot a$ yields
an  $A$-bimodule isomorphism $A\iso M$.
Thus, any bimodule that has a unit is
non-canonically isomorphic
to $A$ as an $A$-bimodule.

Now, assume that one has a  unit $\pi\in H^d_R(A,A\o A)^A$.
Then, the composite map in \eqref{unit}
provides an $A$-bimodule isomorphism
$\SH_\pi:\HH_\bullet^R(A)\iso \HH^{d-\bullet}_R(A).$
We may use this isomorphism to transport 
the Connes differential $\BB: \HH_\bullet^R(A)\to
\HH_{\bullet+1}^R(A)$ to
a differential $\partial_\pi: \HH^\bullet_R(A)\to
\HH^{\bullet-1}_R(A)$.

\begin{claim}\label{BV} \vi For any  $\pi\in H^d_R(A,A\o A)^A,$
the morphism  in \eqref{unit}
intertwines contraction  and  cup-product maps, i.e., we have
$$ \SH_\pi(i_\eta\alpha)=\eta\cup\SH_\pi(\alpha),
\quad\forall \alpha\in \HH_\bullet(A),\,\eta\in \HH^\bullet(A).
$$

\vii Assume in addition that the algebra $A$ has finite
Hochschild dimension and there  is an
 integer $d>0$ such that
one has  $A$-bimodule isomorphisms (Gorenstein type property):
\beq{vdb2}
H^k_R(A,A\otimes A)\cong
\begin{cases} A &\op{if}\enspace k=d\\
0&\text{else}.
\end{cases}
\eeq

Then, for any choice of
{\sf{unit}}  $\pi\in H^d_R(A,A\o A)^A$, the
 differential $\partial_\pi$ makes Hochschild cohomology,
equipped with 
the standard Gerstenhaber algebra  structure, into 
a Batalin-Vilkovisky
algebra, in other words, for the Gerstenhaber bracket
 on Hochschild cohomology
one has
$$\{u,v\}= \partial_\pi(u\cd v)-\partial_\pi(u)\cd v -
(-1)^k u\cd \partial_\pi(v),
\quad \forall u\in \HH^k_R(A), v\in \HH^l_R(A).$$
\end{claim}

We defer the proof of this Claim to a separate publication.

\begin{examp}
Let  $A=\k[X]$ be  the
coordinate ring of a smooth affine algebraic variety $X$
of dimension
$d:=\dim X.$
In this case we have $\HH^k_R(A,A\o A)=0$
unless $k=d.$
A  choice of
unit  $\pi\in H^d_R(A,A\o A)^A=\Gamma(X,\wedge^d\calT_X)$ corresponds to
a choice of a nowhere vanishing section
of the {\em canonical sheaf} $K_X=\Om_X^d$.
This section gives a trivialization of
 $K_X$. Thus, $X$ is a Calabi-Yau manifold,
and the  Batalin-Vilkovisky structure of
 Claim \ref{BV}(ii)
is nothing but the standard   Batalin-Vilkovisky structure
on the space $\bigoplus_{k\geq 0} \Gamma(X,\wedge^k\calT_X)$
associated with the chosen volume form $\pi\inv\in  \Om_X^d$,
cf. e.g. \cite{Sch}.\hfill$\lozenge$
\end{examp} 

The special case where property \eqref{vdb2} holds for $d=2$
leads to a construction of Hamilton operators.
Specifically, one has the
following result that may be thought of as a homological 
analogue of the construction of Poisson brackets
used in \cite{VB2}.
\begin{prop}\label{claim} Let  $A$ be an $R$-algebra of finite
Hochschild dimension such that 
condition \eqref{vdb2} holds for $d=2.$

Then,  for any choice of
{\sf{unit}}  $\pi\in  H^2_R(A,A\o A)^A$,
the corresponding map $\SH_\pi: \HH_1^R(A)\iso\HH^1_R(A)$
in \eqref{unit}
is a  Hamilton operator
for $\HH_\bullet^R(A)$, viewed as an $\HH^1_R(A)$-equivariant complex.

In particular, the space $A/[A,A]$ acquires a natural
Lie bracket $\{-,-\}_\pi$, cf. Proposition \ref{GD}(ii).
\end{prop}

The proposition above easily follows from Claim
\ref{BV}. To see this, let $\SH:=\SH_\pi$, pick
$\al,\be\in \HH_1(A)$, and put $u:=\SH(\al)$ and $v=\SH(\be)$.
We apply  the isomorphism $\SH\inv:\HH^1(A)\iso\HH_1(A)$ to the displayed formula in
 Claim
\ref{BV}(ii). Using part (i) of  Claim
\ref{BV} and writing $[-,-]_\text{Gerst}$ for the
commutator in $\HH^1(A)$,  we obtain 
$$\SH\inv\big([\SH(\al),\SH(\be)]_\text{Gerst}\big)=
\BB(i_{\SH(\al)}\be) + i_{\SH(\al)}\BB(\be)-i_{\SH(\be)}\BB(\al).
$$
The right hand side here is nothing but the
bracket $\{\al,\be\}_{\SH}^1$ introduced in 
\eqref{Hbracket}. We see that the above equation
is exactly the one from  Proposition \ref{GD}(i).\qed

\section{Noncommutative cotangent bundle.}\label{cot}
\subsection{}\label{cot_main} Following an idea of \cite{CB2},
for any $R$-algebra $B$, we define ${T^*B}:=T_B(\dder_RB),$
the tensor algebra of the
$B$-bimodule $\dder_RB$. Thus, ${T^*B}=\oplus_{k\geq 0}\D_kB$,
is a graded $B$-algebra such that $\D_0B=B$. The
algebra ${T^*B}$ may be thought of as the coordinate ring
of the `noncommutative cotangent bundle'
on  $\Spec B$, a `noncommutative space'. Indeed,
we will see later
that the representation functor takes the algebra ${T^*B}$ 
to the cotangent bundle on the representation scheme
for the algebra $B$.

Write $\Delta_B$
 for the $\Delta$-derivation
$\Delta_B: B\to B\o B$; this derivation
 will be frequently
viewed as an element
of $\D_1B=\dder_RB.$
We also have 
the $\Delta$-derivation
$\Delta_{{T^*B}}:  {T^*B}\to {T^*B}\o{T^*B}$, and
 an {\em Euler derivation}
$\Eu: {T^*B}\to {T^*B}$ arising from the grading on ${T^*B}$,
thus, $\Eu(a)=k\cdot a,$ for any $a\in \D_kB$.

\begin{thm}\label{DBthm} Assume that $B$ is smooth. Then,  the algebra
${T^*B}$ is also smooth and it has a canonical
bi-symplectic 2-form $\om\in\DR^2_R({T^*B})$, such that
$L_\Eu\om=\om,$ and such that
the derivation $\Delta_B\in \D_1B=\dder_RB$ is a representative
of the class $\munc(\om)$.
\end{thm} 

\begin{rem} The bi-symplectic structure on $T^*B$ is closely related
to the `double bracket' on  $T^*B$ introduced by
 Van den Bergh
\cite{VB2}.

In more detail, for any $p\in T^*B$,
let $\Th_p: T^*B\to T^*B\o T^*B$ denote the double derivation
corresponding to the 1-form $dp$ under the isomorphism
$\dder_R(T^*B)\cong \Om^1_R(T^*B)$ provided by the
bi-symplectic form of Theorem \ref{DBthm}.
Given $\Phi,\Psi\in \dder_RB,$
view these  double derivations of $B$ as two
degree 1  elements of  $T^*B$, and let $p:=\Phi$. Then, we have a well-defined
element $\Th_\Phi(\Psi)\in  T^*B\o T^*B$.
This
 is again an element
of degree 1 with respect to the induced grading on the tensor product, i.e.,  an element
of $(B\o \dder_RB) \oplus (\dder_RB\o B).$
The map $T^*B\o T^*B\to T^*B\o T^*B,
\,\Phi\o\Psi\mto \Th_\Phi(\Psi)$ can be shown,   see Appendix to \cite{VB2}, to be equal
to  Van den Bergh's
double bracket.
\erem

Theorem \ref{DBthm}  combined with Proposition 
\ref{neck} provides the graded space  $({T^*B})\br$
with  a natural Lie algebra structure
such that the Lie bracket has degree $(-1)$.
Using the explicit construction of the bi-symplectic
2-form  of the Theorem 
given in Sect. \ref{omega} below,
one proves

\begin{cor}\label{deg1} \vi For any $R$-algebra
$A$, the space $(\dder\AR)\br$ has a natural
 Lie algebra structure.

\vii With this Lie algebra structure on $(\dder\AR)\br$
and the standard Lie algebra structure on $\Der\AR$
given by the  commutator of derivations,
 the canonical map $\bm_*: (\dder\AR)\br\to \Der\AR,$
cf. Proposition \ref{form_smooth}, becomes a Lie algebra morphism.\qed
 \end{cor}

Next, following \cite{CB2}, for any $r\in R$ we define an associative 
$R$-algebra
$\Pi^r(B)$ to be the quotient of
${T^*B}$ by the two-sided ideal generated by the element $\Delta_B-r\in
T_1^*B+T_0^*B$.

For $r=0$, the grading on  ${T^*B}$ induces
a grading $\Pi^0(B)=\bigoplus_{k\geq 0}\Pi^0_k(B)$
on the  algebra $\Pi^0(B)$.
Therefore, the Euler derivation
$\Eu$ descends to a derivation 
$\Eu: \Pi^0(B)\to\Pi^0(B).$

For a general $r\in R$,  the grading on  ${T^*B}$ 
gives rise to an increasing filtration
$F_0\Pi^r(B)\sset F_1\Pi^r(B)\sset\ldots.$
Writing $\gr\Pi^r(B)$ for the corresponding associated 
graded algebra, one has a natural surjective
graded algebra map $\Pi^0(B)\onto\gr\Pi^r(B)$.

\begin{rem} It has been shown in \cite{CB2}, that

\npb{There is a canonical isomorphism $\Pi^0(B)\cong
T^\bullet_B(H^1_R(B, B\o B))$.}

\npb{If $B=\k[X]$ is the coordinate ring of 
a smooth affine curve $X$, then
we have $\Pi^0(\k[X])\cong \k[T^*X],$
the coordinate ring of the total space
of the cotangent bundle on $X$.}
\vskip 1pt

\npb{If $B=\k Q$ is the path algebra of a quiver $Q$,
then $\Pi^r(\k Q)\cong\Pi^r(Q)$ is the corresponding 
{\em deformed preprojective algebra}, as defined
in \cite{CBH}.}
\end{rem}

For $r=0$, the grading on the algebra
$ \Pi^0(B)$ induces
a grading $ \Pi^0(B)\br=\bigoplus_{k\geq 0}\Pi^0_k(B)\br.$ 
With respect to this grading, the Lie bracket
on $\Pi^0(B)\br$ has degree $(-1)$, that is,
we have $\{-,-\}: \Pi^0_k(B)\br\times \Pi^0_l(B)\br\to
\Pi^0_{k+l-1}(B)\br,$ for any $k,l\geq 0$.
Therefore, the map
$\Eu_{\Lie}: \Pi^0(B)\br\to\Pi^0(B)\br$, such that 
$\Eu_{\Lie}(a):=(k-1)\cdot a, \, \forall a\in \Pi^0_k(B)\br,\,
k=0,1,\ldots,$
is a derivation of the Lie algebra $ \Pi^0(B)\br$.
We let $\k\cdot\Eu_{\Lie}\ltimes\Pi^0(B)\br$
denote the corresponding semi-direct product Lie algebra.

The following result will play an important role in applications
to preprojective algebras associated with quivers.
\begin{prop}\label{stillinj} 
\vi The Lie bracket on   $({T^*B})\br$
descends to
a Lie algebra structure on $\Pi^r(B)\br$,
and there is
 a natural 
Lie algebra morphism $\Pi^r(B)\br\to\HH^1_R(\Pi^r(B)),\,
p\mapsto \b{\th}_p.$

\vii For $r=0$, the assignment $p\mapsto \b{\th}_p,\,
\Eu_{\Lie}\mapsto \Eu,$ gives rise to a 
Lie algebra morphism
$\k\cdot\Eu_{\Lie}\ltimes(\Pi^0(B)\br/R)\to\HH^1_R(\Pi^0(B)).$
This morphism  is  {\sf{injective}} 
provided one has $(\DR^0_RB)\closed=R$.
\end{prop}

Part (i) of the Proposition is immediate from
Proposition \ref{br_cor} applied to the algebra $A={T^*B}$.
Further, it  is straightforward to see that
the map of part (ii) is a  Lie algebra morphism.
The injectivity statement in part (ii)
is more complicated; it will be proved later, in Sect. \ref{pf_stillinj}.
 
Note that condition $(\DR^0_RB)\closed=R$ in part (ii) may be
interpreted as saying that the algebra $B$ is {\em connected}.

\subsection{The Liouville 1-form.}\label{La}
We recall the structure of  $\Om^1_R(T_BM),$ the bimodule of
noncommutative 1-forms on 
the  tensor algebra $A:=T_BM$ of 
an arbitrary  $B$-bimodule $M$.
To this end,
 introduce the following $A$-bimodule
\beq{DOM}
\dis\wt\Om:= (A\o_B \Om^1_RB\o_B A)\bigoplus 
(A\o_R  M\o_R  A).
\eeq
Abusing the notation slightly,  for any  $a',a''\in A,\,m\in M, \beta\in \Om^1_RB$,
we write $a'\cdot\pa  m\cdot a'':=0\oplus (a'\o m\o a'')\in\wt\Om$ 
and $a'\cdot \pa{\beta}\cdot a'':=(a'\o \beta\o a'')\oplus 0\in\wt\Om.$
 
Let $Q$ denote the  $A$-subbimodule in
$\wt\Om$ 
 generated
by the following set
\beq{dom}
\{
(\pa{b'mb''})
-\pa{db'} \cd (mb'')-b'\cd \pa  m\cd b''-
(b'm)\cd \pa{d  b''}\}_{\{b',b''\in B,m\in M\}}
\eeq

\begin{lem}\label{DBexact}
 Let $M$ be a  projective $B^e$-module, and $A:=T_BM$. Then

\vi There is an $A$-bimodule isomorphism $\Om_R^1A\cong\wt\Om/Q$.

\vii The imbedding of the first direct summand in $\wt\Om$, resp. the
projection to the second  direct summand in  $\wt\Om$, induces,
via the isomorphism in $\sf{(i)}$, a canonical extension of $A$-bimodules
$$ 0\map A\o_B \Om^1_RB\o_B A\stackrel{\eps}\map
\Om_R^1A\stackrel{\nu}\map A\o_B M\o_B A\map 0.
$$

\viii The assignment  
$B\oplus M=T^0_BM\oplus T^1_BM\to \wt\Om,\,
b\oplus m\mapsto d\pa  b+\pa  m$
extends uniquely to a derivation 
$\pa d: A=T^\bullet_BM\to \wt\Om/Q$; this derivation
corresponds, via the isomorphism in  $\sf{(i)}$, to the canonical
universal derivation
$d: A \to \Om_R^1A$.

\iv If $B$ is smooth and $M$ is  finitely generated 
(as a $B^e$-module), then the algebra $A=T_BM$ is also smooth.
\end{lem}
\begin{proof} 
First of all, using the relations \eqref{dom}
it is straightforward to verify that
 the assignment 
$b\oplus m\mto \pa{db}+\pa  m$
gives rise to a derivation
$\pa d : A\to \wt\Om/Q.$ 

Now, let
$E$ be an $A$-bimodule and let $\delta: A\to E$
be a derivation. Again, the relations \eqref{dom}
insure that the assignment 
$\widetilde  {db}+\pa  m\mto
\delta(b)+\delta(m)$ extends to a well-defined
$A$-bimodule map $\wt\Om/Q\to E$.
It follows that the derivation
$\pa d : A\to \wt\Om/Q$ enjoys  the universal property
for  the bimodule of 1-forms. Thus, we conclude that
$\Om^1_RA\cong \wt\Om/Q$. This proves (i) and (iii).

Part (ii) is clear from \eqref{dom}, see also \cite[Corollary
2.10]{CQ};
part (iv) is \cite[Proposition 5.3(3)]{CQ}. 
\end{proof}

Recall that, associated with any 1-form $\al\in\DR^1_RB$,
there is a $B$-bimodule map
$\bi(\al): \dder_RB\to B,\,\Th\ms \bi_\Th\al$.
This map extends uniquely
to an algebra homomorphism $\bi(\al): {T^*B}=T_B(\dder_RB)\map B$,
such that $\bi(\al)|_{T^*_0B}=\Id_B$
 and $\bi(\al)|_{T_1^*B}=\bi(\al)$.
Finally, the  algebra homomorphism $\bi(\al)$ induces, by functoriality,
a DG algebra homomorphism $\bi(\al)_*: \Om^\bullet_R({T^*B})\to\Om^\bullet_RB.$

\begin{prop}\label{liou} Assume that $B$ is smooth. Then,  there
is a canonical {\sf Liouville 1-form} $\la\in \DR^1_R({T^*B})$
such that the following holds:

\vi One has $i_\Eu\la=0,$ and  $L_\Eu\la=\la$.

\vii In $\D_1B=\dder_RB,$
we have $\bi_{\Delta_{T^*B}}\la=\Delta_B.$

\viii For any  1-form $\al\in\DR^1_RB$, in $\DR^1_RB$, we have
$\bi(\al)_*(\la)=\al$.
\end{prop}

\begin{rem} The  Proposition above
is a noncommutative analogue of a similar result  for the standard  
Liouville 1-form on $T^*X$, the cotangent
bundle on a manifold $X$.  To explain this, observe that
any 1-form $\al$ on $X$
is, by definition, a section $\al: X\to T^*X$ of the
 cotangent
bundle. Now, equation $\bi(\al)_*(\la)=\al$ of the Proposition 
is an analogue of the equation $\al^*(\la)=\al$, for the
pull-back of the Liouville 1-form on $T^*X$
via the section $\al: X\to T^*X$.
\erem

To construct the Liouville 1-form,
we consider the following natural 
maps 
$$(\dder_RB)\o_{B^e}\Om^1_RB=\Hom_{B^e}(\Om^1_RB,B\otimes
B)\o_{B^e}\Om^1_RB\to\Hom_{B^e}(\Om^1_RB,\Om^1_RB).
$$
Since $B$ is smooth, the last map is a bijection, cf.
 Proposition \ref{form_smooth}.
We invert this bijection and let $\La
\in(\dder_RB)\o_{B^e}\Om^1_RB$ 
denote the canonical
element corresponding to the
identity $\Id\in\Hom_{B^e}(\Om^1_RB,\Om^1_RB)$
under the composite (from right to left)
of the inverse maps in the displayed formula above.

We can write $\La=\sum_s \Ups_s\o \gamma_s, $ for some
$\Ups_s\in \dder_RB$ and $\gamma_s=db_s,$ where
$b_s\in B$ and $s=1,\ldots,l.$
The map $\Om^1_RB\to\Om^1_RB$ that corresponds
to the element $\sum_s \Ups_s\o {\gamma}_s$
is given by the formula
$\al\mto \sum_s (i'_{\Ups_s}\al)\cd {\gamma}_s\cd (i''_{\Ups_s}\al).$
By construction, this should be the identity map; thus
we must have
\beq{Omid}\sum\nolimits_s (i'_{\Ups_s}\al)\cd {\gamma}_s\cd (i''_{\Ups_s}\al)=\al,
\quad\forall \al\in \Om^1_RB.
\eeq

Now, by Lemma \ref{DBexact}(ii) applied to the
 $B$-bimodule $M=\dder_RB$, there is a canonical
${T^*B}$-bimodule imbedding
$\eps: {T^*B}\o_B\Om^1_RB\o_B{T^*B}\into \Om^1_R({T^*B})$.
Applying the commutator quotient functor
$E\mto E/[{T^*B},\,E]$ to the ${T^*B}$-bimodule map $\eps$
we get the following chain of natural maps
$$
\frac{{T^*B}\o_B\Om^1_RB}{[B,\,{T^*B}\o_B\Om^1_RB]}
=\frac{{T^*B}\o_B\Om^1_RB\o_B{T^*B}}{[{T^*B},\,
{T^*B}\o_B\Om^1_RB\o_B{T^*B}]}
\stackrel{\eps}\map
\frac{\Om^1_R({T^*B})}{[{T^*B},\,\Om^1_R({T^*B})]}.
$$
The leftmost term in this formula is nothing but
${T^*B}\o_{B^e}\Om^1_RB$ and the rightmost term 
equals $\DR^1_R({T^*B}).$
Thus, we have   canonical maps
$$
(\dder_RB)\o_{B^e}\Om^1_RB=T_1^*B\o_{B^e}\Om^1_RB\into
{{T^*B}\o_{B^e}\Om^1_RB\to\DR^1_R({T^*B}).}
$$

We define $\la\in\DR^1_R({T^*B})$ to be the image of the
canonical element $\La\in (\dder_RB)\o_{B^e}\Om^1_RB$ under the composite
map above.
Writing $\La=\sum_s \Ups_s\o {\gamma}_s$
and using the same notation as in formulas \eqref{DOM}-\eqref{dom},
we get the following analogue of the classical expression
`$\la=p\,dq$' for the Liouville 1-form:
\beq{gamma}
\la=\sum_{s=1}^l {\Ups_s}\cd \pa{\gamma}_s=
\sum_{s=1}^l {\Ups_s}\cd \pa{db}_s,\quad \Ups_s\in 
\dder_R B=T_1^*B,\;
\gamma_s=db_s,\; b_s\in B.
\eeq

\subsection{}
Given a $B$-bimodule $M$, let $M^\vee:=
\Hom_{B^e}(M,B\o B)$ denote the dual $B$-bimodule
equipped with the  $B$-bimodule structure
induced by the inner $B$-bimodule structure
on $B\o B$. We have a canonical  $B$-bimodule map
${\tt{bidual}}: M\to (M^\vee)^\vee$. This map  is an isomorphism
for any finitely generated projective $B^e$-module $M$.

In the special case $M=\Om^1_RB$, we
have $(\Om^1_RB)^\vee=\Hom_{B^e}(\Om^1_RB,B\o B)=
\dder_RB$. Therefore, if $\Om^1_RB$ is a finitely generated projective 
bimodule,
one also has a canonical  $B$-bimodule  isomorphism
${\tt{bidual}}:\Om^1_RB\iso(\dder_RB)^\vee,\,\al\mapsto
\al^\vee={\tt{bidual}}(\al)$.
Explicitly, for $\al\in\Om^1_RB,$ 
the element  $\al^\vee$ is a map
given by
\beq{alpha}
\al^\vee:\
 \dder_RB\to B\o B,\quad
\Th\mto -i''_\Th\al\o i'_\Th\al,
\eeq
where the flip of the order of  tensor factors
is due to the fact that the $B$-bimodule
structure on $\Hom_{B^e}(\dder_RB, B\o B)$
comes from the {\em inner} bimodule structure
on $B\o B$ and the negative sign comes from the
sign of the transposition permutation.

Further, we have  canonical bijections
\begin{align*}
{\small 
\dder_RB\underset{^{B^e}}\otimes\Om^1_RB
\underset{^\sim}{\stackrel{\Id\o{\tt{bidual}}}{\tooo}}
\dder_RB\underset{^{B^e}}\otimes(\dder_RB)^\vee\underset{^\sim}{\to}
\Hom_{B^e}(\dder_RB,\dder_RB).}
\end{align*}
The composite of these bijections sends the
canonical element $\La=\sum_s \Ups_s\o\gamma_s
\in (\dder_RB)\o_{B^e}\Om^1_RB$
to the identity map $\Id: \dder_RB\to\dder_RB.$
This yields the following identity, dual in some sense to
the identity \eqref{Omid}:
\beq{Did}
\sum\nolimits_s \Th''(b_s)\cd \Ups_s\cd \Th'(b_s)=
\sum\nolimits_s i''_{\Th}(\gamma_s)\cd \Ups_s\cd i'_{\Th}(\gamma_s)=\Th,
\quad\forall \Th\in\dder_RB.
\eeq

\begin{proof}[Proof of Proposition \ref{liou}.] To prove (i),
we compute
$i_\Eu\la=\sum_s i_\Eu(\Ups_s\cdot\pa\gamma_s)
=\sum_s\Ups_s\cdot (i_\Eu\pa\gamma_s)=0,$
since $i_\Eu\pa\beta=0$ for any $\beta\in\Om^1_RB.$
Observe further that
$\la$ belongs to the image of $\D_1B\otimes_B\Om^1_RB\o_B\D_0B$ in 
$\DR^1_R({T^*B})$, cf. \eqref{gamma}. We deduce that $L_\Eu\la=\la$.

To prove  property (ii) in the Proposition,  we apply
formula \eqref{Did} to $\Th=\Delta_B$ and observe 
that, for any $b\in B$, one has $i_{\Delta_{T^*B}}\pa{db}=
\Delta_{T^*B}(\pa b)= \pa{\Delta_Bb}.$
Thus, in $T^*B$, we  obtain a chain of equalities 
$$
\bi_{\Delta_{{T^*B}}}\la=\sum\nolimits_s (i''_{\Delta_{{T^*B}}}\pa\gamma_s)
\cd \Ups_s\cd
(i'_{\Delta_{{T^*B}}}\pa\gamma_s)
=\sum\nolimits_s \pa{\Delta''_B(b_s)}\cd \Ups_s\cd
\pa{\Delta'_B(b_s)}=\pa{\Delta_B}. 
$$

To prove (iii), fix $\al\in\DR^1_RB$. Then, we find
$$\bi(\al)_*(\la)=
\sum\nolimits_s \bi(\al)_*\bigl(\Ups_s\cd\pa\gamma_s\bigr)=
\sum\nolimits_s (\bi_{\Ups_s}\al)\cd\gamma_s.
$$
But $(\bi_{\Ups_s}\al)\cdot\gamma_s=(i''_{\Ups_s}\al)\cdot
(i'_{\Ups_s}\al)\cdot\gamma_s=
(i'_{\Ups_s}\al)\cdot\gamma_s\cdot (i''_{\Ups_s}\al)\,
\op{mod}\,[\Om^1_RB,B]$.
Hence, using the identity in \eqref{Omid} we get
$\bi(\al)_*(\la)=\al,$ as required.
\end{proof}

\subsection{}\label{omega}
We put $\om:=d\la\in \DR^2_R({T^*B}).$ This is a closed 2-form.

To avoid confusion, we write
$\omh:=\Delta_B$, viewed as an element
of $\D_1B$.

\begin{prop}\label{liuv2} \vi The 2-form $\om$ is  bi-symplectic.

\vii We have $\bi_{\Delta_{{T^*B}}}\om=d\omh$ and $\la=i_\Eu\om$.

\viii
For any  $p\in \D_kB,\,k=0,1,\ldots,$
 in $({T^*B})\br$, we have:
$$i_{\ad p}\la=p\cdot\omh,\quad\text{and}\quad
i_{\th_p}\la=k\cd p.
$$
\end{prop}

It is clear that  the Proposition above
implies Theorem \ref{DBthm}.

In preparation for the proof of Proposition \ref{liuv2}
we put  $A:={T^*B}$ and $M:=\dder_RB$.
We apply the functor $\Hom_{A^e}(-,A\o A)$
to the
short exact sequence
in  Lemma \ref{DBexact}(ii).
For  a smooth algebra $B$, one obtains this way
the following short exact sequence 
$$0\to A\o_B (\dder_RB)^\vee\o_B A
\stackrel{\nu^\vee}\map(\Om^1_RA)^\vee\stackrel{\eps^\vee}\map
A\o_B (\Om^1_RB)^\vee\o_B A\to 0.
$$
 Note that the middle term in this short  exact sequence 
 may be identified with
$(\Om^1_R({T^*B}))^\vee\cong \dder_R({T^*B}).$

Now, contraction with $\om$ gives a
${T^*B}$-bimodule map $\bi(\om): \dder_R(T^*B)\to\Om^1_R(T^*B),$
and we have

\begin{lem}\label{DBdiag} Let $B$ be  a smooth algebra and $A:=T^*B.$ Then,
the following diagram commutes
$$
\xymatrix{
0\to{A}\o_B (\dder_RB)^\vee\o_B {A}\ar[r]^<>(0.5){\nu^\vee}
\ar@{=}[d]_<>(0.5){-\Id\o{\tt{bidual}}\o\Id}&
\dder_R{A}\ar[r]^<>(0.5){\eps^\vee}\ar[d]_<>(0.5){\bi(\om)}&
{A}\o_B (\Om^1_RB)^\vee\o_B {A}\to 0\ar@{=}[d]_<>(0.5){\Id}\\
0\to {A}\o_B (\Om^1_RB)\o_B {A}\ar[r]^<>(0.5){\eps}&
\Om_R^1{A}\ar[r]^<>(0.5){\nu}&
{A}\o_B(\dder_RB)\o_B {A}\to 0
}
$$
\end{lem}
\begin{proof} We use the notation of Sect. \ref{La}. 
Thus, we have the canonical element
$\La=\sum_s \Ups_s\o \gamma_s$,
the Liouville 1-form $\la=\sum_s {\Ups_s}\cdot \pa{\gamma}_s $,
and the 2-form 
$\om=d\la=\sum_s d{\Ups_s}\cdot \pa{\gamma}_s ,$
where  $\gamma_s=db_s,$ see \eqref{gamma}.
Given  $\Th\in \dder_R(T^*B)$, we find
\begin{align}\label{phii}\bi_{\Th}\om&=
\sum\nolimits_s \bi_{\Th}(d{\Ups_s}\cdot \pa{\gamma}_s )\\
&=
-\sum\nolimits_s 
i''_{\Th}(d{\Ups_s})\cd \pa{\gamma}_s \cd
i'_{\Th}(d{\Ups_s}) +\sum\nolimits_s
i''_{\Th}(\pa{\gamma}_s )\cd d{\Ups_s}\cd
i'_{\Th}(\pa{\gamma}_s )\nonumber\\
&=
-\sum\nolimits_s 
\Th''({\Ups_s})\cd \pa{\gamma}_s \cd
\Th'({\Ups_s}) +\sum\nolimits_s
\Th''(\pa{b}_s )\cd d{\Ups_s}\cd
\Th'(\pa{b}_s ).\nonumber
\end{align}

Now, let $\al\in \Om_R^1B$ and put
 $$\overset{_\to}\al:=1\o{\tt{bidual}}(\al)
\o 1\in{T^*B}\o_B (\dder_RB)^\vee\o_B{T^*B}
\sset \dder_R(T^*B).$$
Thus,
$\overset{_\to}\al: {T^*B}\to {T^*B}\o{T^*B}$ is a derivation
 of degree $(-1)$ that annihilates the 
space $\D_0B=B$ and acts on $\D_1B=\dder_RB$ 
as the map $\al^\vee: \dder_RB\to
B\o B=\D_0B\o\D_0B\sset
{T^*B}\o{T^*B}$ given by formula \eqref{alpha}.
Therefore, for $\Th=\overset{_\to}\al$, the second sum in the last line
of formula \eqref{phii} vanishes, and we obtain
\begin{align*}
\bi_{\overset{_\to}\al}\om&=
-\sum\nolimits_s 
(\overset{_\to}\al)''({\Ups_s})\cd \pa{\gamma}_s \cd
(\overset{_\to}\al)'({\Ups_s})\\
&=-\sum\nolimits_s 
(\al^\vee)''({\Ups_s})\cd \pa{\gamma}_s \cd
(\al^\vee)'({\Ups_s})
=\sum\nolimits_s 
i'_{\Ups_s}(\al)
\cd\pa \gamma_s\cd i''_{\Ups_s}(\al)=\pa\al,
\end{align*}
where in the third equality we used  formula \eqref{alpha},
and 
the last equality is due to identity
\eqref{Omid}. This proves that the left square in the
diagram of the Lemma commutes.

Next, we prove the commutativity
of the  right  square in the
diagram of the Lemma. 
Since all maps in the diagram
are ${T^*B}$-bimodule maps, it suffices to
verify commutativity
on any derivation $\Th\in \dder_R(T^*B)$ such that 
$\eps^\vee(\Th)\in$
$ 1\o(\Om^1_RB)\o1$.
This condition means that we have
$\Th(\D_0B)\sset\D_0B\o\D_0B$. Thus, 
$\Th$ gives a derivation $\Th_B: B\to B\o B$,
and we have $\eps^\vee(\Th)=
1\o\Th_B\o1\in {T^*B}\o_B(\Om^1_RB)\o_B{T^*B}$.
We deduce that  one has
 $\Th(\pa b)=\pa{\Th_B(b)},$
for any $b\in B$.

We now take $\Th$ as above and
apply 
the map $\nu$ (in the bottom row of the diagram of the Lemma)
to the expression in the last line of \eqref{phii}.
The first sum in this expression is annihilated by the map $\nu$.
Also, we have $\nu(d{\Ups_s})=\Ups_s.$
Therefore, from \eqref{phii} we get
$$\nu(\bi_\Th\om)=\sum\nolimits_s
\nu\left(\Th''(\pa{b}_s )\cd d{\Ups_s}\cd
\Th'(\pa{b}_s )\right)
=\sum\nolimits_s
\pa{\Th''_B(b_s)}\cd \Ups_s\cd
\pa{\Th'_B(b_s)}=\pa{\Th_B},
$$
where the  last equality is due to identity
\eqref{Did}. This proves that the right square in the
diagram of the Lemma commutes.
\end{proof}

\proof[Proof of Proposition \ref{liuv2}.] To prove that 
$\om=d\la\in \DR^2({T^*B})$ is a bi-symplectic form,
we must show that the vertical map
$\bi(\om)$ in the middle of the diagram
of Lemma \ref{DBdiag} is a bijection. But
the vertical maps on the left and on the
right of the diagram
of Lemma \ref{DBdiag} are both  bijections.
Hence, the statement follows from the commutativity
of the diagram by diagram chase.

Next, by Proposition \ref{liou} we get
$d(\Delta_B)=d(\bi_{\Delta_{T^*B}}\la)=
-\bi_{\Delta_{T^*B}}(d\la)=-\bi_{\Delta_{T^*B}}\om.$
Also, we compute $i_\Eu\om=i_\Eu(d\la)=
L_\Eu\la-d(i_\Eu\la)=\la$, since $i_\Eu\la=0.$
Further,
by Corollary \ref{ipom}, we know that
$i_{\ad p}\om=p\,d\omh.
$ Thus we find
$$i_{\ad p}\la=-i_{\ad p}i_\Eu\om=
i_\Eu i_{\ad p}\om=i_\Eu(p\,d\omh)=p\cd
\Eu(\omh)=p\cd\omh,
$$
where in the last equality we have used that the element $\omh$
has degree 1 with respect to the grading on ${T^*B}$.
Similarly, for any $p\in\D_kB$, we have
 $\Eu(p)=k\cdot p.$ Hence, we obtain
$$i_{\th_p}\la=-i_{\th_p}\,i_\Eu\om=
i_\Eu\, i_{\th_p}\om=i_\Eu\,dp=
\Eu(p)=k\cd p.\qquad\Box
$$

\section{The representation functor}\label{rep_fun}
Throughout this section, $\k$ is an
algebraically closed field of characteristic zero.

\subsection{Representation schemes.}\label{rep_sch}
Let $I$ be a finite set and $R= \k{I}$, the
algebra of  functions $I\to\k$ with pointwise
multiplication.
For each $i\in I$,
let $e_i\in R$ be the characteristic function
of the one point set $\{i\}\sset I$. 
We have $1=\sum_{i\in I}\,e_i,$
furthermore, $\bbr:=\sum_{i\in I}\,e_i\otimes e_i$
is the symmetric separability element for $R$.

Giving a left $R$-module $V$ is the same thing as
giving 
 an $I$-graded vector space $V=\bigoplus_{i\in I}\, V_i,$
where
$V_i=e_i V$.
For any  left $R$-module $V$,
the action of $R$ gives  an  algebra map
$R\to\Hom_\k(V,V)$. This makes $\Hom_\k(V,V)$ an $R$-algebra,
to be denoted $\End:=\Hom_\k(V,V)$ below.

From now on, let $A$ be a {\em finitely presented}  associative
$R$-algebra.
Given a finite dimensional left
$R$-module $V=\bigoplus_{i\in I}\, V_i$, we may consider
the set $\Hom_{\alg{R}}(A,\End)$ of all $R$-algebra maps
$\rho: A\to \End$.
More precisely,  there is
an affine scheme  of finite type
over $\k$,
to be denoted $\Rep(A,V),$
such that the set  $\Hom_{\alg{R}}(A,\End)$
is the set of $\k$-points of  $\Rep(A,V)$.
The latter scheme, by definition, represents the functor
$B\mto \Hom_{\alg{R}}\bigr(A,\End_B(B\otimes V)\bigl),$
from the category of finitely generated commutative $\k$-algebras
to the category $\mathsf{Sets}$.

Let $\GL(V)^R=\prod_i \GL(V_i)$ be the group
of $R$-module automorphisms of $V$, and let
${\mathbb G}_m\sset \GL(V)^R$ be the 1-dimensional torus of scalar
automorphisms of $V$. We put $G:=\GL(V)^R/{\mathbb G}_m,$
and let $\fg=\Lie G$ be the Lie algebra of $G$.
We may (and will) identify both $\fg$
and its dual, $\fg^*$, with the codimension 1 subspace in $\oplus_i\,\End_\k V_i$ formed by 
linear maps with the vanishing total trace.

From now on, we will assume that the scheme
$\Rep (A,V)$ is {\em non-empty}, and
write $\k[\Rep (A,V)]$ for its
coordinate ring.
We consider $\k[\Rep (A,V)]\otimes\End$,
a tensor product of associative
algebras. The map $R\to \End$ makes
 $\k[\Rep (A,V)]\otimes\End$ an $R$-algebra
that
may also be identified with the algebra of regular
maps $\rep\to\End$, equipped with pointwise multiplication.
The action of $G$ on $V$ makes
$\rep$ a $G$-scheme. This gives a $G$-action on $\k[\Rep (A,V)]$ 
by algebra automorphisms.
We also have  a $G$-action on $\End$, by conjugation,
and we let $G$ act diagonally
on  $\k[\Rep (A,V)]\otimes\End$.
We write $\k[\Rep (A,V)]^G$, resp.,
$\bigl(\k[\Rep (A,V)]\otimes\End\bigr)^G$,
for the corresponding subalgebra of $G$-invariants.

\subsection{Evaluation homomorphisms.}\label{eval}
Write $\Sym E=\bigoplus_{j\geq 0}\,\Sym^j_\k E$ for the symmetric
algebra of a $\k$-vector space $E$. 

Given an $I$-tuple $\{d_i\in \k\}_{i\in I}$,
we define a linear function $\bd: R\to\k$ 
by the assignment $e_i\mto d_i$.
This  function
extends by multiplicativity to an algebra
map $\Sym R\to\k,$ and thus makes $\k$ a 1-dimensional
$\Sym R$-module to be denoted $\k_\bd.$

Now let $V=\bigoplus_{i\in I}V_i$
be a finite dimensional left $R$-module.
We put $d_i:=\dim V_i,$ and
write $\bd=\bd(V): e_i\mapsto d_i$ for the corresponding
linear function $\bd: R\to\k$.

Let $A$ be  an $R$-algebra.
The composite map $R\into A\onto A\br$
induces a graded algebra homomorphism
$\Sym R\to\Sym(A\br)$.
We introduce a commutative  $\k$-algebra
$\oo_\bd(A):=\k_\bd\bigotimes_{\Sym R}\Sym(A\br).$

To each element $a\in A$, one associates
the following {\em evaluation} function $\wh{a}: \rep\to\End,$
$\rho\mapsto \wh{a}(\rho):=\rho(a)$. The assignment
$a\mto\wh{a}$ clearly gives an associative $R$-algebra
homomorphism
$$ \ev:\ A\too\bigl(\k[\Rep (A,V)]\otimes\End\bigr)^G,
\quad a\mto\wh{a}.
$$

For any $a\in A$, composing the function
$\wh{a}$ with the trace map $\Tr:\End\to\k$,
applied to the second tensor factor above, we 
obtain a $G$-invariant element
$\Tr\wh{a}\in\k[\Rep (A,V)]^G$.
If $a\in[A,A]$, then
$\Tr\wh{a}=0$, due to symmetry of the trace.
Thus, the assignment
$a\mto\Tr\wh{a}$
gives a well-defined $\k$-linear map
$
\Tr\circ\ev:\, A\br\to\k[\Rep (A,V)]^G.$

It is clear that
 the element
$\wh e_i\in \k[\Rep (A,V)]\otimes\End$
corresponding to an  idempotent
$e_i\in R=\k{I}$
is a constant function on $\Rep(A,V)$ whose value
equals the projector
on the direct summand $V_i$.
Therefore, we have $\Tr\circ\ev(e_i)=\dim V_i=\bd(e_i)$.
Thus, we see that the linear map $\Tr\circ\ev:\,A\br\to\k[\Rep (A,V)]^G$
may be uniquely
extended, by multiplicativity,
to a commutative algebra morphism
\begin{align}\label{trace}
\psi_\bd:\ &\oo_\bd(A)\too\k[\Rep (A,V)]^G,\qquad \bd=\bd(V),\\
&a_1\&\ldots\& a_n\mto
(\Tr\wh{a}_1)\cdot\ldots\cdot (\Tr\wh{a}_n),\quad\forall a_1,\ldots,a_n\in A\br.\nonumber
\end{align}

There is also a natural evaluation map on differential forms
that sends the Karoubi-de Rham complex of $A$ to the ordinary
de Rham complex of the representation scheme. In more detail, write 
$\Om^\bullet(\Rep):=
\Om^\bullet(\rep)=\bigwedge^\bullet_{\k[\rep]}\,\Om^1(\rep)$ for the DG algebra
of algebraic differential forms on
the scheme $\rep$
(in the ordinary sense of commutative algebra).

We have an algebra homomorphism
\begin{align*}\ev:\,
\Om^nA &\map\bigl(\Om^n(\Rep)\otimes\End\bigr)^G,
\al=a_0\,da_1\ldots da_n\mapsto \wh{\al}=\wh{a}_0\,d\wh{a}_1\ldots d\wh{a}_n,
\end{align*}
defined as a composite
\begin{align*}
\Om^nA &=A\otimes (A/\k)^{\otimes n}\stackrel{\ev}\map
\bigl(\k[\Rep]\otimes\End\bigr)\otimes
\bigl(\Om^1(\Rep)\otimes\End\bigr)^{\otimes n}\\
&\to
\left(\bigwedge\nolimits_{\k[\Rep]}^n\,\Om^1(\Rep)\right)
\otimes\End^{\otimes n+1}\;\stackrel{\Id\otimes m}\too\;
\Om^n(\Rep)\otimes\End.
\end{align*}
Observe that since $\wh r$ is a constant function,
for any $r\in R$, we have $d\wh r=0$.
It follows easily that the map above induces
 a well-defined DG algebra morphism $\Om^\bullet\AR\to
\bigl(\Om^\bullet(\rep)\otimes\End\bigr)^G.$
Further, composing  the latter morphism
 with the trace map $\Id\otimes\Tr: \Om^\bullet(\rep)\otimes\End
\map\Om^\bullet(\rep),$ we thus obtain
a  linear map
\beq{trace_om2}
\Tr\circ\ev:\
\DR_R^\bullet A \too \Om^\bullet(\rep)^G,
\quad \al\mto
\Tr\wh{\al}.
\eeq 
The map \eqref{trace_om2} clearly commutes with the de Rham differentials.

\subsection{Double tangent bundle.}\label{KR}
Inspired by  an idea of Kontsevich and Rosenberg \cite{KR},
we are going to discuss the effect of the
action of the representation functor
on  double  derivations.

Fix a finite dimensional left $R$-module $V$.
To simplify notation, write ${\scr R}:=\k[\Rep(A,V)]$.
We identify the algebra ${\scr R}\otimes {\scr R}$
with $\k[\Rep(A,V)\times\Rep(A,V)].$
We will also 
consider the associative $R$-algebra ${\scr R}\otimes\End\otimes {\scr R}$.
There are two algebra maps
$\ev_l,\ev_r: A\to {\scr R}\otimes\End\otimes {\scr R},$
given by the formulas
$\ev_l(a):=\wh a\otimes 1,$ and $\ev_r(a):= 1\otimes\wh a,$
respectively. Here, $\wh a$ is viewed either
as an element
of ${\scr R}\otimes\End$ or as an element of
$\End\otimes {\scr R}$. The two maps $\ev_l,\ev_r,$ make
${\scr R}\otimes\End\otimes {\scr R}$ an $A$-bimodule,
via $a'(u\otimes F\otimes v)a'':=
\ev_l(a')\cdot(u\otimes F\otimes v)\cdot\ev_r(a'').$

We  define the double  tangent ${\scr R}$-bimodule
by
\beq{d_tang}
\calT^e(A,V):=\Der_R(A,{\scr R}\otimes\End\otimes {\scr R}).
\end{equation}
We observe 
that the  $A$-bimodule structure on ${\scr R}\otimes\End\otimes {\scr R}$
commutes with the obvious ${\scr R}$-bimodule structure,
hence, the latter  structure
indeed makes $\calT^e(A,V)$
 an ${\scr R}$-bimodule, that is, a
 $\k[\Rep(A,V)\times\Rep(A,V)]$-module
 (not necessarily free, in general).

Given $\rho\in\rep$, we
let $\k_\rho$  denote the 1-dimensional ${\scr R}$-module
such that $f\in \k[\Rep(A,V)]$ acts in $\k_\rho$ 
via  multiplication by the scalar $f(\rho)\in\k$.
It is clear that the
 geometric  fiber of the
 double  tangent ${\scr R}$-bimodule $\calT^e(A,V)$ at a point
$(\rho,\varphi)\in \Rep(A,V)\times\Rep(A,V)$
is given by
$$\calT^e(A,V)|_{(\rho,\varphi)}:=(\k_\rho\otimes\k_\varphi)
\bigotimes\nolimits_{{\scr R}\otimes
{\scr R}}\calT^e(A,V).
$$

To get a better understanding of the vector space on the right
of this formula,
observe that
any $R$-algebra homomorphism
$\rho: A \to \End V$ makes $V$ a left $A$-module,
to be denoted $V_\rho$. 
The action map $A\otimes{_R} V_\rho\to
V_\rho$ 
is surjective since $A$ has a unit. This gives a surjective map of 
left $A$-modules
$A \otimes{_R} V \onto V_\rho$,
where $ A \otimes{_R} V $ is regarded as a projective left $A$-module 
generated by the vector space $V$.
We set $K_\rho:=$ ${\Ker(A \otimes{_R} V \onto V_\rho),}$
a left $A$-module. 

\begin{lem}\label{double_T} For any
$(\rho,\varphi)\in \Rep(A,V)\times\Rep(A,V)$, one has a canonical
isomorphism
$$\calT^e(A,V)|_{(\rho,\varphi)}\;\simeq\;\Hom_A(K_\rho,V_\varphi).
$$

Furthermore, this  space is finite dimensional if $A$ is finitely
generated.
\end{lem}

\begin{proof}
For any left $A$-module $M$, we have
$\Hom_A(A \otimes{_R} V,\,M)=\Hom_R(V_\rho,M)$ and
$\Ext_A^1(A \otimes{_R} V,\,M)=0$. Hence,
the long exact sequence 
of Ext-groups arising from the short exact sequence
$ K_\rho\into A \otimes{_R} V \onto V_\rho$ reads:
$$\Hom_A(V_\rho,M)\into\Hom_R(V_\rho,M)
\to \Hom_A(K_\rho,M)
\to\Ext^1_A(V_\rho, M)\to 0.
$$

Observe next that
the  ${\scr R}$-bimodule structure  on ${\scr R}\otimes\End\otimes {\scr R}$
induces
an ${\scr R}$-bimodule structure
on each Hochschild cohomology
group ${H^p_R(A, \scr R\otimes\End\otimes \scr R),}$
$p\geq 0.$
Note that
$(\scr R\otimes\End\otimes \scr R)^R=\scr R\otimes\End^R\otimes \scr R,$
and $H^0(A,\,\scr R\otimes\End\otimes \scr R)=
(\scr R\otimes\End\otimes \scr R)^A$.
In particular, we have  the following exact sequence of 
${\scr R}$-bimodules (to be compared with the displayed exact sequence  above):
\begin{align}\label{der_double_HH}
({\scr R}\otimes \End\otimes {\scr R})^A
\into ({\scr R}\otimes\End\otimes {\scr R})^R
&\stackrel{\ad}\too \Der_R(A,{\scr R}\otimes\End\otimes {\scr R})\\
&\too H^1_R(A,{\scr R}\otimes\End\otimes {\scr R})\to 0.\nonumber
\end{align}

We leave to the reader to verify that,
for any two representations $\rho,\varphi\in \Rep(A,V)$,
 the geometric fibers at $(\rho,\varphi)$ of the 
${\scr R}$-bimodules occurring in \eqref{der_double_HH}  are given by
\begin{align*}
&({\scr R}\otimes\End\otimes
{\scr R})^A\big|_{(\rho,\varphi)}=\Hom_A(V_\rho,V_\varphi),\\
&\;({\scr R}\otimes\End\otimes
{\scr R})^R\big|_{(\rho,\varphi)}\enspace=\enspace\Hom_R(V_\rho,V_\varphi),\\
&\Der_R(A,{\scr R}\otimes\End\otimes {\scr R})\big|_{(\rho,\varphi)}=
\Hom_A(K_\rho,V_\varphi),\\
&H^1_R(A,{\scr R}\otimes\End\otimes
{\scr R})\big|_{(\rho,\varphi)}=\Ext_A^1(V_\rho,V_\varphi).
\end{align*}
The statement  of the Lemma is now clear.
\end{proof}

Next, we consider the following maps
$$
\xymatrix{
A\otimes A
\ar[rr]^<>(0.5){\ev\otimes\ev}&&
({\scr R}\otimes\End)\otimes(\End\otimes {\scr R})
\ar[rr]^<>(0.5){{\id_{\scr R}\otimes \bm\otimes\id_{\scr R}}}
&&
{\scr R}\otimes\End\otimes {\scr R},
}
$$
where $\bm: \End\otimes\End\to\End$ denotes the multiplication in the
$R$-algebra $ \End$.
We observe that the two maps above are morphisms of $A$-bimodules,
with respect to the outer bimodule structure
on $A\otimes A$. Let $\ev_{l,r}:
A\otimes A\too {\scr R}\otimes\End\otimes {\scr R}$ denote the composite
map,
which is an $A$-bimodule map again.

The 
``meaning'' of the double  derivation bimodule
$\dder\AR$ is somewhat clarified by the following observation:
{\em The $A$-bimodule morphism $\ev_{l,r}$ gives rise to a canonical
map}
$$\ev_{l,r}:\
\dder\AR\too\calT^e(A,V)=
\Der_R(A,{\scr R}\otimes\End\otimes {\scr R}).
$$

\subsection{Moment maps.}\label{moment_sec} Let $G$ be any algebraic group
with Lie algebra $\fg$. Write
$\fg^*$ for the linear dual of $\fg$;
we often view $x\in\fg=(\fg^*)^*$ as a linear function on $\fg^*$.
The group $G$ acts on $\fg^*$ via the coadjoint action.

Let $X$ be a $G$-scheme
and $\varpi\in\Om^2(X)^G$ be a closed
(not necessarily non-degenerate) $G$-invariant 2-form on $X$. Each  element $x\in\fg$
gives rise to a vector field, $\xi_x,$ on $X$.

\begin{defn} A $G$-equivariant map $\mu: X\to\fg^*$ is said
to be a moment map for $(X,\varpi,G)$ if the pull-back
morphism $\mu^*: \k[\fg^*]\map \k[X]$ has the following property:
$$d(\mu^*(x))=i_{\xi_x}\varpi,\quad\text{for any}\en x\in\fg.
$$
\end{defn}
\begin{rem}
Moment map is {\em not} uniquely determined by the triple
$(X,\varpi,G)$, in general.
\erem

Given a moment map $\mu: X\to\fg^*$, we write $\mu\inv(0)\sset X$
for its scheme-theoretic zero fiber, a $G$-stable subscheme in $X$.
If $X$ is  smooth and the 2-form $\om$ is nondegenerate,
then  $\k[X]$ acquires a natural
 Poisson algebra structure.
The latter  is known to descend to a well-defined
 Poisson  algebra structure on
$\k[\mu\inv(0)]^G.$

Now, let $A$ be a finitely generated associative $R$-algebra
equipped with 
a closed 2-form $\om\in\DR^2\AR$.
Fix a finite dimensional $R$-module $V$.
On $\rep$, we get a closed $G$-invariant
2-form $\Tr\wh{\om}\in(\Om^2\rep)^G$, cf. \eqref{trace_om2}.
We are interested in moment maps for the triple
$(\rep,\,\Tr\wh{\om},\,G)$. 

Observe that
the trace pairing $(x,y)\mapsto\Tr(x\cdot y)$ on $\fg$
 provides a $G$-invariant
isomorphism $\fg^*\cong\fg$ and that the Lie
algebra $\fg=\Hom_R(V,V)$ may be viewed as a subspace
in the associative algebra $\End=\Hom_\k(V,V)$.
 Thus, we have an imbedding $\fg^*\into\End$ and any moment map
$\mu$ for $(\rep,\Tr\wh{\om},G)$ may be identified with a map
$\rep\map\End$, that is, with an element
$\bmu\in\bigl(\k[\rep]\otimes\End\bigr)^G.$

Next, recall that we have defined a noncommutative moment map
$\munc: (\DR^2\AR)\closed\to (A/R)^R.$ Note that since $R=\k{I}$ is
commutative and
semisimple,
we have $R\sset A^R$ and, moreover, $(A/R)^R=A^R/R.$
Given a closed 2-form $\om\in\DR^2\AR$ and
 a representative  $\omh\in A^R$
of the class $\munc(\om)\in A/R$, 
we form the algebra $A_\omh=A/A\omh A$. The projection
$A\onto A_\omh$ gives a closed imbedding of schemes
$\Rep(A_\omh,V)\into\rep$, hence, a natural
restriction map
$\k[\rep]\onto\k[\Rep(A_\omh,V)].$

Recall that if the form $\om$ is bi-symplectic then
Propositions \ref{neck} and \ref{br_cor}
provide  the spaces $A\br$ and $(A_\omh)\br$ with
compatible Lie brackets. These Lie brackets
make the algebras $\Sym(A\br),\Sym((A_\omh)\br)$ into
Poisson algebras, with $\Sym R$ being a central
subalgebra in each of them. Therefore, for any linear map
$\bd:$
$ R\to\k$, the corresponding Poisson algebra structures
descend to the quotients $\oo_\bd(A)$ and $\oo_\bd(A_\omh),$
respectively.

Our main result about representation functors is the following
\begin{thm}\label{moment_thm} Fix a closed 2-form
$\om\in\DR^2\AR$.

\vi For any $\omh\in\aar$ in  the preimage of the
element $\munc(\om)\in\aar/R,$  the element 
$$\bmu:=\wh{\omh}\in (\k[\rep]\otimes\End)^G$$
gives a moment map for the triple $(\rep,\,\Tr\wh{\om},\,G)$.

Write $\mu: \rep\to\fg^*$ for the corresponding map,
and put $\bd=\bd(V)$. Then,
we have $\Rep(A_\omh,V)=\mu\inv(0)$, as subschemes in $\rep$;
thus, there is a  commutative diagram
$$
\xymatrix{\oo_\bd(A)\ar@{->>}[rr]^<>(0.5){\ev}\ar@{->>}[d]_<>(0.5){\text{\em
projection}}&&
\k[\rep]^G\ar@{->>}[d]^<>(0.5){\text{\em restriction to }\mu\inv(0)}\\
\oo_\bd(A_\omh)\ar@{->>}[rr]^<>(0.5){\ev}&&
\k[\mu\inv(0)]^G.
}
$$

\vii Assume in addition that $A$ is smooth and  $\om$ is a bi-symplectic form.

Then, the scheme $\rep$ is  smooth;
furthermore,
the 2-form $\Tr\wh{\om}$ makes $\rep$ a symplectic manifold,
and all maps in the diagram above are Poisson algebra morphisms.
\end{thm}
\begin{rem}
Different choices
of elements $\omh$ corresponding
to the same class $\munc(\om)\in\aar/ R$
give rise to {\em different}
moment maps, in general.
\erem

\begin{proof}[Proof of Theorem \ref{moment_thm}.] 
We set
 $\wh A:=\k[\rep]\otimes\End,$ and view it as 
an $R$-algebra. The map
$a\mapsto \wh a$ is an $R$-algebra morphism $A\to\wh A$.

Given any
$x\in\fg$ we may (and will) view it as the element $1\otimes x\in\wh A$.
In particular, for $u\in \wh A$, we write $x\cdot u:=(1\otimes x)\cdot u$,
resp., $u\cdot x:=u\cdot(1\otimes x).$ 
Thus, associated to $x\in\fg$, we have an inner derivation
$\ad x: \wh A\to\wh A,\,
u\mapsto\ad x(u):=(1\otimes x)\cdot u-u\cdot(1\otimes x).
$ 

Further, differentiating the $G$-action on $\rep$,
one associates to any element $x\in\fg$
a vector field $\xi_x$ on $\rep$, and we let
$\xi_x(f)$ denote the result of the $\xi_x$-action
on a  function $f\in\k[\rep]$. 
We will also consider
$\xi_x$-action on elements of
$\wh A$; 
for an $\End$-valued function
$f\otimes y\in \k[\rep]\otimes\End=\wh A,$
the action is defined 
by the formula $\xi_x(f\otimes y)=\xi_x(f)\otimes y.$
This way, one gets a derivation
$\xi_x: \wh A\to\wh A.$

Observe next that, for 
any $a\in A$, the function
$\wh{a}:\rep\to\End$ is $G$-equivariant. Hence, in $\wh A$,
one has an equality
\beq{xa}
\xi_x\wh{a}= \ad x (\wh{a}),\quad\forall x\in\fg,\,a\in A.
\eeq
This implies a similar equation for differential forms.
For instance, given a 2-form $\om=a_0\,da_1\,da_2\in\DR^2\AR$,
we get
\begin{align}\label{itr}
i_{\xi_x}\wh{\om}&=
i_{\xi_x}(\wh a_0\,d\wh a_1\,d\wh a_2)=
\wh a_0\xi_x(\wh a_1)\,\,d\wh a_2-\wh a_0\,\,d\wh a_1\,\,\xi_x(\wh a_2)
\\
&\stackrel{\eqref{xa}}{\eqq}
\wh a_0\ad x(\wh a_1)\,\,d\wh a_2-\wh a_0\,\,d\wh a_1\,\,\ad x(\wh a_2)=
i_{\ad x}(\wh a_0\,\,d\wh a_1\,\,d\wh a_2)=i_{\ad x}\wh\om.\nonumber
\end{align}

Now, with the notation introduced above,
proving the
theorem amounts to showing the identity:
\beq{moment_ident}
i_{\xi_x}(\Tr\wh{\om})=d\Tr(x\cd\wh{\omh}),
\quad\forall x\in \fg.
\eeq

We have
$$i_{\xi_x}(\Tr\wh{\om})=\Tr(i_{\xi_x}\wh{\om})=
\Tr(i_{\ad x}\wh\om),
$$
by \eqref{itr}.
Further, from Proposition \ref{omh} applied to the algebra
$\wh A$ and to the 2-form $\wh\om$ we deduce that,
in $\DR^1(\wh A)$, one has
$$
\Tr(i_{\ad x}\wh\om)=
\Tr\bigl(x\cd d\munc(\wh\om)\bigr)=
\Tr\bigl(x\cd \wh{d\omh}\bigr)=
\Tr\bigl(d(x\cd \wh{\omh})\bigr)=d\Tr(x\cd\wh{\omh}).
$$
 Hence, we obtain
$i_{\xi_x}\Tr\wh{\om}=\Tr(i_{\ad x}\wh\om)=d\Tr(x\cdot\wh{\omh}).$
This proves formula \eqref{moment_ident}. 
We conclude that $\wh{\omh}$ is indeed a moment map for
the triple $(\rep,\Tr\wh{\om},G).$

Finally, surjectivity of the horizontal maps in the
diagram of the Theorem follows from the well-known
result due to Le Bruyn-Procesi \cite{LBP}.
Commutativity of the diagram is clear.

We now turn to part (ii) of the Theorem. 
First, recall that if $B$ is any
$R$-algebra then, for any
left $B$-module $Q$, the space
$Q^*=\Hom_\k(Q,\k)$ acquires a natural
structure of right $B$-module,
 and vice versa.
Further, given a right $B$-module $P$
and a left  $B$-module $Q,$ 
one has a canonical linear map
\begin{align}\label{Psi}
&\Psi:\ P\otimes_B Q\too\Hom_B(P,Q^*)^*,\quad p\otimes q
\mto\Psi(p\otimes q),\\
&\qquad\text{where}\en
\Psi(p\otimes q):\ \Hom_B(P,Q^*)\map\k,
\en f\mapsto \langle f(p), q\rangle.\nonumber
\end{align}
It is immediate to see that the map $\Psi$
is a bijection if $P$ is a finite
rank free $B$-module and $\dim_\k Q<\infty$. We deduce
that:
{\em For any finitely generated projective 
$B$-module $P$ and finite
dimensional module $Q$, the map $\Psi$ becomes
an isomorphism of (finite dimensional)
vector spaces}.

We  now resume the proof of part (ii) of the Theorem.
The implication
$A$ {\em is smooth} $\Longrightarrow$
$\rep$ {\em is smooth} is standard, cf. e.g. \cite{Gi3}.

Next, fix $\rho\in\rep$, view $V=V_\rho$ as
a left $A$-module and view $\End_\k V_\rho$ as an $A$-bimodule.
We put
$\End_\rho:=\End_\k V_\rho$ and $\End_\rho^*:=(\End_\k V_\rho)^*$.
The trace paring $(u,v)\mto \Tr(u\cdot v)$ gives
an $A$-bimodule isomorphism
$tr: \End_\rho^*\iso \End_\rho$.

  By a standard infinitesimal
computation, the Zariski tangent space to
the scheme $\rep$ at the point $\rho\in\rep$ equals
$$T_\rho\rep=\Der_R(A, \End_\rho)=\Hom_{A^e}(\Om^1\AR,\End_\rho).$$
Therefore, for the cotangent space at $\rho\in\rep$, we get
$T^*_\rho\rep=\Der_R(A, \End_\rho)^*=\Hom_{A^e}(\Om^1\AR,\End_\rho)^*$.

Further, let  $\om\in\DR^2\AR$.
The 2-form $\om$ gives rise to the `contraction map'
$\bi(\om):\dder\AR\map\Om^1\AR,$ and
 the corresponding
2-form $\Tr\wh{\om}\in\Om^2(\rep)$, induces
a similar contraction  map $\wh i=i(\Tr\wh{\om}): T_\rho\rep\to T^*_\rho\rep$.

It is straightforward to verify that
the  following diagram commutes:
$$
\xymatrix{
\dder\!A\otimes_{A^e}\End_\rho^*
\ar[d]^<>(.5){\bi(\om)\otimes\Id}
\ar[r]^<>(.5){\eqref{m}}&
\Der(A,\End_\rho^*)\ar@{=}[r]^<>(.5){tr}&
\Der(A,\End_\rho)\ar@{=}[r]&T_\rho\ar[d]_<>(.5){\wh i}\\
\Om^1\otimes_{A^e}\End_\rho^*
\ar[r]^<>(.5){\Psi}&
\Hom_{A^e}(\Om^1,\End_\rho)^*\ar@{=}[r]&
\Der(A,\End_\rho)^*\ar@{=}[r]&T^*_\rho,
}
$$
where we have used simplified notation
$T_\rho=T_\rho\rep,$  $\Om^1=\Om^1\AR,$ and
$\Der=\Der_R$.

Assume now that $A$ is smooth, so
 $P:=\Om^1\AR$ is a finitely generated projective $A^e$-module. Then,
the map $\Psi$ in \eqref{Psi} becomes an isomorphism,
and the map \eqref{m} in the top row of the commutative diagram
above also becomes an isomorphism, by Proposition \ref{form_smooth}.
Further, if $\om$ is bi-nondegenerate, then the 
vertical map on the left of the diagram is a bijection.
It follows that  the 
vertical map on the right is a bijection as well.
Thus, we have proved that the 2-form $\Tr\wh\om$ is nondegenerate.

The last statement of the Theorem,  compatibility
with Poisson brackets, is the main result of \cite{Gi}.
\end{proof}

\section{Hamiltonian reduction in noncommutative geometry}
Throughout this section, $R$ stands for a finite dimensional semisimple
$\k$-algebra.
\subsection{Proof of Proposition \ref{br_cor}.}\label{proof}
Let $A$ be an $R$-algebra with a bi-symplectic form
$\om\in\DR^2\AR$, and
   $\bi(\om): \dder\AR\iso
\Om^1\AR$  the bijection induced by $\om$.
For any $v\in A$, let $\Th_v:=\bi(\om)\inv(dv)$ denote the
corresponding double derivation and, using Sweedler's
notation, write $\Th_v=\Th_v'\o \Th_v''.$

With this understood, we have the following explicit formula
for the map
$\SH_\om : \DR^1\AR\to\Der\AR,$ introduced in \eqref{SH}.

\begin{lem}\label{SH_formula}  
Let $u',\,u'',\,v\in A,$ and
put
$\th:=\SH_\om(u'\cdot dv\cdot u'')\in\Der\AR.$  Then,

\vi We have
$\dis\th(q)=\Th_v'(q)\cd u''\cd u'\cd \Th_v''(q),\quad\forall q\in A.$

\vii Let $\omh\in A$ be a representative of the
class $\munc(\om)\in (A/R)^R$, let $u=u''\o u'\in (A\o A)^R,$
and put  $\th:=\SH_\om(u'\cd d\omh\cd u'')$. Then,
$$\bi(\om)\inv(u'\cd d\omh\cd u'')=u'\cd \Delta\cd
  u''=\ad u,\quad\text{and}\quad
\th=\ad(u''u').$$
\end{lem}
\begin{proof} The map $\bi:=\bi(\om)$ is an $A$-bimodule isomorphism,
hence, so is $\bi\inv$. We deduce
$$\bi\inv(u'\cd  dv\cd  u'')=
u'\cd \bi\inv(dv)\cd  u''=
u'\cd \Th_v\cd  u''
=(\Th_v'\cd  u'')\o (u'\cd \Th_v'').
$$
Further,  by definition one has  $\th=\bm\br\ccirc \bi\inv(u'\cdot dv\cdot
u'')$. Hence, we get $\th(q)=$
${(\Th_v'(q)\cdot u'')\cdot} (u'\cdot\Th_v''(q)),$
and (i) follows.
 Part (ii) is  a reformulation of the first
equation in Proposition  \ref{omh}(iii).
\end{proof}

Next, 
we set $I:=A\omh A\sset A$, so $A_\omh=A/I.$

\begin{lem}\label{deck}
\vi For any $p\in A\br$,
one has $\th_p(I)\sset I.$
Hence, $\th_p$ descends to a well-defined derivation
$\bar{\th}_p: A/I\to A/I$.

\vii If $p\in I + [A,A],$ then $\bar{\th}_p: A/I\to A/I$ is an inner derivation.
\end{lem}

\begin{proof}
To prove (i), we 
compute
$\th_p(A\omh A)\sset \th_p(A)\,\omh\, A+ A\,\th_p(\omh)\, A +A\,\omh \,\th_p(A).$
The  summand in the middle vanishes
since $\th_p(\omh)=0$, by Lemma  \ref{H_1}(ii). The other two  summands clearly
belong to $I=A\omh A$, and (i) follows.

We turn to part (ii).
It suffices to prove the result for elements
 of the form $p=a  \omh b\in A\omh A=I.$
For such an element, in $\Om^1\AR$, we have
$\dis
dp=d(a  \omh b)=  (da  )\,\omh\, b  +a  \,(d\omh)\, b+a  \,\omh\, (db). 
$
Hence, using the formulas and notation of
Lemma \ref{SH_formula}, for $\th_p=\th$ and any $q\in A$, we find
$$
\th_p(q)=\Th'_a(q)\cd \omh\cdot b\cd\Th''_a(q)
+\ad(ba)(q)+\Th'_b\cd a  \cd\omh\cd\Th''_b(q).
$$
In this sum, the first and third summands
clearly  belong to the ideal
$A\omh A=I$. Hence, these terms give no contribution
to the induced  map $A/I\to A/I,\,q\mapsto\bth_p(q).$
The summand in the middle  gives an inner derivation
$\ad(ba)$.
We conclude
that the induced derivation $\bar{\th}_p$
is an inner derivation $A/I \to A/I$
corresponding to the image of $ba$ in $A/I$.
\end{proof}

\begin{proof}[Proof of Proposition \ref{br_cor}.]
By  part (i) of Lemma \ref{deck}
we get a map $\varphi: A\to \Der_RA_\omh,\,a\mto\bar{\th}_a.$
This map vanishes on $[A,A]$.
Furthermore, part (ii) of Lemma \ref{deck}
implies that the composite map 
$A\stackrel{\varphi}\too \Der_RA_\omh\to \HH^1_R(A_\omh)$ vanishes
on $I=A\omh A$. Thus, we conclude
that  this composite
 vanishes
on $I+[A,A]$, hence, descends to a well-defined
 map 
\beq{act}
\bar{\varphi}:\
A_\omh/[A_\omh,A_\omh]=A/(I+[A,A])\too \HH^1_R(A_\omh).
\eeq

We may now combine the map $\bar{\varphi}$
with the natural action of  $\HH^1_R(A_\omh)$
on $(A_\omh)\br,$
see Sect. \ref{ham_hoch},
to obtain a map
$$
(A_\omh)\br\map \Hom_\k\bigl((A_\omh)\br,\,(A_\omh)\br\bigr),\quad a\mto \bar{\varphi}_a.
$$
It is immediate from the construction that,
writing $\pi: A\br\onto (A_\omh)\br$ for the natural projection,
one has
$$
\bar{\varphi}_{\pi(a)}(\pi(a'))=\pi\bigl(\th_a(a')\bigr)=\pi\bigl(\{a,a'\}\bigr),
\quad \forall a,a'\in A\br.
$$
Therefore, skew-symmetry and Jacobi identity
for the Poisson bracket on $A\br$ imply
similar properties for the bracket
$(A_\omh)\br\times (A_\omh)\br\to (A_\omh)\br,$
$a\times a'\mapsto\bar{\varphi}_{a}(a').$

Finally, it is clear that  the map $\varphi: A\br\to \Der_RA_\omh$,
considered above
is a Lie algebra homomorphism. Therefore,
we deduce that
the map $\bar\varphi:(A_\omh)\br\map
\HH^1_R(A_\omh)$ in \eqref{act} is a Lie algebra
homomorphism, cf. also Theorem \ref{br_thm}(ii).
\end{proof}

\subsection{Main result.}\label{K}
Let 
$\DR^\bullet_{R }B$ be the 
 Karoubi-de Rham complex of an $R$-algebra $B$.
We define
$$K^\bullet_R B:=\Ker\big[\bi_\Delta :\DR_R^\bullet B\map\Om_R^{\bullet-1}B\big].
$$
The graded subspace $K^\bullet_R  B\sset\DR_R^\bullet B$ is stable under
 $d$, as well as under the maps 
$L_\th,i_\th$, for any $\th\in\Der_RB$, since
 $\bi_\Delta$ either commutes or anti-commutes with these maps.

For any inner derivation $\ad a,\,a\in B^R,$
by Corollary  \ref{ipom2} we have $i_{\ad a}\om=a\cdot \bi_\Delta\om$.
Hence, we deduce $i_{\ad a}\om=0,$ for any $\om\in K^\bullet_RB.$
In this case, we also have
$L_{\ad a}\om=d\ccirc 
i_{\ad a}\om+i_{\ad a}\ccirc d\om=0.$ We see that, when restricted to 
$K^\bullet_RB,$  the operations $i_\th, L_\th$ depend only on the cohomology
class of
$\th$ in $\HH^1_R(B)=\Der_RB/\Inn_RB$. Thus,
the complex $K^\bullet_RB$ acquires a natural structure
of $\HH^1_R(B)$-equivariant complex in the sense of Sect.
\ref{sec_GD}.

\begin{rem} According to \cite{Gi2},
there is a natural graded space isomorphism
$K^\bullet_RB\cong\HH_\bullet^RB$.
The above defined action of $\HH^1_R(B)$ on
$K^\bullet_RB$ then  becomes 
part of the standard action
of Hochschild cohomology on
Hochschild homology mentioned in Sect. 4.5.
This homological interpretation is not quite convenient, however, for 
performing Hamiltonian reduction.
Thus, in the present paper, we will neither use nor prove the isomorphism
$K^\bullet_RB\cong\HH_\bullet^RB$.
\erem

Note that since $\DR^{-1}_RB=0$, for any $\om\in\DR^0_RB=B\br,$ we have
$\bi_\Delta\om=0$. Thus, $ K^0_RB=B\br$, and the
 $\HH^1_R(B)$-action on $B\br$ is the action
mentioned in Sect. \ref{ham_hoch}, that has been already used in the
proof of Proposition \ref{br_cor}.
\smallskip

Now, let $A$ be an $R$-algebra,  $\om\in\DR^2\AR$ a closed
 $2$-form, 
 $\omh\in A$  a representative of the class $\munc(\om)\in A/R$. Below,
we will consider the complex $K^\bullet_RB$ for the algebra
$B=A_\omh:=A/A\omh A$, the Hamiltonian reduction of $A$ at $\omh.$
Let $\om_\omh\in\DR^2_RA_\omh$ denote the image of
$\om$ under the projection $\DR^2_RA\onto\DR^2_RA_\omh.$

\begin{lem}\label{bbii} For any $a\in A_\omh$ and  $\th\in\Der_RA_\omh,$
we have $da,\,i_{\th}\om_\omh\in K^1_RA_\omh.$ 

Furthermore,
the map $\th\mapsto i_{\th}\om_\omh$ kills inner derivations,
hence, descends to a well-defined map
$\bbi: \HH^1_R(A_\omh)\to K^1_R(A_\omh).$
\end{lem}
\begin{proof} It is immediate that $\bi_\Delta(da)=0$,
hence $da\in K^1_RA_\omh.$ 
It is also clear that we have  $\bi_\Delta\om_\omh=0$ since
$\omh=0$ in $A_\omh$. Therefore,
for any derivation $\th\in\Der_RA_\omh,$
we find $\bi_\Delta(i_{\th}\om_\omh)=-i_{\th}(\bi_\Delta\om_\omh)=0.$
Thus, we have proved that $i_{\th}\om_\omh\in K^1_RA_\omh.$
Further, for $a\in A_\omh$, by Proposition \ref{omh}(iii)
we get
$i_{\ad a}\om_\omh= a\cdot d\omh$, which is equal to zero in
$\DR^1_RA_\omh$.
\end{proof}

Assume now that $\om$ is a bi-symplectic 2-form.
Then we have the Lie bracket on $A\br$ provided
by Proposition \ref{neck} and we consider the projection
$A\br\onto(A_\omh)\br$.

The following Theorem, which is a strengthening of
 Proposition \ref{br_cor}, is our
 main result about Hamiltonian reduction in
noncommutative geometry.
\begin{thm}\label{br_thm} \vi The  bi-symplectic form $\om$
gives rise, canonically, to a Hamilton
operator $\SH_\omh:  K^1_R(A_\omh)\to \HH^1_R(A_\omh)$
such that $i_{\SH_\omh(\al)}\om=\al,$
for any $\al\in K^1_R(A_\omh)$, cf. \eqref{eqn}.

\vii With the  Lie bracket on $(A_\omh)\br$ induced by the Hamilton
operator $\SH_\omh$ from {\sf{(i)}}, the
natural projection
$A\br\onto(A_\omh)\br$ becomes a Lie algebra morphism.

\viii If the algebra $A$ is smooth then
the maps $\bbi$ and $\SH_\omh$ are mutually inverse
bijections.
\end{thm}

From part (i) we deduce
\begin{cor} The map $\SH_\omh$ is injective and the map $\bbi$ 
is surjective.\qed
\end{cor}

\subsection{Proof of Theorem \ref{br_thm}.}
First of all,  Lemma \ref{SH_formula}(ii)
yields the following result.

\begin{lem}\label{hh_map} Let $\om\in\DR^2\AR$ be
a bi-symplectic form on $A$. Then, the
bijection $\bi: \dder_RA\map \Om^1\AR$ induced by $\om$ 
restricts to an $A$-bimodule isomorphism
$\bi: \Inn_R(A,A\otimes A)\iso A\,(d\omh)\,A.$\qed
\end{lem}

Next, in the setup of Theorem \ref{br_thm}, we introduce
the two-sided ideal
$I:=A\omh A\sset A,$ and put
\beq{deri_def}
\deri:=\{\th\in\Der\AR\mid \th(I)\sset I\}=\{\th\in\Der\AR\mid \th(\omh)\in I\}.
\eeq

Further, consider the canonical surjection
$\varrho: \DR^\bullet_RA\onto \DR^\bullet_R(\AA)$,
cf. \eqref{A/I},
and put $\DR^\bullet_R(A,I):=\varrho\inv(\KK^\bullet_\RR(\AA)).$
Clearly, $\DR^\bullet_R(A,I)$ is a subcomplex of
the de Rham complex $\DR^\bullet_{R}A$.

\begin{lem}\label{Blem0} For 
the map 
$\SH_\om:
\DR^1\AR\to\Der\AR$
in \eqref{SH} and any $\al\in \DR^1\AR$, we have
$\al\in \DR^1_R(A,I)$ if and only if
$\SH_\om(\al)\in\deri.$

In particular, the map $\SH_\om$
 restricts to an 
imbedding 
$\DR^1_R(A,I)\into\deri$.
\end{lem}
\proof
The map  $\varrho$ used in the definition of
$\DR^1_{R}(A,I)$ is induced by  the
natural DG algebra projection $\varrho_\Om:
\Om^\bullet\AR\onto\Om^\bullet_R(A_\omh)=\Om^\bullet_R(A/I).$
Furthermore, for any $k=1,2,\ldots,$ we have a
commutative diagram
\beq{II_diag}
 \xymatrix{
\DR^{k}_RA\ar[r]^<>(0.5){\bi_\Delta}\ar[d]^<>(0.5){\varrho} &\Om^{k-1}\AR
\ar[d]^<>(0.5){\varrho_\Om}\\
\DR^{k}_R(A/I)\ar[r]^<>(0.5){\bi_\Delta}&\Om^{k-1}_R(A/I)
}
\eeq

For $k=1$, from the commutativity of the diagram, we
deduce 
\beq{keyI}
\DR^1_{R}(A,I)=
\varrho\inv\bigl(\KK^1_\RR(\AA)\bigr)=
\{\al\in\DR^1\AR\mid \bi_\Delta\al\in I\}.
\eeq

Now, let $\al\in\DR^1_RA$ and let
$\th=\SH_\om(\al)\in \Der\AR$ be the derivation corresponding to
$\al$ via the bi-symplectic form. 
Thus,
$\al=i_\th\om$, by \eqref{eqn}. Further, by
Proposition \ref{thomh}(i) we have $\th(\omh)=-\bi_\Delta\al$.
The proof of the Lemma is now completed by 
the following equivalences
$$ \al\in\DR^1_R(A,I) 
\en\Leftrightarrow\en \bi_\Delta\al\in I
\en\Leftrightarrow\en 
\th(\omh)\in I\en\Leftrightarrow\en\th\in\Der_R(A,I).\qquad\Box
$$
\smallskip

Any derivation
$\th\in\deri$ 
  induces  a well-defined
derivation $\bth:\AA\to\AA$.
It is clear  that $\deri$ is a Lie subalgebra in
$\Der\AR$ and
 the assignment $\th\mapsto\bth$
gives a well-defined Lie algebra morphism 
$\deri\to\Der_{\RR}\AA.$
We form the following composite
$$\Phi: \
\deri\map\Der_{\RR}\AA\onto\Der_{\RR}\AA/\Inn_{\RR}\AA=\HH^1_R(\AA),
$$
where the second map is the natural
projection. Also, let $\Psi$ be
 the restriction of
the map $\varrho,$ cf. \eqref{II_diag},
to the subspace $\DR^1_R(A,I)\sset \DR^1_{R}A$.

According to Lemma \ref{Blem0} we obtain a
diagram
\beq{key_diag} 
 \xymatrix{
\KK^1_\RR(\AA) &
\DR^1_R(A,I)\ar@{->>}[l]_<>(0.5){\Psi}
\ar@{^{(}->}[r]^<>(0.5){\SH_\om}&
\deri\ar[r]^<>(0.5){\Phi}&
\HH^1_\RR(\AA).
}
\eeq

\begin{lem}\label{Blem} For $\al\in \DR^1_R(A,I),$ we have 
$\Psi(\al)=0\en\Longleftrightarrow\en\Phi\ccirc\SH_\om(\al)=0$.
\end{lem}

\begin{proof} We need to introduce some auxiliary objects. Let 
$$\II:=\Om^\bullet\AR\cd(I+dI)\cd \Om^\bullet\AR
\sset\Om^\bullet\AR$$
 be the two-sided graded
 ideal in the algebra $\Om^\bullet\AR$ generated by the spaces $I\sset \Om^0\AR$
and $dI\sset \Om^1\AR$. Note that 
\beq{dI}
\II\cap \Om^0\AR=I,\quad\text{and}
\quad \II\cap \Om^1\AR=A(d\omh)A+ I\cd dA\cd A+ A\cd dA\cd I.
\eeq
Clearly, $\II$ is the smallest
$d$-stable  two-sided ideal in  $\Om^\bullet\AR$ 
that contains $I$.
Further, formula \eqref{A/I} shows that, in the notation
of diagram \eqref{II_diag}, we have
$\II=\Ker \varrho_\Om.$ 
Thus, from  the commutativity of 
the diagram, as in the proof of the previous lemma, we
deduce 
\beq{dr}
\DR^\bullet_{R}(A,I)=
\varrho\inv\bigl(\KK^\bullet_\RR(\AA)\bigr)=
\{\al\in\DR^\bullet\AR\mid \bi_\Delta\al\in\II\}.
\eeq

 Given a vector subspace $E\sset\Om^\bullet_{R}A$,
let $\lll E\rrr$ denote the image of
$E$ under the projection
$\Om^\bullet_{R}A\onto\DR^\bullet_{R}A$.
With this notation, we obtain
\beq{AW}
\KK^\bullet_\RR(\AA)\cong\DR^\bullet_R(A,I)/\lll\II\rrr=
\{\al\in\DR^\bullet\AR\mid \bi_\Delta\al\in\II\}/\lll\II\rrr.
\eeq
In particular,  we see that
 $\lll\II\rrr\sset\DR^\bullet_R(A,I)$,
in other words, one has
$\bi_\Delta(\lll\II\rrr)\sset\II.$

We can now turn to the proof of the implication `$\Rightarrow$'
of the Lemma.
The isomorphism in  \eqref{AW} and formula \eqref{dI} show that
\beq{ker}
\Ker\Psi=\lll\II\cap \Om^1\AR\rrr=\lll A (d\omh)A+I\cdot dA\cdot A +A\cdot
dA\cdot I\rrr=\lll I\cdot dA +A\cdot d\omh\rrr.
\eeq
Thus,
it suffices to  show that if
$\al\in \lll I\cdot dA +A\cdot d\omh\rrr$,
then $\Phi\ccirc \SH_\om(\al)=0$.

Assume first that $\al \in A\cd d\omh$. Then,
Lemma \ref{hh_map} implies that 
$\SH_\om(\al)$ is an inner derivation.
It follows that the induced map
$A/I\to A/I$  is also an inner derivation,
hence $\Phi\ccirc\SH_\om(\al)=0$.

Next let
$\al=u\, dv\in I\cdot dA$ and put $\th:=\SH_\om(\al)$.
The explicit formula for
$\th=\SH_\om(\al)$ provided by  Lemma \ref{SH_formula}(i)  says that,
 for any $a\in A$,
we have $\th(a)=\Th_v'(a)\cdot u\cdot \Th_v''(a).$
We see that  $u\in I$ implies $\th(a)\in I$,
for any $a\in A$. Hence, for $\al\in I\cdot dA,$ the corresponding
derivation
$\th$ induces the zero map $\bth: A/I\to A/I$,
and we are done.

We prove the implication `$\Leftarrow$'.
Let $\th\in \Der_R(A,I)$ be such that the induced
derivation $\bth: \AA\to\AA$ gives the zero class in
$\HH^1_\RR(\AA)$. This means that $\th=\xi+\ad a$,
where $a\in A^R$ and $\xi\in \Der_R(A,I)$ is such that
$\xi(A)\sset I$. (Note that for any $a\in A$,
we have $\ad a(I)=[a,I]\sset I,$ hence
$\Inn_RA\sset \Der_R(A,I)$.)

Let $\al=i_\th\om\in \DR^1\AR$ be the
1-form that corresponds to $\th$.
We can write
$\al=i_\xi\om+i_{\ad a}\om.$
By Lemma \ref{hh_map}, we get
 $i_{\ad a}\om\in \lll A(d\omh)A\rrr$.
Since  $\lll A(d\omh)A\rrr\sset\Ker\Psi$, by \eqref{ker},
we conclude that $i_{\ad a}\om\in\Ker\Psi$.

Next, given any 2-form $\be=x\,dy\,dz$,
we deduce from the inclusion $\xi(A)\sset I$ 
that in $\DR^1\AR$ one has
$i_\xi(x\,dy\,dz)=x\,\xi(y)\,dz-x\,dy\,\xi(z)\in
I\cdot dA$. Therefore, we obtain
$i_\xi\om\in \lll I\cdot dA\rrr\sset \Ker\Psi$.
Thus, we have $\al=i_\xi\om+i_{\ad a}\om\in\Ker\Psi$,
and the Lemma is proved.
\end{proof}

\begin{proof}[Proof of Theorem \ref{br_thm}.] Lemma 
\ref{Blem} implies that the map
$\SH_\om: \DR^1_R(A,I)\to\Der_R(A,I)$
sends the subspace $\Ker\Psi\sset \DR^1_R(A,I)$ into
the subspace $\Ker\Phi\sset \Der_R(A,I)$.
We see that this map 
descends to a well-defined and injective map,
see diagram \eqref{key_diag}: 
\beq{sshh}
\xymatrix{
\KK^1_\RR(\AA)=\DR^1_R(A,I)/\Ker\Psi\;
\ar@{^{(}->}[r]^<>(0.5){\SH_\om}&
\;\Der_R(A,I)/\Ker\Phi\cong\op{Im}\Phi.
}
\eeq
Composing the resulting map with the 
imbedding $\op{Im}\Phi\into\HH^1_\RR(\AA),$
we thus obtain an injective map $\SH_\omh: \KK^1_\RR(\AA)\into\HH^1_\RR(\AA)$.

The map $\SH_\omh$ satisfies all the identities
required for a Hamilton operator. To check these identities,
 one chooses representatives in $\DR^1_R(A,I)$
for classes in $\KK^1_R(\AA)$ and verifies the corresponding
identities for those  representatives. The latter
satisfy the required  identities due to the fact
that the map $\SH_\om: \DR^1\AR\to\Der\AR$ in
\eqref{SH} is known to be a  Hamilton operator.
This completes the proof of parts (i)-(ii) of the Theorem.

Assume now that the algebra $A$ is smooth.
Then, the map $\bm\br: (\dder\AR)\br\to\Der\AR$
is a bijection, by Proposition \ref{form_smooth}.
It follows that  the map $\SH_\om$ in
\eqref{SH}
is also a bijection. Now, Lemma \ref{Blem} shows
that  the latter bijection
induces a bijection $\SH_\om: \DR^1_R(A,I)\iso\deri.$
We see that, for the map $\SH_\omh:
 \KK^1_\RR(\AA)\into\HH^1_\RR(\AA)$ defined above, we have
$\op{Im}(\SH_\omh)=\op{Im}(\Phi)$.
Thus, we are reduced to proving surjectivity of
the map $\Phi$, which is equivalent
to surjectivity of the map
$\Der_R(A,I)\to \Der_\RR\AA,\,\th\mto \bth.$

To prove this last statement, we let 
$f: A\onto A/I=\AA$ denote the projection and consider the
following diagram
$$\xymatrix{
\Der_R(A/I,A/I)\ar[r]^<>(0.5){f^*}&
\Der_R(A,A/I)&\Der_R(A,A)\ar[l]_<>(0.5){f_*}
}
$$

For a smooth algebra $A$ the
 bimodule $\Om^1\AR$ is projective.
Hence, the functor
$\Der_R(A,-)=\Hom_{A^e}(\Om^1\AR,-)$ is exact.
We deduce that the map $f_*$ in the diagram above
is surjective. Hence, for any
$\bth\in \Der_R(A/I,A/I)$ there exists an $R$-linear derivation
$\th: A\to A$ such that
$f^*\bth=f_*\th$. From this equation we get
$f(\th(I))=(f^*\bth)(I)=0$.
Therefore,  $\th\in\Der_R(A,I)$
and, furthermore, the map
$\Der_R(A,I)\to \Der_\RR\AA$ sends $\th\mto \bth.$
This completes the proof of surjectivity,
hence, the proof of the Theorem.
\end{proof}

\subsection{Proof of Proposition \ref{stillinj}.}\label{pf_stillinj}
Set  $\Pi:=\Pi^0(B)$. This is a graded algebra with
degree zero component being equal to $B$. 

To prove the injectivity statement of part (ii) of the Proposition,
we use the homotopy invariance of 
de Rham cohomology, see Lemma \ref{homotopy}.
We deduce that $(\DR^0_R\Pi)\closed=(\DR^0_RB)\closed=R$,
where the last equality holds  by assumption.

Next, write $\om_\omh$ for the image of $\om$
in $\DR^2_R\Pi$.
Since ${T^*B}$ is a smooth algebra, Theorem
\ref{br_thm} implies that the map $\bbi: \HH^1_R(\Pi)\to
K_R^1(\Pi),\, \b\th\mapsto i_{\b\th}\om_\omh$ is a bijection.
The  bijection $\bbi$ clearly sends
the derivation $\b{\th}_p$ to the 1-form
$dp\in \DR^1(\Pi)$. The latter is zero if and only
if $p\in R,$ since $(\DR^0_R\Pi)\closed=R.$
Thus, we have proved that  $\b\th_p$ is a nonzero
element in $\HH^1_R(\Pi),$ for any
nonzero element $p\in \Pi\br/R$.

Further, by Proposition \ref{liuv2}(ii),
we have $\bbi(\Eu)=i_\Eu\om_\omh=\la_\omh$,
the image in $\DR^1_RA_\omh$ of the  Liouville 1-form $\la$. 
The 
form $\la_\omh$ is not exact. Indeed, we have $d\la_\omh=\om_\omh$,
and the latter form is nonzero since the
map $\b\th\mapsto i_{\b\th}\om_\omh$ is a bijection.

Thus, any nonzero element
of $\k\cdot\Eu_{\Lie}\ltimes(\Pi\br/R)$
goes to a class $u\in \HH^1_R(\Pi)$ such that
$\bbi(u)\neq 0.$ Hence, $u\neq 0,$ and injectivity
of the map in part (ii) of Proposition \ref{stillinj}
follows.
\qed

\section{The necklace Lie algebra}
\subsection{}\label{ncsymplectic}
Let $Q$ be a quiver with  vertex set $I$, 
and write $\overline{Q}$
for the double of $Q$. 
Let $\PP:=\k\overline{Q}$ be the path algebra of $\overline{Q}$.
The subalgebra $R\sset\PP$ formed by trivial paths
 may be (and will be) identified
with $\k{I}$.

Given an edge $a\in \overline{Q}$
we write $a^*$ for the corresponding reverse edge.
We set $\om:=\sum_{a\in Q}\, da\,da^*\in\DR^2_R\PP$.
It is clear that $d\om=0$.

The following result allows to apply
the general machinery developed in the previous sections
to algebras associated with quivers.

\begin{prop}\label{basic}  The path algebra $\PP$ is
smooth. Furthermore, we have 

\vi The assumptions of
Proposition \ref{ham_lemma} hold for $A=\PP$, more precisely,
one has
\beq{poinc}
H^k(\DR^\bullet_R\PP)=
\begin{cases}
R & \text{if}\enspace k=0\\
0 & \text{if}\enspace k>0.
\end{cases}
\eeq

\vii The 2-form $\om$ is  bi-symplectic.

\viii The following element gives a representative of the class $\munc(\om)\in\PP/R$:
$$
\omh:=\sum_{a\in Q}\, (a\cdot a^* - a^*\cdot a)
=\sum_{a\in Q}\,[a,a^*]\in [\PP,\PP]^R,
$$
more precisely, we have $\omh=\widetilde{\munc}(\om).$
\end{prop}

It is convenient to introduce
 a function $\epsilon: \overline{Q}\to \{\pm 1\}$
by setting  $\epsilon(a)$ to be 1 if $a\in Q$
and $-1$ if $a\in \overline{Q}\setminus Q$.
Then, we may write $
\omh=\sum_{a\in \overline{Q}}\, \epsilon(a)\cdot
 a a^*$.

\begin{proof}[Proof of Proposition \ref{basic}.]
For any two vertices  $i,j\in I$, 
let $E_{ij}$ be  a vector space with basis $\{a\}_{a: i\rightarrow j},$ formed
by the set of edges (from  $i$ to $j$)
of the quiver $\overline{Q}$. 
The 
assignment sending a pair of edges $a,b\in \overline{Q}$ to
$\epsilon(a)$ if $b=a^*$ and to zero otherwise
extends to a non-degenerate 
 pairing $E_{ij}\times E_{ji}\to \k$.
In particular,
the pairing provides, for any $i\in I$,
the space $E_{ii}$ with  a  non-degenerate skew-symmetric bilinear
form,
and also gives rise to vector
space isomorphisms $E_{ji}\cong E_{ij}^*$, for any
pair $i,j\in I,\,i\neq j$. We set $E:=\bigoplus_{i,j\in I}\,E_{ij}$
and extend the pairings above to a symplectic
form $\langle-,-\rangle: E\times E\to\k$. 

The space $E$ has a natural $R$-bimodule structure, and
the symplectic form  $\langle-,-\rangle$
yields an $R$-bimodule isomorphism
$S: E^*\iso E,\, \langle-,a\rangle\mto a^*$,
where $E^*=\Hom_\k(E,\k)$, the $\k$-linear dual
of $E$, is equipped with a natural 
 $R$-bimodule  structure.

We have $\PP=T_RE$, the tensor algebra of an $R$-bimodule.
This implies, in particular, that the algebra $\PP$ is smooth,
as claimed in Proposition \ref{basic}.

To prove Proposition \ref{basic}(ii),
we consider the following
 isomorphisms of $\PP$-bimodules
\begin{equation}\label{double1}
\Der_R(\PP, \PP\otimes\PP)\underset{^{F}}\iso
\PP\otimes_R E^*\otimes_R\PP
\,\underset{^{\Id\otimes S\otimes\Id}}\iso\,
\PP\otimes_R E\otimes_R\PP
\underset{^{G}}\iso\ncO^1_R(\PP),
\end{equation}
where the  maps $F,G$ are given, respectively, by
$$\Th
\stackrel{F}\mto \sum_{a\in \overline{Q}} (\Th''a \otimes \langle a,-\rangle\otimes \Th'a),
\quad\text{resp.},\quad
\sum_{a\in \overline{Q}}
f_a\otimes a\otimes g_a\stackrel{G}\mto\sum_{a\in \overline{Q}}
f_a\,da\,g_a.
$$

It is straightforward to verify that the composite
bijection in \eqref{double1} is nothing
but the map $\Th\mapsto\bi_\Th\om$.
We conclude that $\om$ is a
bi-symplectic 2-form on~$\PP.$ 

Since $\PP$ is smooth, we deduce
that the map $i: \Der_R\PP\to\DR^1_R\PP,\,\th\mapsto i_\th\om,$
is a bijection as well. This can also be seen directly
from the following diagram

\beq{derdr}{\xymatrix{
{\Der_R\PP}\quad
\ar[r]_<>(.5){F\br}^<>(.5){\sim}
\ar@/^2pc/[rrrr]|-{\,{\theta\mto i_\th\om}\,}
&\PP\otimes_{R^e} E^*\ar[rr]_<>(.5){\Id_\PP\otimes S}^<>(.5){\sim}
&&\PP\otimes_{R^e} E\ar[r]_<>(.5){G\br}^<>(.5){\sim}
&\quad \DR^1_R\PP,}}
\eeq
where the maps $F\br$ and $G\br$ are given by
$$
F\br:\;\th\mto \sum_{a\in \overline{Q}} \th(a)\otimes \langle a,-\rangle,
\quad
G\br:\;\sum_{a\in \overline{Q}}
f_a\otimes a\mto\sum_{a\in \overline{Q}}
f_a\,da
$$
In the formulas above, we
view $\PP$ as a {\em right} $R^e$-module, hence
we have $\PP\otimes_{R^e} E
=(\PP\otimes_R E)/[R,\PP\otimes_R E]$.
The inverse isomorphism in \eqref{derdr}
sends a 1-form $\sum f_a\,da\in \DR^1_R\PP$ to
the derivation $\th$ such that $\th(a)=f_{a^*}.$

The  path algebra 
comes equipped with a natural grading
$\PP=\bigoplus_{k\geq 0}\, \PP_k,$ by length of the
path. Using Lemma
\ref{homotopy} on homotopy invariance of
de Rham cohomology, we obtain formula \eqref{poinc}.
Observe further that  $\PP_0=R$ and $[\PP,\PP]\sset \bigoplus_{k\geq 1}\, \PP_k.$
 Therefore,
we have $[\PP,\PP]\cap R=0$, hence
the sequence \eqref{exact} in Proposition
 \ref{ham_lemma} is exact for $A=\PP.$
Thus, Proposition \ref{ham_lemma} holds for $A=\PP$.
This  completes the proof of part (i)
of Proposition \ref{basic}. 

Part  (iii) 
 is verified by
a straightforward computation based on
the last formula in Proposition \ref{omh}(iii).
\end{proof}

\subsection{Preprojective algebras.}
For any  $c\in \k I$,
one clearly has $(\omh-c)\,\op{mod}\,R=\munc(\om).$
Write $\PP(\omh-c)\PP$
for  the two-sided ideal in $\PP$
generated by $\omh-c$.
The corresponding Hamiltonian reduction
$\PP/\PP(\omh-c)\PP$
is known as {\em deformed preprojective
algebra},
$\Pi^c(Q):=\PP/\PP(\omh-c)\PP$.

\begin{rem} Let $B=\k{Q}$ be the path algebra of the quiver $Q$
(not of $\overline{Q}$). It is easy to see
 that there is a natural isomorphism $\PP\cong{T^*B},$ 
where ${T^*B}=\D(\k{Q})$ is the algebra considered in Sect. \ref{cot}.
Furthermore, the bi-symplectic form $\sum_{a\in Q} da\,da^*$, on $\PP,$
agrees with the  bi-symplectic form on $T^*(\k{Q})$
constructed in  Sect. \ref{cot}.
Note, however, that the grading on $\PP$ by length of the path 
{\em differs} from the grading on the algebra ${T^*B}$ used
 in Sect. \ref{cot} (this latter grading assigns grade degree
{\em zero} to any edge $a\in  Q$).

Further, 
it has been shown in \cite{CB2}
that, in the notation of  Sect. \ref{cot}
one has an algebra isomorphism
$\Pi^c(Q)\cong\Pi^c(\k{Q})$.
\erem

For any dimension vector $\mathbf d=(d_i\in \Z_{\geq 0},\, i\in I)$,
one associates with $\PP$, resp., with $\Pi^c(Q)$, the corresponding
commutative  algebra $\oo_\bd(\PP)$, resp., $\oo_\bd(\Pi^c(Q)),$
see Sect. 6.
The necklace Lie bracket induces Poison algebra structures 
on both $\oo_\bd(\PP)$ and  $\oo_\bd(\Pi^c(Q)),$
 see Sect. \ref{moment_sec}.
Furthermore, according to Theorem \ref{moment_sec}, the
scheme
  $\Rep (\PP,V)$, $\dim V=\bd,$  acquires a symplectic structure. This
makes $\k[\Rep (\PP,V)]$ a Poisson algebra, and
 each of the following subalgebras
is stable under the Poisson bracket
$$
 \k[\Rep (\PP,V)]\supset
\k[\Rep (\PP,V)]^G\onto\k[\Rep(\Pi^c(Q), V)]^G=\k[\mu\inv(c)]^{G}.
$$

In the special case $c=0$, 
we put $\Pi:=\PP/\PP\omh\PP$, the ordinary  preprojective
algebra. In Section \ref{proof_hh0} we will prove

\begin{prop}\label{hh0}
If $Q$ is connected and not Dynkin or extended
Dynkin, then its preprojective algebra $\Pi$ is
prime, and has center $Z(\Pi)=\k 1$.
\end{prop}

Next, recall the $R$-bimodule $E$ spanned by the edges of $\overline{Q}$.
The degree 2 component of 
the algebra $\PP$ is equal to $E\otimes_R E$.
Thus, using Sweedler's notation
we may write $\omh=\omh'\otimes \omh''\in E\otimes_R E.$
With this notation, one has the following explicit
description of the space of 1-forms, resp., of derivations, for the
algebra $\Pi.$
\begin{lem}
\vi The assignment 
$p'\,(da)\,p''\mapsto p'\otimes a\otimes p''$ gives a  $\Pi$-bimodule isomorphism
$$\Om^1_R\Pi\cong\frac{\Pi\otimes_R E
\otimes_R \Pi}{\Pi\cd(\omh'\otimes \omh''\otimes 1-
1\otimes \omh'\otimes \omh'')\cd\Pi}
$$

\noindent
\vii The assignment 
$\th\mto \sum_a\,\th(a^*)\otimes a$
induces a vector space  isomorphism
$$\Der_R\Pi\iso 
\big\{\sum\nolimits_{a\in \overline{Q}}\, p_a\otimes a\in
(\Pi\otimes_R E)^R\enspace\big|\enspace
 \sum\nolimits_{a\in \overline{Q}}\,[p_a,a]=0\big\}.
$$
\end{lem}

Proof of the Lemma is straightforward and is left to the reader.

\subsection{Necklace bracket.} 
The bi-symplectic form $\om$ gives a Lie bracket
on $\PP\br$, known as  the {\em necklace Lie bracket},
see \cite{BLB},\cite{Gi}.

To write explicit formulas for this bracket it is convenient to
introduce the following notation.
Given  a path $p$ in $\overline{Q}$ of length $\ell(p)$,
for each
 integer $1\le i\le \ell(p)$, let $p_{<i}$, $a_i$
and $p_{>i}$ be the paths of lengths $i-1$, $1$ and $\ell(p)-i$
respectively, with $p=(p_{<i}) a_i (p_{>i})$.
For each edge $a \in\overline{Q}$
we introduce a map

\beq{explicit}
\partial_a: \PP\br\map \PP,
\enspace
p\mapsto \partial_ap:= \sum_{i=1}^{\ell(p)} 
\delta_{a,a_i}\cdot  (p_{>i})(p_{<i})
\end{equation}
  (where  $\delta$ is the Kronecker delta function).

The proof of the next  lemma  will be given at the end of this section.

\begin{lem}\label{part2} 
We have the following formulas:
\begin{align*}
&dp=\sum\nolimits_{a\in  \overline{Q}}
(\partial_a p)\,da\quad\op{holds}\en\op{in}\en\DR^1_R\PP,\quad
\forall p\in\PP\br.\\
&\th_p(a)=\epsilon(a)\cd\partial_{a^*} p\quad\op{holds}\en\op{in}\en\PP,
\quad\forall p\in \PP\br;\\
&\{p,q\}=\sum\nolimits_{a\in  \overline{Q}}\epsilon(a^*)\cdot 
(\partial_a p)\,(\partial_{a^*}q),\quad
\forall p,q\in\PP\br.
\end{align*}
\end{lem}

By  Proposition \ref{br_cor}(ii),
we know that, for any $c\in R$,  the necklace Lie bracket
on $\PP\br$
descends to
a well-defined  Lie bracket
 on $\Pi^c(Q)\br$,
cf. also 
 Proposition \ref{stillinj}.

From now on, we will only consider the case $c=0$,
so $\Pi=\Pi^0(Q)$. We set
 $L:=\Pi\br=\Pi/[\Pi,\Pi].$

The grading on the path algebra $\PP=\bigoplus_{k\geq 0}\, \PP_k,$ by length of the
path, descends to  natural gradings
$\Pi=\bigoplus_{k\geq 0}\, \Pi_k$ and  $L=\bigoplus_{k\geq 0}\, L_k$.
These gradings make $\PP\br$ and $L$  graded Lie algebras,
with necklace Lie bracket {\em of degree} $(-2)$, e.g., for any $k,m\geq 0$, we have
$\{L_k,L_m\}\sset L_{k+m-2}$ (the bracket vanishes whenever $k+l <0$).

Associated with the grading on $\Pi$,
we have the corresponding {\em Euler derivation} $\eu: \Pi\to\Pi$
and a Lie algebra derivation
$\eu_{_{\Lie}}: L\to L$ defined by  $\eu|_{L_k}=(k-2)\cdot\Id_{L_k},$
for any $k=0,1,\ldots.$ Our present
notation for Euler derivations is different from the
one used in  Sect. \ref{cot}, since the
grading on $\Pi$ that we are using now
is {\em not} the one used in  Sect. \ref{cot}.

\subsection{Hochschild cohomology.}\label{HH}
We use the notation $\HH^\bullet(\Pi):=H^\bullet_R(\Pi,\Pi)$
for relative Hochschild cohomology of the preprojective
algebra $\Pi$, cf. Remark \ref{H_R}. According to Proposition \ref{br_cor},
we have a Lie algebra map $L=\Pi\br\to\HH^1(\Pi).$
As in Proposition \ref{stillinj},
this map gives rise to  a Lie algebra map $\k\cdot\eu_{_{\Lie}}\ltimes (L/R)\to\HH^1(\Pi),$
where the semidirect
product $\k\cdot\eu_{_{\Lie}}\ltimes (L/R)$ is
viewed as a graded Lie algebra with
$\eu_{_{\Lie}}$ assigned grade degree zero.

One of the  main results of this paper,
to be proved in \S9, is the following
\begin{thm}\label{main} Assume that $Q$ is neither Dynkin
nor extended Dynkin. Then 

\vi The map
$p\mapsto \bar{\th}_p,\,\eu_{_{\Lie}}\mapsto \eu,$ 
induces a graded  Lie algebra isomorphism
$$\k\cdot\eu_{_{\Lie}}\ltimes (L/R)\iso \Der_R\Pi/\Inn_R\Pi
=\HH^1(\Pi).$$

\vii The Hochschild cohomology of 
the preprojective algebra is given by
\beq{formula}
\HH^i(\Pi)\cong\begin{cases}\k &\op{if}\enspace i=0\\
\k\cdot\eu_{_{\Lie}}\ltimes (L/R)&\op{if}\enspace i=1\\
L&\op{if}\enspace i=2\\
0&\op{if}\enspace i>2.
\end{cases}
\end{equation}
\end{thm}

The isomorphism $\HH^1(\Pi)\cong
\k\cdot\eu_{_{\Lie}}\ltimes (L/R)$ in \eqref{formula} 
can be shown to map the (Gerstenhaber)
Lie bracket on $\HH^1(\Pi)$ to the bracket
in the Lie algebra $\k\cdot\eu_{_{\Lie}}\ltimes (L/R)$.
Furthermore, one verifies that
the  Gerstenhaber 
bracket $\HH^1(\Pi)\times \HH^2(\Pi)\map\HH^2(\Pi)$
corresponds, via  formula \eqref{formula} for $i=1,2$,
to the adjoint action of the Lie algebra
$L/R$ on $L$.
Thus, Theorem \ref{main}
completely describes the Gerstenhaber 
bracket on the Hochschild cohomology of $\Pi$.

 There is also an associative 
and graded-commutative {\em cup-product}
$\HH^i\times \HH^j\to\HH^{i+j}$
on  Hochschild cohomology.
In the case of the preprojective algebra,
the only non-trivial cup-product
is the skew-symmetric pairing 
$\cup: \HH^1(\Pi)\times\HH^1(\Pi)\to\HH^2(\Pi).$
It is easy to see that, under the isomorphism
in \eqref{formula}, this pairing is given by
the formulas
$$\eu          \cup\eu          =0,
\quad \eu\cup p =(\deg p)\cd p,\quad
p\cup q=\{p,q\},
$$
for any homogeneous elements $p,q\in L/R$,
(since $R$ is central
in $L$, the necklace bracket in the last formula is
viewed as a map $L/R\times L/R\to L$).

\subsection{Dynkin case.}\label{dyn}
 In the case of (finite) Dynkin  quiver $Q$,
the corresponding preprojective algebra $\Pi$ is
{\em not} Koszul, in general, and there are nonzero
cohomology groups $\HH^i(\Pi)$ for $i>2$,
see \cite{ES}, \cite{BBK}.

Now let $Q$ be an extended Dynkin quiver and
let $\G\sset SL_2(\k)$ be the finite
group associated to $Q$ via the McKay correspondence.
It is convenient to introduce a 2-dimensional
vector space $V$ and write
$\G\sset SL(V)=SL_2(\k)$.
Also, put ${\mathscr S}:=\G\sminus\{1\}$ and
let $\k[{\mathscr S}]^\G$ denote the space of class-functions on
${\mathscr S}$.

According to \cite{CBH}, the
preprojective algebra $\Pi=\Pi^0(Q)$
is Morita equivalent to  $\k[V]\sharp\G $, the
cross product of $\G$ with the polynomial algebra on $V$.
Thus, by Morita invariance of  Hochschild cohomology,
we have $\HH^\bullet(\Pi)\cong\HH^\bullet(\k[V]\sharp\G ).$
Further,
the  algebra $\k[V]$, hence the cross product algebra  $\k[V]\sharp\G $,
is Koszul. Therefore, for  Hochschild cohomology
 one easily finds
\beq{cbh}
\HH^i(\Pi)\cong\HH^i(\k[V]\sharp\G )=
\begin{cases}\k[V]^\G &\op{if}\enspace i=0\\
(\k[V]\otimes V)^\G&\op{if}\enspace i=1\\
(\k[{\mathscr S}]\oplus\k[V])^\G\otimes \wedge^2V
&\op{if}\enspace i=2\\
0&\op{if}\enspace i>2.
\end{cases}
\eeq

In formula \eqref{cbh}, the space $(\k[V]\otimes V)^\G$ stands
for the Lie algebra of $\G$-invariant polynomial
vector fields on $V$. Note that
we may interpret $L/R=\Pi\br/R=(\k[V]\sharp\G)\br/R$
as the Lie algebra of $\G$-invariant
{\em Hamiltonian} vector fields on $V$.
Thus, we see that the semidirect product
$\k\cdot\eu\ltimes (L/R)$,
occurring in formula \eqref{formula},
is a {\em proper} Lie subalgebra in
 $(\k[V]\otimes V)^\G$.
Observe also that the group  $\G$ acts trivially on $\wedge^2V$.
Thus,
$\k[V]^\G\otimes \wedge^2V=(\k[V]\otimes \wedge^2V)^\G,$
is the space  of $\G$-invariant polynomial
bivector fields on $V$. This  is clearly a rank 1 free
$\k[V]^\G$-module with generator
$\pi\in\wedge^2V$, the constant volume
element.

Further, by the general duality between Hochschild homology and
cohomology due to
Van den Bergh \cite{VB1}, see \eqref{vdb_iso}-\eqref{unit},
one has a canonical isomorphism
$\wedge^2V\otimes\HH_\bullet(\k[V]\sharp\G)\cong\HH^{2-\bullet}(\k[V]\sharp\G).$
Therefore, the isomorphisms in \eqref{cbh} yield formulas
for  Hochschild homology of the algebra $\Pi$ as well. In particular,
we get a chain of isomorphisms
\beq{CCC}
(\k[{\mathscr S}]\oplus\k[V])^\G\iso\wedge^2V^*\otimes \HH^2(\Pi)\iso
\HH_0(\Pi)=\Pi\br.
\eeq

The composite isomorphism in \eqref{CCC} can be described explicitly as
follows.
The restriction of  this map  to the direct summand
$\k[{\mathscr S}]^\G$ is obtained as a composition
 $\k[{\mathscr S}]^\G\iso
R_0\into R\into \Pi\br$, where we have used the notation
 $R_0:=\{r=\sum r_ie_i\in R\mid
\sum r_i=0\}$ and  the first of the above maps is 
 explained e.g. in  \cite{CBH}.
The  restriction of  the composite in \eqref{CCC} to the
 direct summand $\k[V]^\G$ is equal to the composition
$$\k[V]^\G\into \k[V]\sharp\G\onto
(\k[V]\sharp\G)\br\iso \Pi\br,
$$
where the first two maps are the natural maps
and the rightmost isomorphism is provided by Morita equivalence.

\subsection{Poisson algebras.}
Let $R[\omh]\sset\PP$ be the subalgebra
generated by $R$ and $\omh$, and let
 $\overline{R[\omh]}\sset \PP\br$
denote the image of $R[\omh]$
 under the projection
 $\eta: R[\omh]\into\PP\onto\PP/[\PP,\PP]=\PP\br$.
The imbedding of vector spaces
$\overline{R[\omh]}\into \PP\br$ extends, by multiplicativity
to an injective graded algebra homomorphism
$\Sym(\overline{R[\omh]})\into\Sym(\PP\br),$
of the corresponding symmetric algebras (over $\k$).

The Lie algebra structure on $\PP\br$, resp., on
 $\Pi\br$, gives the symmetric algebra
$\Sym(\PP\br),$ resp., $\Sym(\Pi\br)$,  the structure of a
Poisson algebra.
By Lemma \ref{hhh}, the space $\overline{R[\omh]}$
is contained
in the center of the Lie algebra $\PP\br$,
hence, 
$\Sym(\overline{R[\omh]})$ is a central subalgebra of
the Poisson algebra
$\Sym(\PP\br)$.

Similarly, the imbedding $R\into \Pi\br$
 makes $\Sym R$ a central Poisson subalgebra in
$\Sym(\Pi\br)$.

The third important result of this paper is
 the following theorem, to be proved
 in Section \ref{pc} below.

\begin{thm}\label{cen}
Assume that $Q$ is neither Dynkin nor extended Dynkin. Then
  
\vi The Poisson center of the algebra $\Sym(\PP\br)$ is equal to 
$\Sym(\overline{R[\omh]})$.

Furthermore, the
kernel of the natural projection $\Sym(R[\omh])\onto \Sym(\overline{R[\omh]})$ 
is a principal ideal in the algebra
$\Sym(R[\omh])$ generated by the element $\omh\in \Sym^1(R[\omh])=R[\omh].$

\vii  The Poisson center of  the algebra $\Sym(\Pi\br)$ is equal to $\Sym R$. 
\end{thm}

\begin{cor}\label{cent_prop} If $Q$ is neither Dynkin
nor extended Dynkin, then we have

\vi The center of the Lie algebra $\PP\br$ 
equals $\overline{R[\omh]}$; moreover, the kernel
of the projection $\eta: R[\omh]\onto \overline{R[\omh]}$ is
the line spanned by $\omh$.

\vii  The center of the Lie algebra $\Pi\br$ equals $R$.\qed
\end{cor}

\subsection{Proof of Lemma \ref{part2}.} First, 
fix an arbitrary derivation $\th$.
We compute
\begin{align}\label{ip}
i_\th\om=i_\th(\sum\nolimits_{a\in Q}\, da\,da^*)&=\sum\nolimits_{a\in Q}\,(i_\th da)\,da^*
-\sum\nolimits_{a\in Q}\, da\,(i_\th da^*)\nonumber\\
&=\sum\nolimits_{a\in Q}\,(\th(a)\,da^*-da\,\th(a^*))\\
&=\sum\nolimits_{a\in \overline{Q}}\,\epsilon(a)\cdot \th(a)\,da^*
\;\modu[\PP,\ncO^1_R\PP].\nonumber
\end{align}
 
Next, let  $p=(p_{<i}) a_i (p_{>i})$ be a path in $\overline{Q}$.
For the
1-form $dp\in\DR^1_R\PP=\ncO^1_R\PP/[\PP,\ncO^1_R\PP]$, we 
find 
\beq{dp}
dp=\sum_{i=1}^{\ell(p)}(p_{<i})\,da_i\,(p_{>i})=
\sum_{i=1}^{\ell(p)} 
 (p_{>i})(p_{<i})\,da_i\;\modu[\PP,\ncO^1_R\PP].
\end{equation}

Now,
 let  $\th : \PP\to\PP$ be the  derivation that
annihilates $\k I\sset \PP$
and acts on generators of $\PP$ by the 
formula $a\mto \epsilon(a)\cdot\partial_{a^*}p.$
For this derivation, we obtain
\begin{align*}
i_{\th}\om\;\stackrel{\eqref{ip}}{=\!=}\;
\sum_{a\in \overline{Q}}\sum_{i=1}^{\ell(p)} 
\epsilon(a^*)\cd\delta_{a,a_i}\cdot  (p_{>i})(p_{<i})\,da^*
=
\sum_{i=1}^{\ell(p)} 
 (p_{>i})(p_{<i})\,da_i\;\stackrel{\eqref{dp}}{=\!=}\;dp.
\end{align*}

We see that
 the derivation $\th$
satisfies the same equation $i_\th\om=dp$ as 
the derivation $\th_p$
 attached to the  image of $p$ in $\PP\br$ via the map~\eqref{th_p}.
Hence, $\th=\th_p$ and, therefore, $\th_p(a)=\epsilon(a)\cdot\partial_{a^*}p$.

The rest of the proof  is straightforward and is left
to the reader.
\qed

\section{Proof of Theorem \ref{main}}
\subsection{} Given a $\Z$-graded vector space  $M=\oplus_{k\in\Z}M_k$, 
let $P(M,t)=$\break
$\sum_{k\in\Z}\dim M_k\cdot t^k$
$\in \k[[t,t\inv]]$  denote its  Poincar\'e series.
If each graded piece $M_k$ is an 
{$R$-bimodule,} then one can define
a refined matrix-valued  Poincar\'e series\break
$\pmat(M,t)$, that is an $I\times I$-matrix whose
entries are elements of $\k[[t,t\inv]]$  given by
$\pmat(M,t)_{ij}:=\sum_{k\in\Z}\dim (e_i M e_j)\cdot t^k,\,i,j\in I.$
With this notation, we have
\beq{pmatrix}
P(M^R,t)=\Tr\pmat(M,t).
\eeq
\subsection{} From now on, we
assume that $Q$ is neither Dynkin
nor extended Dynkin.

 We are interested in
the Hochschild cohomology $\HH^\bullet(\Pi)$
of the algebra $\Pi$.
Observe that the natural  grading on $\Pi$ induces a grading on
each  Hochschild cohomology group.

\begin{lem}\label{hoch}
 The groups $\HH^i(\Pi)$ vanish for all $i> 2$, and
we have
$$P(\HH^0(\Pi),t)
-P(\HH^1(\Pi),t)
+P(\HH^2(\Pi),t)=|I|.
$$
\end{lem}
\begin{proof} Since
 $Q$ is neither Dynkin
nor extended Dynkin, the  algebra $\Pi$ is known
to be a Koszul algebra, see [MV], [MOV]. Specifically,
a
$\Pi$-bimodule  resolution of $\Pi$
is provided by the following  Koszul complex
\beq{kos0}
\Bigl(0\to\Pi\otimes_R \Pi\map \Pi \otimes_R E\otimes_R \Pi\map
 \Pi\otimes_R \Pi\Bigr)\stackrel{m}\onto\Pi,
\end{equation} 
where $E$ denotes the $R$-bimodule generated by the edges
of $\overline{Q}$.
Therefore, the Hochschild cohomology of $\Pi$ may be computed
by applying
 $\Hom_{\bimod \Pi}(-,\Pi)$  to the complex above.
This way, we  get the following complex
\beq{kos}\Pi^R\map  (\Pi\otimes E)^R\map  \Pi^R.
\end{equation} 
The vanishing part of the Lemma is now clear.

Let $\chimat$ denote the  matrix-valued graded 
Euler characteristic, that is, 
 the alternating sum 
 of the matrix-valued  Poincar\'e series
of the complex $\Pi\to\Pi\o E\to\Pi.$
The matrix $\chimat$ 
is easily expressed in terms  of the matrix $P:=\pmat(\Pi,t)$
and the adjacency matrix $C$  of the graph $\overline{Q}$
by the formula $\chimat=t^2\cdot P-$
${t\cdot P\cdot C+P}.$

Further,
by formula \eqref{pmatrix},
the Euler characteristic
of the complex \eqref{kos} is an element of $\k[[t]]$
which is equal to $\Tr(\chimat)$.
 Now, according to  [MOV], one has a matrix identity
$P=1/(1-C\cdot t+t^2)$. Thus, we compute
$$\Tr(\chimat)=
\Tr(t^2\cdot P-t\cdot P\cdot C+P)=
\Tr (P\cdot(1-C\cdot t+t^2))=\Tr(\Id)=|I|,
$$
and we are done.
\end{proof}

\subsection{Proof of formula \eqref{formula}.}
Case $i=0$ of the formula follows from Proposition \ref{hh0}
(to be proved later, in Sect.\,\ref{proof_hh0}).

Next, one observes  that the complex
 inside the parenthesis in \eqref{kos0} is {\em self-dual}.
In other words, applying the functor $\Hom_{\bimod{\Pi}}(-,\Pi\o\Pi)$
to this complex one obtains the same complex again. Therefore,
using   \eqref{kos0} as a resolution for the computation of
$\Ext^\bullet_{\bimod{\Pi}}(\Pi,\Pi\o\Pi)$, we see that the only
nonzero cohomology group sits in degree two and we have
 $\Ext^2_{\bimod{\Pi}}(\Pi,\Pi\o\Pi)\cong\Pi.$
We conclude that formula \eqref{vdb2} holds for
the algebra
$A=\Pi$ and $d=2$.
Thus, from  Van den Bergh's  isomorphism
\eqref{vdb_iso} we deduce
$\HH^i(\Pi)$ $=\HH_{2-i}(\Pi),$ for any $i$.
In particular, we get
$\HH^2(\Pi)=\HH_0(\Pi)={\Pi/[\Pi,\Pi]=L.}$

We have already computed $\HH^i(\Pi)$ for $i=0,2$.
Furthermore, we know that $\HH^i(\Pi)=0$ for all $i>2$,
by the first claim of  Lemma \ref{hoch}.
Therefore, from the formula for the Poincar\'e series
given in Lemma \ref{hoch} we deduce
\begin{align*}
P(\HH^1(\Pi),t)&=-|I|+P(\HH^0(\Pi),t)+P(\HH^2(\Pi),t)\\
&=-|I|+P(Z(\Pi),t)+P(L,t)=P(L,t)-|I|+1\\
&=
P(L/R,t)+1=P(\k\cdot\eu_\text{Lie}\ltimes(L/R),\,t).
\end{align*}
 By comparing the two sides, we see that the
case $i=1$ of \eqref{formula},
as well as part (i) of the Theorem,  follows by dimension argument from
the injectivity statement in Proposition
 \ref{stillinj}. This completes the proof of formula
\eqref{formula}.
\qed

\section{Deformations of the preprojective algebra}
\label{pif} 
\subsection{}\label{pifpif} It is well-known that first order
infinitesimal deformations of an associative algebra
are controlled by the second  Hochschild cohomology
group, and obstructions to  deformations 
 are controlled by the third Hochschild cohomology
group. Thus, according to formula
\eqref{formula}, first order
infinitesimal deformations of the preprojective algebra
$\Pi$ are unobstructed, and its versal deformation is parametrized
by the vector space $L$.

We will show that, in the non-Dynkin case, {\em any}
formal (infinite order) one-parameter deformation of the algebra $\Pi$
 is equivalent to one obtained by deforming 
the defining relation in the
 preprojective algebra as follows:
 $$
\dis\sum\nolimits_{a\in Q}\,[a,a^*]=t\cdot f,
$$
 where $t$ is the  deformation parameter
 and $f$ is an element of $\PP[[t]] $
that commutes with $R$.

In more detail, from the general
formula \eqref{MR} applied to the
algebra $A=\Pi$ we deduce that one can find
a  graded subspace
$\PP_L\sset \PP^R$ which
is complementary
to the graded space
$\Ker[\PP\onto L]=[\PP,\PP]+\PP\omh\PP\sset\PP$,
i.e., such that we have
\beq{gr_comp}
\PP=\PP_L\bigoplus ([\PP,\PP]+\PP\omh\PP),
\quad\text{and}\quad
\PP_L\sset \PP^R.
\eeq
We fix such a subspace $\PP_L$ once and for all, and also
choose and fix a graded   $\k$-basis
$\{f_j\}_{j\in \N}$  of  $\PP_L$.
 For each $j\in \N$,
introduce a formal variable $t_j$ that
may be thought of as the $j$-th coordinate function
on  $\PP_L$ with respect to the basis $\{f_j\}_{j\in \N}$.

Given a  $\k$-vector space $V,$ write  $V[[\bt]]:=
\limp V[[t_1,...,t_n]] $
for the vector space of $V$-valued formal power series in
the (infinitely many)  formal variables $t_j$. Thus, 
$\k[[\bt]]$ is a complete topological local $\k$-algebra
with maximal ideal $(\bt)\sset\k[[\bt]]$ formed by the 
series without constant term.
For any $V$, the space $V[[\bt]]$
acquires a natural  structure
of complete topological $\k[[\bt]]$-module.
In particular, we have a  $\k[[\bt]]$-module
$\PP_L[[\bt]]$.
Observe that
the identity map
$\Id:\PP_L\to\PP_L$ may be viewed as an element of
$\PP_L[[\bt]]$ that has the form
$\sum_{j\in\N}\,t_j\cdot f_j$.

We define a formal deformation, $\pif$, of
the preprojective algebra $\Pi$
 as  the quotient of the algebra $\PP[[\bt]]$ modulo the {\em closed}
two-sided ideal (topologically) generated by the element
\begin{equation}\label{fj}
\sum\nolimits_{a\in Q}\, [a,a^*]-\sum\nolimits_{j\in\N}\, t_j\cdot f_j.
\end{equation}
Thus,  $\pif$ is a complete  topological $\k[[\bt]]$-algebra
that may be thought of as a deformation of $\Pi$
 with formal parameters 
$\{t_j\}$, i.e., the base of the deformation is 
$\Spec\bigl(\k[[\bt]]\bigr)$, a formal scheme.

\begin{thm}\label{def1} If the quiver $Q$ is not Dynkin (but possibly extended 
Dynkin) then the deformation $\pif$ 
is flat over  $\k[[\bt]]$, and 
is a versal deformation of $\Pi.$ 
\end{thm}

\subsection{} 
We begin the proof of Theorem \ref{def1} with an easy lemma. 

Let the vector space $\widetilde\Pi=\Pi\otimes\k[t]/(t^2)$ be 
equipped with a $\k[t]/(t^2)$-algebra structure given
by a star-product  of the form 
$p,q\mto p\star q=pq+ t\cdot\beta(p,q)$.
Thus, $\widetilde\Pi$ is a first order deformation of the
algebra $\Pi$.
We
write $[\beta]\in\HH^2(\Pi)$ for the  Hochschild cohomology
class  corresponding to this  deformation.

Recall the notation $\Pi^R$
for the centralizer of $R$ in $\Pi$.

\begin{lem}\label{class} Assume that
 in   $\widetilde\Pi$
we have 
$$\sum_{a\in Q} (a\star a^*-a^*\star a)= t\cdot c\enspace \modu (t^2),\quad
\text{for some}\quad c\in\Pi^R.$$
Then, the class $[\beta]\in\HH^2(\Pi)$ goes, via the isomorphism
$\HH^2(\Pi)\iso L$ of Theorem \ref{main}, cf. also Sect.
\ref{dyn},
to the image of $c\in\Pi$ under the projection $\Pi\onto \Pi/[\Pi,\Pi]=L$.
\end{lem}
\begin{proof} 
Recall that the Koszul resolution \eqref{kos0} of the algebra $\Pi$
admits a standard imbedding, as a subcomplex,
into the Bar-resolution of $\Pi$:
\begin{equation}\label{bar}
\xymatrix{
\text{Koszul:}&\Pi\otimes_R\Pi\,\ar[r]\ar@{^{(}->}[d]^<>(0.5){\imath_2}
&\Pi\otimes_R E\otimes_R\Pi\,\ar[r]\ar@{^{(}->}[d]^<>(0.5){\imath_1}
&\Pi\otimes_R\Pi\ar@{=}[d]^<>(0.5){\Id}\\
\text{Bar:}&\Pi\otimes_RT_R^2\Pi\otimes_R\Pi\ar[r]&
\Pi\otimes_R\Pi\otimes_R\Pi\ar[r]&
\Pi\otimes_R\Pi,
}
\end{equation}
In this diagram, $\imath_1$ is a  $\Pi$-bimodule map
induced by the natural imbedding $E\into \Pi$
and $\imath_2$ is a  $\Pi$-bimodule map
defined 
by the formula
$$
\imath_2
:\, p\otimes q\mto \sum_{a\in Q} (p\otimes a\otimes
a^*\otimes q-p\otimes a^*\otimes a\otimes q).
$$

Now, given a deformation of $\Pi$ as in the lemma,
the class $[\beta]\in\HH^2(\Pi)$ corresponds
to a 2-cocycle in 
 $\Hom_{\bimod{\Pi}}(\Pi\otimes_RT_R^2\Pi
\otimes_R\Pi,\Pi)$ given by 
$$\widetilde{\beta}: \,
\Pi\otimes_RT_R^2\Pi\otimes_R\Pi\too
\Pi,\quad
p\otimes u\otimes v\otimes q\mto p\cdot\beta(u,v)\cdot q.
$$
It is clear that, under the imbedding of complexes
\eqref{bar}, the cocycle $\widetilde{\beta}$ restricts
to a map $\imath_2^*(\widetilde{\beta})\in
\Hom_{\bimod{\Pi}}(\Pi\otimes_R\Pi,\Pi)$
such that $\imath_2^*(\widetilde{\beta}):\,
p\otimes q\mto p\cdot c\cdot q$. 
Hence, the cohomology class
of $\imath_2^*(\widetilde{\beta})$ in the complex \eqref{kos}
is represented by the element
$c\in\Pi^R$. The lemma follows.
\end{proof}

\subsection{} By Theorem \ref{main}, we know that
 deformations of $\Pi$ are unobstructed. 
Thus, general deformation theory implies that there exists a 
formal flat multi-parameter deformation $\widetilde\Pi$ of 
the algebra
$\Pi$, parametrized by the formal neighborhood of zero in 
$\HH^2(\Pi)=L.$ 

Recall that we have chosen
a vector space $\PP_L\sset\PP$.
It is clear from  \eqref{gr_comp}
that the projection  $\PP\onto L$ induces an isomorphism
$\PP_L\iso  L$.

We use the last isomorphism to identify
$\HH^2(\Pi)$ with $\PP_L$.
Thus, the deformation  $\widetilde\Pi$
becomes  parametrized by the formal neighborhood of zero in 
$\PP_L$. Therefore, there exists
 an
 $R [[\bt]]$-bimodule isomorphism
$\widetilde\Pi\cong \Pi[[\bt]]$, equal to the identity modulo 
the ideal $(\bt)$, in the
notation of \S\ref{pifpif}.
We transport the algebra structure
on  $\widetilde\Pi$ to   $\Pi[[\bt]]$ via that isomorphism.
This way,
we make  $\Pi[[\bt]]$ a complete topological $\k[[\bt]]$-algebra,
to be denoted ${{}^\star\Pi[[\bt]]}$.

Further, for each $k\geq 0$, let $\Pi[[\bt]]_k$ denote
the space of $\Pi$-valued series
homogeneous in  $t_j$'s of degree $k$.
Writing a power series in the variables $t_j$
as  $p=\sum_{k\geq 0}\,p_k, \,p_k\in\Pi[[\bt]]_k$, yields
a direct  product decomposition
$\Pi[[\bt]]=\prod_{k\geq 0}\,\Pi[[\bt]]_k$.
Thus, multiplication in the algebra  ${{}^\star\Pi[[\bt]]}$
takes the form of 
 {\em star-product}:
$$(p,q)\mto p\star q= pq +\beta_1(p,q)+\beta_2(p,q)+\ldots,
\quad\beta_k: \Pi\times\Pi\map \Pi[[\bt]]_k.
$$

\subsection{} We return now to the setup of Theorem \ref{def1}
and consider the $\k[[\bt]]$-algebra $\pif$.
In Sect. \ref{below} we will prove the following

\begin{lem}\label{phi} There exists a continuous algebra automorphism
$\varphi: \k[[\bt]]\iso$
$\k[[\bt]]$ and a continuous
$\k$-algebra isomorphism
$\Phi: \pif\iso{{}^\star\Pi[[\bt]]}$ such that

\npb{The map $\Phi$ is $\varphi$-semilinear,
i.e., for any $f\in \k[[\bt]]$ and
$p\in \pif$, one has
$\Phi(f\cdot p)=\varphi(f)\cdot\Phi(p)$;}

\npb{The differential of $\varphi$ at the origin
$0\in\PP_L$ equals ${d\varphi=\Id:
\PP_L\to\PP_L,}$  the identity map.} 
\end{lem}

This Lemma easily implies Theorem \ref{def1}.
Indeed,
the deformation ${}^\star\Pi[[\bt]]$ being flat by construction,
 it follows from
Lemma \ref{phi}  that the  deformation  $\pif$ is also flat 
and provides a versal deformation. 
The theorem is proved. \qed

\subsection{Proof of Lemma \ref{phi}.}\label{below}
Recall the vector space $\PP_L\sset \PP$
and let $\Pi_L$ denote its image under
the projection $\PP\onto\Pi$. By \eqref{gr_comp},
we have 
\beq{complement}
\Pi=\Pi_L\bigoplus[\Pi,\Pi],\quad
\text{and}\quad
\Pi_L\sset\Pi^R.
\end{equation}
Furthermore  we have natural isomorphisms
$\PP_L\iso\Pi_L\iso L,$ where the first map is induced by the
projection $\PP\onto\Pi$ and the second map
by the projection $\Pi\onto L$.

Next, we introduce the notation $\bbm$ for a multi-index, that is
for a function $\bbm:\N\to \{0,1,\ldots\}$ such that 
$\bbm (j)=0$ for all but finitely many $j$'s.
Given such a multi-index $\bbm$,
write $t^\bbm:=t_1^{\bbm(1)}t_2^{\bbm(2)}\ldots$,
for the corresponding monomial in the  variables $t_j$.
Put $|\bbm|=\sum_{j\in\N}\,\bbm(j).$

To prove the lemma,
we must show that there exist algebra homomorphisms 

$$\Phi: \pif\to {{}^\star\Pi[[\bt]]},\quad
\text{and}\quad
\varphi: \k[[\bt]]\to \k[[\bt]],\, \bt\mto\varphi(\bt)=:\bt'=\{t'_j\}_{j\in{\mathbb{N}}},
 $$
where $t'_j=t_j+\sum_{|\bbm|>1}
\lambda_{i\bbm}t^\bbm, \lambda_{i\bbm}\in \k$
(formal change of variables), 
such that $\Phi$ is the identity in degree $0$ and 
$\Phi(a)=a+\sum_{|\bbm|>0} t^\bbm z_\bbm(a)$ for each edge $a
\in\overline Q$. 
(here if $a\in e_i\Pi e_j$ then $z_\bbm(a)$ also must 
belong to $e_i\Pi e_j$
to ensure the relations between projectors $e_i$ and $a$
are respected by $\Phi$).

We construct $\Phi$ by induction in $|\bbm|$. First we
consider the case $|\bbm|=1$, i.e., 
$t^\bbm=t_j$ for some $j$.

In the algebra ${{}^\star\Pi[[\bt]]}$, we have an expansion
$$\sum_{a\in Q}\, (a\star a^*-a^*\star a)=\sum_{|\bbm|>0} 
t^\bbm\cdot y_\bbm,\quad y_\bbm\in \Pi^R. $$
Hence,  Lemma \ref{class}
implies that the image of $y_\bbm\in \Pi^R $
under the projection to $L$ is equal to $f_j. $
Therefore, the element 
$z_j(a)$ in the expansion
$\Phi(a)=a+\sum_{|\bbm|>0} t^\bbm z_\bbm(a)$
must satisfy the equation 
$$f_j-y_j=\sum_{a\in Q} ([z_j(a),a^*]+[a,z_j(a^*)]).
$$
Now, since $y_j$ projects to $f_j\in L$, we have 
 $f_j-y_j\in [\Pi,\Pi]$.
Since  $[\Pi,\Pi]$ is spanned by elements of the form
$[a,z],[a^*,z]$ we see that the above equation
 does have some solution ~$z_j(a)$.

Next, let $|\bbm|\ge 2.$ Then, the elements
$z_\bbm(a)$ are determined by the equation 

$$y_\bbm+\sum_a([a,z_\bbm(a^*)]+[z_\bbm(a),a^*])+
\sum_{{\mathbf{k}}+{\mathbf{l}}=
\bbm,|{\mathbf{k}}|>0,|{\mathbf{l}}|>0}
\sum_a [z_{\mathbf{k}}(a),z_{\mathbf{l}}(a^*)]\in \Pi_L. 
$$
This equation again has solutions, since $\Pi=[\Pi,\Pi]\oplus \Pi_L$
and $\Pi_L\sset\Pi^R,$ by  \eqref{complement}. 
For such a solution, we have 
$$y_\bbm+
\sum_a([a,z_\bbm(a^*)]+[z_\bbm(a),a^*])+\sum_{{\mathbf{k}}+{\mathbf{l}}=
\bbm,|{\mathbf{k}}|>0,|{\mathbf{l}}|>0}
\sum_a [z_{\mathbf{k}}(a),z_{\mathbf{l}}(a^*)]=
\sum \lambda_{j\bbm}f_j.  
$$
Continuing so, we obtain $z_\bbm(a)$ 
for all $\bbm$, and we see 
that because of the construction of $z_\bbm(a), \Phi$
 indeed defines a homomorphism as required, with  
$t_j'=t_j+\sum_{|\bbm|>1}\lambda_{i\bbm}t^\bbm.$

\section{Representation schemes}
\label{pc}
\newcommand{\OO}{{\mathsf{O}}}

\subsection{Stabilization.} Fix a quiver $Q$ with vertex set $I$
and let $R=\k I$.

Let  $a_1,\ldots,a_p\in \k Q$  be
a finite (possibly empty) collection of homogeneous elements
in the path algebra of $Q$, of some 
degrees $\deg a_j\geq 2.$ Let $\BI\sset\k Q$ denote the
two-sided ideal generated by the elements $a_1,\ldots,a_p$,
and set $A=\k Q/\BI.$
Thus, $\BI$ is a graded ideal and $A=\bigoplus_{r\geq 0} A(r)$ is
a graded algebra  such that $A(0)=R$.

Let $\mathbf d=(d_i\in \Z_{\geq 0},\, i\in I)$
be  an $I$-tuple, to be referred
to as {\em dimension vector}.
Let  $V=\bigoplus_{i\in I}\, V_i$ be
an $R$-module such that
$\dim V_i=d_i,\,\forall i\in I$.
We will write
 $\Rep (A,\mathbf d):=\Rep(A,V)$
for the corresponding  representation scheme
and, following Sect. \ref{rep_sch},
put $G_\bd:=\left(\prod_i \GL(V_i)\right)\big/{\mathbb G}_m.$

We consider the algebras
$\oo_\bd(A)$ and $\k[\Rep(A,\bd)]^{G_\bd}$, cf. Sect. \ref{eval}.
The grading on $A$ makes $\oo_\bd(A)$, resp.,
$\k[\Rep(A,\bd)]^{G_\bd}$,
a graded algebra with respect to the {\em total} grading
$\oo_\bd(A)=\bigoplus_{r\geq 0} \oo_\bd(A)(r)$,
resp., $\k[\Rep(A,\bd)]^{G_\bd}=\bigoplus_{r\geq 0} \k[\Rep(A,\bd)]^{G_\bd}(r)$.

In  \eqref{trace}, we
have defined an algebra map
$\psi_\bd=\Tr\ccirc\ev_\bd: \oo_\bd(A)\too
\k[\Rep(A,\bd)]^{G_\bd}$.
This map 
is clearly a graded algebra homomorphism.

The following result is a refined version of
the `stabilization phenomenon' observed in \cite[\S4]{Gi}.

\begin{prop}\label{iso1}
For any  positive integer  $r$, there exists $N(r)\gg0$ such that 
the following holds:

{\em If ${\mathbf d}=\{d_i\}_{i\in I}$ is such
that $d_i\ge N(r)$ for all $i\in I$, then
the map  $\psi_\bd: \oo(A)(s)\to\k[\Rep_\bd A]^{G_\bd}(s),$
the restriction of $\psi_\bd$ to
the homogeneous component of degree $s,$ is a bijection
 for all $0\leq s\leq r$.}
\end{prop}

First, we are going to prove the proposition in the special
case where $\BI=0$,    that is, for $A=\k Q.$
In this case, we have $\Rep(A,\bd)=\Rep(Q,\bd)$ is the variety of 
$\mathbf d$-dimensional representations 
of the  quiver $Q$. This variety is a vector space.

\begin{lem}\label{iso} Proposition \ref{iso1} holds for $A=\k Q$.
\end{lem}

Let $\BL$ denote the direct sum of homogeneous components of $(\k Q)\br$
of strictly positive
degrees.
Thus,  we have a vector space direct sum decomposition
$(\k Q)\br=R\bigoplus\BL$.
The composite $\Sym(\BL)\into\Sym((\k Q)\br)
\onto\oo_\bd(\k Q)$ is clearly a graded algebra  isomorphism.

 The proof of Lemma \ref{iso} given below copies the argument
in the proof of \cite[Proposition 4.2]{Gi}; it is based on the fact that 
there are no polynomial identities which are satisfied in matrix
algebras of all sizes. 

\begin{proof}[Proof of Lemma \ref{iso}.]
To simplify the exposition, we will present the proof 
in the case where $Q$ has one vertex and two edge-loops at that vertex; the general case is 
entirely similar. 
Thus, $\mathbf d=d$ is a single integer,
so $\Rep (Q,\mathbf d)$ is the space 
of pairs of $d\times d$-matrices $X,Y$,
and the vector space
 $\BL$ is spanned by nonempty 
cyclic words $w$ in the alphabet with  two letters $x,y$.

Assume the statement of the lemma is not true.
Thus, there exists an element $g\in \Sym\BL$
such that ${\psi}_\bd(g)=0$ for all $d$. In other words, 
we have a relation of degree $r$ of the form 
$$
\sum\nolimits_J c_J\cdot \prod_w \Tr(w(X,Y))^{J(w)}=0,\ (c_J\in \k)
$$
which is satisfied for any matrices $X$ and $Y$.
Here $J$ is a function on the set of cyclic words with nonnegative 
integer values and finite support. 

We may assume without loss of generality
that at least one cyclic word present in this relation involves $x$. 
Let us differentiate the relation with respect to $X$. 
We get 
$$
\sum_J \sum_u J(u)\cdot c_J\cdot  \prod_w \Tr(w(X,Y))^{J(w)-\delta_{wu}}\cdot
\sum_{k: u_k=x}(u_{k+1}...u_{k-1})(X,Y)=0,
$$ 
where $u_k$ is the k-th letter of the word $u$ (we choose a 
representation of  every cyclic word by a usual word, and 
agree that $m\pm 1=\pm 1$ 
if $u$ has length $m$). But for large enough $d$, 
all words of length $\le r$ of matrices $X,Y$ 
are linearly independent for generic $X,Y$ 
(for example, take $X,Y$ to be the operators by which $x,y$ act
in the quotient of the free algebra in $x,y$ by the $r$-th power of
the augmentation ideal). Thus, we have 
$$
\sum_J J(u)\cdot c_J \cdot \prod_w \Tr(w(X,Y))^{J(w)-\delta_{wu}}=0
$$
for each $u$ which contains $x$. 
This gives a relation of degree smaller than $r$, 
so we are done by using induction. 
\end{proof}

\begin{proof}[Proof of Proposition \ref{iso1}.]
The algebra projection $\k Q\onto A$ makes
$\Rep(A,\mathbf d)$ a closed  subscheme
in
 $\Rep(\k Q,\mathbf d),$ not necessarily reduced
in general.
The defining ideal of this subscheme
 is  generated by 
the matrix elements of the matrix valued functions $\wh{a}_1,\ldots,
\wh{a}_p$.

To avoid confusion, we will use boldface notation
for the map $\bpsi_\bd: \oo_\bd(\k Q)\to\k[\Rep(Q,\bd)]^{G_\bd}$
and reserve the notation $\psi_\bd$
for the map $\oo_\bd(A)\to\k[\Rep(A,\bd)]^{G_\bd}.$

Assume the statement of the proposition does not hold. Then there exists a homogeneous element
$g\in \oo_\bd(A)$, of degree $0<s\leq r$,
and a dimension vector $\mathbf d$ such that  $d_i\ge N(r),\,\forall i,$ 
and such that one has $\psi_\bd(g)=0$. 
Let $\pa{g}$ be a representative of
$g$ in $\Sym(\BL)$. Thus,
${\bpsi}_{\mathbf d}(\pa{g})$ is a  degree $r$ homogeneous
 polynomial on the vector space $\Rep(\k Q,\mathbf d)$.
 
By construction, the polynomial ${\bpsi}_{\mathbf d}(\pa{g})$
vanishes on the subscheme $\Rep (A,\mathbf d)$.  Hence, 
${\bpsi}_{\mathbf d}(\pa{g})$ is a $G_{\mathbf d}$-invariant
homogeneous polynomial that belongs to the ideal
 generated by 
the matrix elements of the functions $\wh{a}_1,\ldots,
\wh{a}_p$.
Therefore, using Weyl's first fundamental theorem 
of invariant theory we deduce that 
${\bpsi}_{\mathbf d}(\pa{g})$ is a linear combination of products of the form
$\Tr (\wh{w}_1)\cdot\ldots\cdot\Tr (\wh{w}_m)\cdot
\Tr (\wh{w}_{m+1}\wh{a}_j)$,
where $j\in [1,p]$ and  $w_1,\ldots,w_{m+1}$ are certain  paths in $Q$
of total length $\ell(w_1)+\ldots+\ell(w_{m+1})+\ell(a_j)=r$. 

Now, Lemma \ref{iso} says
 that the element $\pa{g}$ must be equal to the corresponding 
linear combination of the  elements
$w_1\&\ldots\& w_m\& (w_{m+1}a_j)$
$\in \Sym^{m+1}(\BL)$.
It is clear that all elements of this form generate
the kernel of the algebra projection
$\Sym(\BL)\onto\oo_\bd(A).$ 
It follows  that $g=0$, as desired. 
\end{proof}

\subsection{Quiver varieties.}\label{quiver} 
Let $\overline Q$ be the double of a quiver $Q$,
and write $\PP:=\k\overline Q$ for the corresponding path algebra.

For any dimension vector $\bd\geq 0$,
the vector space $\Rep (\overline Q,\mathbf d)=T^*\Rep(Q,\mathbf d)$ is
equipped with a natural symplectic
structure. It is easy to see that the corresponding
symplectic 2-form on $\Rep (\overline Q,\mathbf d)$
equals $\Tr\wh{\om}$, the image of the 2-form $\om\in\DR^2_R\PP$
under the evaluation map, see Sect. \ref{eval}.

Fix  a dimension vector $\bd$. Recall that 
the Lie algebra $\fg_\bd=\Lie G_\bd$, as well as  its dual $\fg_\bd^*$, are
both identified with the
a codimension 1 hyperplane formed by the $I$-tuples
$(x_i\in \Lie\GL(d_i))_{i\in I}$ such that $\sum_i \Tr x_i=0$.

The group $G_\bd$ acts
 linearly
on $\Rep (\overline Q,\mathbf d)$. 
This action 
is Hamiltonian, see  Sect. \ref{eval}, and the corresponding
 moment map reads, cf. \cite{CB1}:
\beq{mu_map}\mu: \Rep (\overline Q,\mathbf d)\too \fg_\bd,
\quad\rho\mto\mu(\rho)=\sum_{a\in Q} [\rho(a),\rho(a^*)].
\eeq

In this section, we are interested in  the algebra $\Pi$, the preprojective algebra
of $Q$. We have $\Pi=\PP/\BI$,
where $\BI$ is the two-sided ideal generated by the single element
$\omh=\sum_{a\in Q}[a,a^*]$.
The algebra $\Pi$ is the Hamiltonian reduction of $\PP$
with respect to this element.
The projection $\PP\onto\Pi$
makes $\Rep (\Pi,\mathbf d)$, the scheme
of all $\mathbf d$-dimensional representations 
of the algebra $\Pi$, a closed subscheme in
$\Rep (\overline Q,\mathbf d)$.
We have 
$\Rep (\Pi,\mathbf d)=\mu\inv(0)$, the scheme theoretic  zero fiber of the moment map.

We let
$M(Q,\mathbf d):=\Spec \k[\Rep (\Pi,\mathbf d)]^{G_\bd}$ be the corresponding
categorical quotient, an affine subscheme of $\Rep (\overline Q,\mathbf d)/\!/{G_\bd}
:=\Spec \k[\Rep (\overline{Q},\mathbf d)]^{G_\bd}$
 of finite type.
The scheme  $M(Q,\mathbf d)$ is called {\em quiver variety},
it is a Poisson scheme that may be thought of as a hamiltonian reduction of 
$\Rep (\overline Q,\mathbf d)$.  
By Theorem \ref{moment_thm}, the
map $\psi_\bd: \oo_\bd(\PP)\map\k[\Rep (\overline Q,\mathbf d)]^{G_\bd}$
is a surjective morphism of graded Poisson algebras
that descends to a surjective
graded
Poisson algebra  homomorphism 
$\psi_\bd: \oo_\bd(\Pi) \onto \k[M(Q,\mathbf d)]^{G_\bd}$.

\subsection{}
We need to recall some results from \cite{CB1} about quiver varieties. 
Let $C$ denote the adjacency matrix
of the quiver $Q$ and write
$A=2{\rm Id}-C$ for the corresponding Cartan matrix.
Define the Tits quadratic form on the vector space 
$\k^{|I|}$
 by the formula
$$
q(\mathbf d)=\frac{1}{2}(\mathbf d,A\mathbf d).
$$
Let $p(\mathbf d)=1-q(\mathbf d)$. 
 
Let $R_+$ be the set of positive roots of $Q$.
Let $\Sigma_0$ denote the set of positive 
dimension vectors $\beta$
such that whenever $\beta$ is written as a sum of two or more positive roots, 
$\beta=\beta_1+...+\beta_m$, one has $p(\beta)>p(\beta_1)+...+p(\beta_m)$.

\begin{thm} \label{quivar} (\cite{CB1}). 
\vi   $\mathbf d\in \Sigma_0$ if and only if there is a simple representation 
of $\Pi$ with dimension vector $\mathbf d$. 

If $\mathbf d\in \Sigma_0$ then:

\vii  the schemes $\Rep (\Pi,\mathbf d)$ and 
$M(Q,\mathbf d)$ are reduced, i.e., they are 
affine algebraic varieties. Moreover, these varieties are irreducible,
and their generic points correspond to simple representations;

\viii $\Rep (\Pi,\mathbf d)$ is a complete intersection
in $\Rep (\overline{Q},\mathbf d)$.

\iv The Poisson structure on $M(Q,\mathbf d)$
 is generically symplectic.\qed
\end{thm}
 
 One can   pull-back polynomial functions via the moment  map
\eqref{mu_map}. This yields an algebra homomorphism
$\mu^*:\k[\fg_\bd]\to\k[\Rep (\overline{Q},\mathbf d)],$ and we have

\begin{cor}\label{pasha_cor}
Assume that
 $\mathbf d\in \Sigma_0$. Then the Poisson center
of the algebra $ \k[\Rep (\overline{Q},\mathbf d)]^{G_\bd}$
equals $\mu^*\left(\k[\fg_\bd]^{G_\bd}\right),$
the pull-back of the algebra of $\Ad G_\bd$-invariant polynomials
on the Lie algebra $\fg_\bd.$
\end{cor}
\begin{proof} The assumption on $\bd$ insures
that the moment map
$\mu: \Rep (\overline Q,\mathbf d)\to\fg_\bd$ is flat, in particular,
it is surjective. This follows from 
Theorem \ref{quivar}(iii), since the dimension
of any fiber of $\mu$ is not less the dimension
of the zero fiber, and the latter is a complete
intersection. Furthermore, all fibers  of $\mu$
are irreducible. To see this, let $x\in\fg_\bd.$
The standard increasing filtration on
the algebra $ \k[\Rep (\overline{Q},\mathbf d)]$,
by degree of the polynomial, induces
an  increasing filtration on the quotient algebra
$\k[\mu^{-1}(x)]$. It is well known that,
since $\mu$ is flat and the zero fiber of $\mu$ is reduced,
for the corresponding associated graded
algebra, one has
$\gr\k[\mu^{-1}(x)]\cong\k[\mu^{-1}(0)]$.
But, the scheme $\mu^{-1}(0)$ is irreducible by
Theorem \ref{quivar}(ii).
It follows that the algebra $\k[\mu^{-1}(0)]$,
hence also $\k[\mu^{-1}(x)]$,
has no zero-divisors. Thus, the scheme
$\mu^{-1}(x)$ is irreducible as well.

Next, we claim that for any sufficiently general  $\Ad G_\bd$-conjugacy
class $\OO\sset\fg_\bd$, the preimage $\mu^{-1}(\OO)$ is a smooth submanifold in
$\Rep (\overline Q,\mathbf d)$ and, moreover, the action of $G_\bd$ on
 $\mu^{-1}(\OO)$ is free. Indeed, the former statement follows from
the latter by general properties of moment maps. To prove the freeness,
 we let the conjugacy class $\OO$ be such that any
$x\in \OO$ is a diagonalizable endomorphism, say $x=\text{diag}(s_1,\ldots,s_d)$,
and, in addition,  such that for any proper subcollection
$S\sset \{s_1,\ldots,s_d\}$ one has
$0\neq\sum_{s\in S} s.$  Let $\rho\in \mu^{-1}(\OO)$   be a
 representation of 
$\overline Q$ and assume that $E\sset\k^{\bd}$ is a nontrivial
subrepresentation. Therefore, the vector space $E$ is $\rho(a)$-stable
for any endomorphism $\rho(a), a\in \overline Q$.
Hence, we must have
$\Tr x|_E= \sum_{a\in Q} \Tr\left([\rho(a)|_E,\rho(a^*)|_E]\right)=0,$
cf. \eqref{mu_map}. On the other hand, it is clear that
$E$ is the span of certain eigen-spaces of $x$ and thus
$\Tr x|_E=\sum_{s\in S} s,$ for some
 subcollection
$S\sset \{s_1,\ldots,s_d\},$ where $|S|=\dim E$.
This contradicts our condition on the eigenvalues of $x$.
We conclude that any point
$\rho\in \mu^{-1}(\OO)$ is an irreducible representation of 
$\overline Q$, and our claim follows by the Schur lemma.

Thus, we have proved that for any sufficiently  general 
  $\Ad G_\bd$-conjugacy
class $\OO\sset\fg_\bd,$ the preimage $\mu^{-1}(\OO)$
is a smooth submanifold with a free
$G_\bd$-action. This implies that all $G_\bd$-orbits in
 $\mu^{-1}(\OO)$ are closed and correspond bijectively
to points of the categorical  quotient  $\mu^{-1}(\OO)/\!/{G_\bd}=
\Spec \k[\mu^{-1}(\OO)]^{G_\bd}$. Furthermore, the categorical
quotient is a smooth affine symplectic algebraic variety.

Now, let $F$ be a central element of the Poisson
algebra $\k[\Rep (\overline{Q},\mathbf d)]^{G_\bd}$.
Then, for any general  conjugacy classes $\OO$,
the restriction of $F$ to $\mu^{-1}(\OO)$ is
a central element of the Poisson
algebra $\k[\mu^{-1}(\OO)]^{G_\bd}=\k[\mu^{-1}(\OO)/\!/{G_\bd}].$
The  variety  $\mu^{-1}(\OO)/\!/{G_\bd}$ being
 symplectic, we deduce that
the polynomial $F$ must be constant on connected components of
any  general fiber of the
moment map. Further, we know that any such fiber
is irreducible, hence connected.
Thus, $F$ is constant on general fibers.
We deduce, since  $\fg_\bd$ is smooth, that
$F=\mu^*(\bar F)$ for some
$\bar F\in\k[\fg_\bd].$
Finally, since $F$ is  $G_\bd$-invariant,
one may arrange  that $F\in
 \mu^*\left(\k[\fg_\bd]^{G_\bd}\right).$
\end{proof}

We will also use the following Lemma, which follows from \cite{CB1}: 

\begin{lem}\label{bigd}
Assume that $Q$ is neither Dynkin nor extended Dynkin.
Then for any integer $N\ge 1$, the set of vectors $\mathbf d\in \Sigma_0$ 
such that $d_i\ge N$ for all $i$ is Zariski dense in 
$\mathbb \k I$.  
\end{lem}
\begin{proof}
According to Corollary 5.7 of \cite{CB1},
$\Sigma_0$ is the set of all nonzero dimension vectors $\mathbf d$
such that for any nonzero dimension vectors
$\alpha,\beta$ such that $\alpha+\beta=\mathbf d$, one has
$(\alpha,A\beta)\le -2$.
                                                                                
Recall that the fundamental region is the set of
nonzero dimension vectors $\mathbf d$ with connected support
such that $\sum_j a_{ij}d_j\le 0$ for all $i$.

Assume that $\mathbf d$ is in the fundamental region and
$(\mathbf d,A\mathbf d)<0$.
Then for any nonzero dimension vectors
$\alpha,\beta$ such that $\alpha+\beta=\mathbf d$, by Lemma 8.2 of \cite{CB1}
one has
$(\alpha,A\beta)=(\mathbf d-\beta,A\beta)<0$. 
But it is clear  from the definitions that $(\beta,A\beta)$ is  an {\em even}
integer for any integral vector $\beta$. This implies that
if in addition all $d_i$ are even then $(\alpha,A\beta)=(\mathbf
d-\beta,A\beta)
\le -2$.
Thus, if $\mathbf d$ is in the fundamental region,
$(\mathbf d,A\mathbf d)<0$, and $d_i$ are even, then $\mathbf d\in \Sigma_0$.
                                                                                
Now consider the Perron-Frobenius eigenvector $\mathbf v=(v_i)$
of the matrix $C$, with $\sum v_i=1$. Since $Q$ is not Dynkin and not
affine, one has $C\mathbf v=\lambda\mathbf  v$,
where $\lambda>2$. Thus $(\mathbf v,A\mathbf v)=(2-\lambda)|\mathbf v|^2<0$, and
$\sum_j a_{ij}v_j=(2-\lambda)v_i<0$.
This means that there exists an open cone $K$ in $\mathbb R^I$
with axis of symmetry
going through $\mathbf v$ and such that any vector $\mathbf d\in K$ with even
integer coordinates $\mathbf d_i$ is contained in $\Sigma_0$.
This implies the lemma.
\end{proof}                                                                                

\subsection{Proof of Theorem \ref{cen}.}
Given an element $x\in\Sym\PP$, we will denote by
$\bar x$ the image of that element 
under the projection $\Sym\PP\onto\Sym(\PP\br).$

Let  $\k[ w_{im}]$ be
 the polynomial algebra
in an infinite set of  variables $ w_{im},\,i\in I,\,m=1,2,\ldots.$
The assignment $ w_{im}\mto e_i\omh^m\in\PP$
clearly induces a graded algebra isomorphism
$\Sym R\otimes \k[ w_{im}]\iso$
$\Sym(R[\omh]).$
Abusing the notation, we will often identify
$ w_{im}$ with the corresponding element
of $\Sym(R[\omh]).$

Let $\bd$ be a dimension vector and
$\psi_\bd: \oo(\PP)\to\k[\Rep (\overline Q,\mathbf d)]^{G_\bd}$
the corresponding algebra homomorphism.
We define the following functions
$$
\rho_{im}:=
{\psi}_{\mathbf d}(\bar w_{im})=\Tr\wh{e_i\omh^m}\in
\k[\Rep (\overline Q,\mathbf d)]^{G_\bd},\quad\forall i\in I,\,m=1,2,\ldots. 
$$

The relation between $\mu$ and $\omh$ provided by Theorem \ref{moment_thm} 
shows that the subalgebra
$\mu^*\left(\k[\fg_\bd]\right)\sset \k[\Rep (\overline Q,\mathbf d)]$
is generated by the matrix elements of the matrix valued
function $\wh\omh$. Restricting attention to
 $G_\bd$-invariants and using surjectivity of the homomorphism
$\psi_\bd$, we obtain
\beq{muinv}
\mu^*\big(\k[\fg_\bd]^{G_\bd}\big)=
\psi_\bd\big(\Sym(R[\omh])\big)=\k[\rho_{im}].
\eeq

We now  prove part (i) of the Theorem. Let ${\bar c}\in\Sym(\PP\br)$ be
a  homogeneous  element 
 of  degree $r>0$ which is a central element with respect
to the Poisson bracket. 
Then, the element ${\psi}_{\mathbf d}({\bar c})$ has to be central in the  Poisson 
algebra $\k[\Rep (\overline Q,\mathbf d)]^{G_\bd}$,
since the map ${\psi}_{\mathbf d}$ is surjective.
It follows from Corollary \ref{pasha_cor} that,
for any $\mathbf d\in \Sigma_0$, we have
${\psi}_{\mathbf d}({\bar c})\in\mu^*\left(\k[\fg_\bd]^{G_\bd}\right).$
Hence, for such a $\bd$, by \eqref{muinv},
 there exists a degree $r$ homogeneous
element $z\in \k[\bar w_{im}]$ such
that ${\psi}_{\mathbf d}({\bar c})=
\psi_{\mathbf d}(z).
$
Furthermore,
if the dimension vector $\mathbf d\in \Sigma_0$ is such that, in addition, one has
$d_i\ge N(r)$, then, by Proposition \ref{iso1}, in $\oo_\bd(\PP)$,
we get ${\bar c}=z.$ Thus, we have shown that
\beq{equality}
{\bar c}\in \k[\bar w_{im}](r)\sset \oo_{\mathbf d}(\PP)(r),
\quad\forall\mathbf d\in \Sigma_0,\,
d_i\ge N(r), i\in I.
\eeq

Next,  identify  $\Sym R$ with $\k[R^*]$,
the polynomial algebra on the dual vector space
 $R^*=\Hom_\k(R,\k)$. Any  dimension vector may be
viewed as a point $\bd\in R^*.$ Thus, for any $f\in\Sym R=\k[R^*]$  there is a well-defined
scalar $f(\bd)\in\k.$

Now, write
 the direct sum decomposition $\PP\br=R\oplus\BL,$
where $\BL$ is the sum of all homogeneous components of positive
degrees. Thus, we have 
an algebra isomorphism $\Sym R \o\Sym\BL\iso \Sym\PP\br,$
and we can write our element ${\bar c}\in \Sym\PP\br$  
in the form 
${\bar c}=\sum_{j=1}^l {f}_j\cdot  {\bar c}_j$, where ${f}_j$ are elements of $\Sym R$ and 
${\bar c}_j$ are linearly independent elements of $\Sym \BL$.
With this
notation, the image of ${\bar c}$ in $\oo_\bd(\PP)$
equals $\sum_{j=1}^l {f}_j(\bd)\cdot  {\bar c}_j.$

Recall next that
the imbedding $\Sym\BL\into \Sym\PP\br$ induces an isomorphism
$\Sym\BL\iso \oo_\bd(\PP).$ 
We see that formula \eqref{equality} may be intrpreted as saying that
the polynomial map
$R^*\to \Sym \BL,\,\bd\mto \sum_{j=1}^l
{f}_j(\bd)\cdot {\bar c}_j$, takes values in the subspace
$\k[\bar w_{im}](r),$ for every $\bd\in \Sigma_0$ such that
$d_i\ge N(r), \forall i\in I.$
By Lemma \ref{bigd}, this set of dimension vectors is Zariski dense.
Hence, we deduce that $\sum_{j=1}^l
{f}_j(\bd)\cdot {\bar c}_j\in\k[\bar w_{im}](r)$ for all $\bd\in R^*$,
which means that 
 $\sum_{j=1}^l {f}_j\cdot  {\bar c}_j\in 
\Sym R\o\k[\bar w_{im}]=\Sym(\overline{R[\omh]}).$ 
This completes the proof of
the first statement of part (i) of the theorem.

To prove the second  statement of part (i), 
let $c\in \Sym(R[\omh])$ be
a  nonzero homogeneous element of degree $r>0$ such  that
its image, $\bar c\in\Sym(\overline{R[\omh]})$, vanishes.
We can write $c=\sum_j f_j\cdot p_j(w_{im})$ where
$f_j\in\Sym R$ and  where $p_j=p_j(w_{im})\in \k[w_{im}]$ are some 
homogeneous linearly 
independent elements of positive degrees.
Then,  since  $\bar c=0$, applying the homomorphism $\psi_\bd$,
 in  $\k[\Rep(\overline{Q},\bd)]$, we find 
$$0=\psi_\bd(\bar c)=
\sum\nolimits_j f_j(\bd)\cdot p(\rho_{im}).$$

Now, let
 $\mathbf d\in \Sigma_0$ be such that
$d_i\ge N(r), \forall i\in I.$
Then
the only independent relation between the functions
$\rho_{im}\in \Rep(\overline{Q},\bd)$ with $m<N$ is
the relation $\sum_i
\rho_{i1}=\psi_\bd(\omh)=0$.
Thus, $\sum_j f_j(\bd)\cdot p_j$, thought of
as an element of an abstract polynomial
algebra $\k[w_{im}],$ must belong to the principal ideal generated
by the element $\sum_i
w_{i1}$. Now, Lemma \ref{bigd}
and the linear independendence of the polynomials $p_j=p_j(w_{im})$
allows us to conclude that 
the element $\sum_j f_j\o p_j\in\Sym R\o\k[w_{im}]$
is divisible by $\sum_i w_{i1}$. Part (i) of the theorem
follows.


Now we prove part (ii) of the theorem
by a very similar argument. 

Let $\bar c$ be a homogeneous central element 
of $\Sym \Pi\br$ of degree $r>0$.
The map $\psi_\bd: \oo_\bd(\Pi)\to\k[\Rep(\Pi,\bd)]^{G_\bd}$
being surjective, we conclude that
the element $\bar c$ has to be a central element in the Poisson algebra
$\k[M(Q,\mathbf d)]$. 
By Theorem \ref{quivar}, if $\mathbf d\in \Sigma_0$, 
this scheme is in fact an irreducible algebraic variety
which is, moreover, generically symplectic. 
Therefore, any Casimir on $M(Q,\mathbf d)$ has to be a constant function, and hence
$\psi_{\mathbf d}(\bar c)=0$.

Write $\Pi\br=R\oplus L$,
where $L$ is the positive degree part of $\Pi\br$. We can write $\bar c$ in the form 
$\bar c=\sum_{j=1}^l {f}_j\cdot \bar c_j$, where ${f}_j$ are elements of $\Sym R$ and 
$\bar c_j$ are linearly independent elements of $\Sym L$. Again, 
applying the homomorphism $\psi_\bd$ in our $\Pi$-setting, we get
 $0=\psi_\bd(\bar c)=\sum_{j=1}^l {f}_j(\bd)\cdot  \psi_\bd(\bar c_j)$.

 Let $d_i\ge N(r)$ for all $i$.
Then, by Proposition  \ref{iso1}, we deduce that $\sum_j f_j(\mathbf d)\cdot \bar c_j=0$. 
This implies that $f_j(\mathbf d)=0$
for any $j$ and any $\mathbf d\in \Sigma_0$ such that 
$d_i\ge N(r)$. We know from Lemma \ref{bigd}
that the set of such $\mathbf d$ is Zariski dense. 
Thus the polynomials $f_j$ are all identically zero, 
i.e. $\bar c=0$, which proves (ii).

\subsection{The center of $\Pi$.}\label{proof_hh0}
In this section we will prove Proposition \ref{hh0}.
The argument is based on Lemma \ref{bigd}.

Let $z$ be a central element of $\Pi$.
To show that $z$ is a constant it suffices to show that
this is so in every finite dimensional representation of $\Pi$.
Indeed, let ${\scr J}_N$ be the ideal of elements of degree $\ge N$ in $\Pi$; then
any element of $\Pi$ acting by a scalar in $\Pi/{\scr J}^N$ for all $N$ is necessarily
a scalar.
                                                                                
Let $Y$ be a finite dimensional representation of $\Pi$ with dimension vector
$d$. We want to show that $z$ acts by a scalar in $Y.$
By Lemma \ref{bigd}, there exists a nonnegative dimension vector $d'$
such that $d+d'\in \Sigma_0$ (stable region). Let $Y'$ be the direct sum of $Y$
with the augmentation representation of dimension $d'$. Then
$Y'\in \Rep(\bar Q,d+d')$, which is an irreducible variety
by \cite{CB2}.    So it suffices to show that $z$ acts
by a scalar on the generic representation of dimension $d+d'$.
But the generic representation is irreducible, again by \cite{CB2}.
So we are done by Schur's lemma.

We now prove similarly that the algebra $\Pi$ is prime.
Suppose $a,b\in\Pi$ are nonzero, but $a\Pi b=0$.
Choose a dimension vector $\mathbf d$ such that
there are representations of dimension $\mathbf d$
which are not annihilated by $a$ and representations
which are not annihilated by $b$. By enlarging $d$ if
necessary we may assume that $\mathbf d\in\Sigma_0$.
Let $R_a$ and $R_b$ be the closed subsets consisting
of the representations in $\Rep(\Pi,\mathbf d)$
annihilated by $a$ and $b$ respectively. By assumption
they are proper subsets. If $X$ is a simple
representation of $\Pi$, then $a\Pi b X = 0$, but
$\Pi b X = 0$ or $X$, so $aX =0$ or $bX = 0$. Thus
$R_a\cup R_b$ contains all simple representations, and
since they are dense in $\Rep(\Pi,\mathbf d)$ we have
$R_a\cup R_b = \Rep(\Pi,\mathbf d)$. But this contradicts
the fact that $\Rep(\Pi,\mathbf d)$ is irreducible.
\qed

{\small{

}}

\footnotesize{
{\bf W. C-B.}: Department of Pure Mathematics, University of Leeds,
Leeds LS2 9JT, UK;\\
\hphantom{x}\quad\, {\tt W.Crawley-Boevey@leeds.ac.uk}
\smallskip

{\bf P.E.}: Department of Mathematics, Rm 2-165, MIT,
77 Mass. Ave, Cambridge,\\
\hphantom{x}\quad\,  MA 02139, USA;
\quad {\tt etingof@math.mit.edu}
\smallskip

{\bf V.G.}: Department of Mathematics, University of Chicago,
Chicago, IL
60637, USA;\\
\hphantom{x}\quad\, {\tt ginzburg@math.uchicago.edu}}

\end{document}